\documentclass[a4paper, 10pt]{amsart}
\usepackage{amscd, amssymb, color, booktabs}
\usepackage{myown01}
\title[Uniform Projection of the Weil Character]
{Decomposition of the Uniform Projection of the Weil Character}
\author{Shu-Yen Pan}
\address{Department of Mathematics,
National Tsing Hua University, Hsinchu 300, Taiwan}
\email{sypan@math.nthu.edu.tw}

\thanks{This work is partially supported by Taiwan MOST-grant 107-2115-M-007-010-MY2.}

\keywords{Howe correspondence, unipotent character, Lusztig series, reductive dual pair}
\subjclass[2010]{Primary: 20C33}

\date{\today}

\begin{document}

\begin{abstract}
In this paper we obtain a decomposition formula of the uniform projection
of the Weil character of a finite reductive dual pair
consisting of a symplectic group and an even orthogonal group.
This is the first and major step to give an explicit description of the 
Howe correspondence on unipotent characters confirming a conjecture by Aubert-Michel-Rouquier.

\end{abstract}

\maketitle
\tableofcontents


\section{Introduction}

\subsection{}
Let $\omega^\psi_{\Sp_{2N}}$ be the character of the Weil representation (\cf.~\cite{howe-finite}, \cite{gerardin})
of a finite symplectic group $\Sp_{2N}(q)$ with respect to a nontrivial additive character $\psi$ of a finite field
$\bff_q$ of characteristic $p\neq 2$.
Let $(G,G')$ be a reductive dual pair in $\Sp_{2N}(q)$.
It is well known that there are the following three basic types of reductive dual pairs:
\begin{enumerate}
\item[(1)] two general linear groups $(\GL_n,\GL_{n'})$;

\item[(2)] two unitary groups $(\rmU_n,\rmU_{n'})$;

\item[(3)] one symplectic group and one orthogonal group $(\Sp_{2n},\rmO^\epsilon_{n'})$
\end{enumerate}
where $\epsilon=\pm$.
Now $\omega^\psi_{\Sp_{2N}}$ is regarded as a character of $G\times G'$ and denoted by
$\omega^\psi_{\bfG,\bfG'}$ via the homomorphisms $G\times G'\rightarrow G\cdot G'\hookrightarrow \Sp_{2N}(q)$.
Then $\omega^\psi_{\bfG,\bfG'}$ is decomposed as a sum of irreducible characters
\[
\omega^\psi_{\bfG,\bfG'}=\sum_{\rho\in\cale(G),\ \rho'\in\cale(G')}m_{\rho,\rho'}\rho\otimes\rho'
\]
where each $m_{\rho,\rho'}$ is a non-negative integer,
and $\cale(G)$ (resp.~$\cale(G')$) denotes the set of irreducible characters of $G$ (resp.~$G'$).
We say that $\rho\otimes\rho'$ \emph{occurs} in $\omega^\psi_{\bfG,\bfG'}$
if $m_{\rho,\rho'}\neq 0$,
and then we establish a relation $\Theta_{\bfG,\bfG'}$ between
$\cale(G)$ and $\cale(G')$ called the \emph{Howe correspondence} (or \emph{$\Theta$-correspondence})
for the dual pair $(G,G')$.
The main task is to describe the correspondence explicitly.

\subsection{}
Let $G$ be a member of a finite reductive dual pair.
Recall that for a maximal torus $T$ and $\theta\in\cale(T)$,
Deligne and Lusztig define a virtual character $R_{T,\theta}$ of $G$ (\cf.~\cite{dl}, \cite{carter-finite}).
A class function on $G$ is called \emph{uniform} if it is a linear combination of the $R_{T,\theta}$'s.
For a class function $f$,
let $f^\sharp$ denote its projection on the subspace of uniform class functions.

Let $(G,G')$ be a reductive dual pair.
If we assume that the cardinality $q$ of the base field is sufficiently large,
a decomposition of the uniform projection $\omega^\sharp_{\bfG,\bfG'}$ as a linear combination
of Deligne-Lusztig virtual characters is obtained in \cite{srinivasan} p.148
and \cite{pan-odd} theorem 3.8.
Based on Srinivasan's proof, we show that in fact the decomposition holds for any
finite field of odd characteristic (\cf.~Theorem~\ref{0338}).
This is the first step to study Howe correspondence by using Deligne-Lusztig's theory.

From the decomposition of $\omega^\sharp_{\bfG,\bfG'}$ by Srinivasan,
except for the case $(\Sp_{2n},\rmO_{2n'+1})$,
it is known that the unipotent characters are preserved by the $\Theta$-correspondence.
Hence we have the decomposition
\[
\omega_{\bfG,\bfG',1}=\sum_{\rho\in\cale(G)_1,\ \rho'\in\cale(G')_1}m_{\rho,\rho'}\rho\otimes\rho'
\]
where $\omega^\psi_{\bfG,\bfG',1}$ denotes the unipotent part of $\omega^\psi_{\bfG,\bfG'}$,
and $\cale(G)_1$ denotes the subset of irreducible unipotent characters.

In \cite{amr} th\'eor\`eme 5.5, th\'eor\`eme 3.10 and conjecture 3.11,
Aubert, Michel and Rouquier give an explicit description (in terms of partitions or bi-partitions)
of the correspondence of unipotent characters for a dual pair of the first two cases
and provide a conjecture on the third case (when the orthogonal group is even).
To describe the conjecture, we need more notations.

\subsection{}
So now we focus on the correspondence of irreducible unipotent characters
for a dual pair of a symplectic group and an even orthogonal group.
First we review Lusztig's theory on the classification of the irreducible unipotent characters
in \cite{lg}, \cite{lg-symplectic} and \cite{lg-orthogonal}.
Let
\[
\Lambda=\binom{A}{B}=\binom{a_1,a_2,\ldots,a_{m_1}}{b_1,b_2,\ldots,b_{m_2}}
\]
denote a (reduced) \emph{symbol} where $A,B$ are finite subsets of non-negative integers such that $0\not\in A\cap B$.
Note that we always assume that $a_1>a_2>\cdots>a_{m_1}$ and $b_1>b_2>\cdots>b_{m_2}$.
The \emph{rank} and \emph{defect} of a symbol $\Lambda$ are defined by
\begin{align*}
\rank(\Lambda)
&= \sum_{a_i\in A}a_i+\sum_{b_j\in B}b_j-\left\lfloor\biggl(\frac{|A|+|B|-1}{2}\biggr)^2\right\rfloor,\\
{\rm def}(\Lambda) &=|A|-|B|
\end{align*}
where $|X|$ denotes the cardinality of a finite set $X$.
Then Lusztig gives a parametrization of irreducible unipotent characters $\cale(G)_1$
by symbols $\cals_\bfG$:
\[
\begin{cases}
\cale(\Sp_{2n}(q))_1 \\
\cale(\rmO_{2n}^+(q))_1 \\
\cale(\rmO_{2n}^-(q))_1 \\
\end{cases}
\text{parametrized by}
\begin{cases}
\cals_{\Sp_{2n}}=\{\,\Lambda\mid{\rm rank}(\Lambda)=n,\ {\rm def}(\Lambda)\equiv 1\pmod 4\,\};\\
\cals_{\rmO^+_{2n}}=\{\,\Lambda\mid{\rm rank}(\Lambda)=n,\ {\rm def}(\Lambda)\equiv 0\pmod 4\,\};\\
\cals_{\rmO^-_{2n}}=\{\,\Lambda\mid{\rm rank}(\Lambda)=n,\ {\rm def}(\Lambda)\equiv 2\pmod 4\,\}.
\end{cases}
\]
The irreducible character parametrized by $\Lambda$ is denoted by $\rho_\Lambda$.

For a symbol $\Lambda=\binom{a_1,a_2,\ldots,a_{m_1}}{b_1,b_2,\ldots,b_{m_2}}$,
we associate it a \emph{bi-partition}
\[
\Upsilon(\Lambda)
=\sqbinom{a_1-(m_1-1),a_2-(m_1-2),\ldots,a_{m_1-1}-1,a_{m_1}}{b_1-(m_2-1),b_2-(m_2-2),\ldots,b_{m_2-1}-1,b_{m_2}}.
\]
For $(\bfG,\bfG')=(\Sp_{2n},\rmO^\epsilon_{2n'})$,
we define a relation on $\cals_\bfG\times\cals_{\bfG'}$ by
\begin{align*}
\calb_{\Sp_{2n},\rmO^+_{2n'}}
&=\{\,(\Lambda,\Lambda')\in\cals_\bfG\times\cals_{\bfG'}
\mid\mu^\rmT\preccurlyeq\lambda'^\rmT,\ \mu'^\rmT\preccurlyeq\lambda^\rmT,
\ {\rm def}(\Lambda')=-{\rm def}(\Lambda)+1\,\}; \\
\calb_{\Sp_{2n},\rmO^-_{2n'}}
&=\{\,(\Lambda,\Lambda')\in\cals_\bfG\times\cals_{\bfG'}
\mid\lambda'^\rmT\preccurlyeq\mu^\rmT,\ \lambda^\rmT\preccurlyeq\mu'^\rmT,\
{\rm def}(\Lambda')=-{\rm def}(\Lambda)-1\,\}
\end{align*}
where $\sqbinom{\lambda}{\mu}=\Upsilon(\Lambda)$, $\sqbinom{\lambda'}{\mu'}=\Upsilon(\Lambda')$,
$\lambda^\rmT$ denotes the \emph{dual partition} of a partition $\lambda$,
and for two partitions $\lambda=[\lambda_i]$ and $\mu=[\mu_i]$
(with $\{\lambda_i\},\{\mu_i\}$ are written in decreasing order),
we denote $\lambda\preccurlyeq\mu$ if $\mu_i-1\leq\lambda_i\leq\mu_i$ for each $i$.

Then the conjecture by Aubert-Michel-Rouquier in \cite{amr} can be reformulated as follows:

\begin{conj*}
Let $(\bfG,\bfG')=(\Sp_{2n},\rmO_{2n'}^\epsilon)$.
Then
\[
\omega_{\bfG,\bfG',1}
=\sum_{(\Lambda,\Lambda')\in\calb_{\bfG,\bfG'}}\rho_\Lambda\otimes\rho_{\Lambda'}.
\]
\end{conj*}

\subsection{}
In this article, we prove that in the inner product space of class functions on $G\times G'$,
two vectors $\omega_{\bfG,\bfG',1}$ and
$\sum_{(\Lambda,\Lambda')\in\calb_{\bfG,\bfG''}}\rho_\Lambda\otimes\rho_{\Lambda'}$
have the same uniform projection (Theorem~\ref{0310}):

\begin{thm*}
Let $(\bfG,\bfG')=(\Sp_{2n},\rmO_{2n'}^\epsilon)$.
Then
\[
\omega_{\bfG,\bfG',1}^\sharp
=\sum_{(\Lambda,\Lambda')\in\calb_{\bfG,\bfG'}}\rho_\Lambda^\sharp\otimes\rho_{\Lambda'}^\sharp.
\]
\end{thm*}
Based on the theorem, in the sequential article \cite{pan-finite-unipotent},
we will show that two vectors $\omega_{\bfG,\bfG',1}$ and
$\sum_{(\Lambda,\Lambda')\in\calb_{\bfG,\bfG''}}\rho_\Lambda\otimes\rho_{\Lambda'}$
are in fact equal, i.e.,
the above conjecture is confirmed.
Moreover, in another article \cite{pan-Lusztig-correspondence},
we will show that the Howe correspondence and the Lusztig correspondence are commutative,
and so the Howe correspondence of irreducible characters can be reduced via Lusztig correspondence
to the correspondence on irreducible unipotent characters.
Combining all these results together,
we obtain a complete description of the whole Howe correspondence of irreducible characters
for any finite reductive dual pair of a symplectic group and an orthogonal group.

\subsection{}
The contents of the paper are organized as follows.
In Section~\ref{0236}, we recall the definition and basic properties of symbols introduced by
Lusztig.
Then we discuss the relations $\cald_{Z,Z'}$ and $\calb^\epsilon_{Z,Z'}$ which play the
important roles in our main results.
In particular, Lemma~\ref{0210} is used everywhere.
In Section~\ref{0343},
we recall the Lusztig's parametrization of irreducible unipotent characters of
a symplectic group or an even orthogonal group.
And we prove that the decomposition of $\omega^\sharp_{\bfG,\bfG'}$ holds for any finite field of odd
characteristic.
Then we state our main theorem (Theorem~\ref{0310}) in Subsection~\ref{0322}.
In Section~\ref{0403}, we prove the main result for the cases which involve the correspondence
of irreducible unipotent cuspidal characters.
In Section~\ref{0514}, we analyze the basic properties of the relation $\cald_{Z,Z'}$.
In particular, we show that $(Z,Z')$ occurs in the relation if $\cald_{Z,Z'}$ is non-empty.
In Section~\ref{0706}, we prove our first main result for the case that $\epsilon=+$, $Z,Z'$ are regular and
$\cald_{Z,Z'}$ is one-to-one.
In Section~\ref{0802}, we prove Theorem~\ref{0310} for the case $\epsilon=+$ by reducing it to the case proved
in Section~\ref{0706}.
In Section~\ref{0803}, we prove Theorem~\ref{0310} for the case $\epsilon=-$.


\section{Symbols and Bi-partitions}\label{0236}

\subsection{Symbols}\label{0233}
A \emph{symbol} is an ordered pair
\[
\Lambda=\binom{A}{B}=\binom{a_1,a_2,\ldots,a_{m_1}}{b_1,b_2,\ldots,b_{m_2}}
\]
of two finite subsets $A,B$ (possibly empty) of non-negative integers.
We always assume that elements in $A,B$ are written respectively in strictly decreasing order, i.e.,
$a_1>a_2>\cdots>a_{m_1}$ and $b_1>b_2>\cdots>b_{m_2}$.
The pair $(m_1,m_2)$ is called the \emph{size} of $\Lambda$ and is denoted by
${\rm size}(\Lambda)=(m_1,m_2)$.
Each $a_i$ or $b_i$ is called an \emph{entry} of $\Lambda$.
A symbol is called \emph{degenerate} if $A=B$, and it is called \emph{non-degenerate} otherwise.

The \emph{rank} of a symbol $\Lambda=\binom{A}{B}$ is a non-negative integer defined by
\begin{equation}
\rank(\Lambda)
=\sum_{a_i\in A}a_i+\sum_{b_j\in B}b_j-\left\lfloor\biggl(\frac{|A|+|B|-1}{2}\biggr)^2\right\rfloor,
\end{equation}
and the \emph{defect} of $\Lambda$ is defined to be the difference
\begin{equation}
{\rm def}(\Lambda)=|A|-|B|
\end{equation}
where $|X|$ denotes the cardinality of the finite set $X$.
Note that we do not take the absolute value of the difference to define the defect,
unlike the original definition in \cite{lg} p.133.
We think the definition here will be more convenient for our formulation.
It is easy to see that
\[
{\rm rank}(\Lambda)\geq\left\lfloor\biggl(\frac{{\rm def}(\Lambda)}{2}\biggr)^2\right\rfloor.
\]

For a symbol $\Lambda$,
let $\Lambda^*$ (resp.~$\Lambda_*$) denote the first row (resp.~second row) of $\Lambda$, i.e.,
$\Lambda=\binom{\Lambda^*}{\Lambda_*}$.
For two symbols $\Lambda,\Lambda'$, we define the difference
$\Lambda\smallsetminus\Lambda'=\binom{\Lambda^*\smallsetminus\Lambda'^*}{\Lambda_*\smallsetminus\Lambda'_*}$.
If $\Lambda^*\cap\Lambda'^*=\emptyset$ and $\Lambda_*\cap\Lambda'_*=\emptyset$,
then we define the union
$\Lambda\cup\Lambda'=\binom{\Lambda^*\cup\Lambda'^*}{\Lambda_*\cup\Lambda'_*}$.
It is clear that ${\rm def}(\Lambda\cup\Lambda')={\rm def}(\Lambda)+{\rm def}(\Lambda')$.
For a symbol $\Lambda=\binom{A}{B}$,
we define its \emph{transpose} $\Lambda^\rmt=\binom{B}{A}$.
Clearly,
${\rm rank}(\Lambda^\rmt)={\rm rank}(\Lambda)$ and
${\rm def}(\Lambda^\rmt)=-\,{\rm def}(\Lambda)$.

Define an equivalence relation generated by the rule
\begin{equation}\label{0230}
\binom{a_1,a_2,\ldots,a_{m_1}}{b_1,b_2,\ldots,b_{m_2}}\sim
\binom{a_1+1,a_2+1,\ldots,a_{m_1}+1,0}{b_1+1,b_2+1,\ldots,b_{m_1}+1,0}.
\end{equation}
It is easy to see that two equivalent symbols have the same rank and the same defect.
A symbol $\binom{A}{B}$ is called \emph{reduced} if $0\not\in A\cap B$.
It is clear that a symbol is equivalent to a unique reduced symbol.
Let $\cals_{n,\beta}$ denote the set of reduced symbols of rank $n$ and defect $\beta$.

In the remaining part of this article,
a symbol is always assumed to be reduced unless specified otherwise.

\subsection{Special symbols}\label{0232}
A symbol
\begin{equation}\label{0217}
Z=\binom{a_1,a_2,\ldots,a_{m+1}}{b_1,b_2,\ldots,b_m}
\end{equation}
of defect $1$ is called \emph{special} if
\[
a_1\geq b_1\geq a_2\geq b_2\geq\cdots\geq a_m\geq b_m\geq a_{m+1};
\]
similarly a symbol
\begin{equation}\label{0218}
Z=\binom{a_1,a_2,\ldots,a_m}{b_1,b_2,\ldots,b_m}
\end{equation}
of defect $0$ is called \emph{special} if
\[
a_1\geq b_1\geq a_2\geq b_2\geq\cdots\geq b_{m-1}\geq a_m\geq b_m.
\]

For a special symbol $Z=\binom{A}{B}$,
define
\[
Z_\rmI=\binom{A\smallsetminus(A\cap B)}{B\smallsetminus(A\cap B)}.
\]
Entries in $Z_\rmI$ are called \emph{singles} of $Z$.
Clearly, $Z$ and $Z_\rmI$ have the same defect,
and the \emph{degree} of $Z$ is defined to be
\[
\deg(Z)=\begin{cases}
\frac{|Z_\rmI|-1}{2}, & \text{if $Z$ has defect $1$;}\\
\frac{|Z_\rmI|}{2}, & \text{if $Z$ has defect $0$}.
\end{cases}
\]
A special symbol $Z$ is called \emph{regular} if $Z=Z_\rmI$.
Entries in $Z_{\rm II}:=A\cap B$ are called ``doubles'' in $Z$.
If $a\in A$ and $b\in B$ such that $a=b$,
then the pair $\binom{a}{b}$ is called a \emph{pair of doubles}.

A symbol $M$ is a \emph{subsymbol} of $Z_\rmI$ and denoted by $M\subset Z_\rmI$ if $M^*\subset(Z_\rmI)^*$ and $M_*\subset(Z_\rmI)_*$.
Because we can naturally construct a subsymbol $\binom{M\cap(Z_\rmI)^*}{M\cap(Z_\rmI)_*}$ of $Z_\rmI$ from each subset $M\subset(Z_\rmI)^*\cup(Z_\rmI)_*$.
So we will make no distinguish between a subsymbol of $Z_\rmI$ and a subset of $(Z_\rmI)^*\cup(Z_\rmI)_*$.
Now for a subset $M$ of $Z_\rmI$,
define
\begin{equation}\label{0201}
\Lambda_M=(Z\smallsetminus M)\cup M^\rmt,
\end{equation}
i.e.,
$\Lambda_M$ is obtained from $Z$ by switching the row position of elements in $M$
and keeping other elements unchanged.
It is clear that 
\begin{equation}\label{0221}
(\Lambda_M)^\rmt=\Lambda_{Z_\rmI\smallsetminus M}.
\end{equation}
In particular, $Z=\Lambda_\emptyset$, $Z^\rmt=\Lambda_{Z_\rmI}$.
An entry $z\in\Lambda_M$ is said to be in its \emph{natural position} if it is in the same position in $Z$.
By definition, an entry $z\in\Lambda_M$ is not in its natural position if $z\in M$.

For a special symbol $Z$ and an integer $\delta$,
we define
\begin{align}
\begin{split}
\overline\cals_Z &=\{\,\Lambda_M\mid M\subset Z_\rmI\,\},\\
\cals_{Z,\delta} &=\{\,\Lambda\in\overline\cals_Z\mid {\rm def}(\Lambda)=\delta\,\}.
\end{split}
\end{align}
Note that if $\Lambda\in\overline\cals_Z$,
then ${\rm def}(\Lambda)\equiv{\rm def}(Z)\pmod 2$.
For $\Lambda_M,\Lambda_{M'}\in\overline\cals_Z$,
we define
\[
\Lambda_M+\Lambda_{M'}=\Lambda_N
\]
where $N=(M\cup M')\smallsetminus(M\cap M')$.
This gives $\overline\cals_Z$ an $\bff_2$-vector space structure with identity element $\Lambda_\emptyset=Z$.

\begin{lemma}\label{0231}
Let $M,N$ be subsets of $Z_\rmI$.
Then $(\Lambda_M+\Lambda_N)^\rmt=(\Lambda_M)^\rmt+\Lambda_N$.
\end{lemma}
\begin{proof}
We have a disjoint decomposition
\[
Z_\rmI=(Z_\rmI\smallsetminus(M\cup N))\cup(M\smallsetminus N)\cup(N\smallsetminus M)\cup(M\cap N).
\]
Then by (\ref{0221}), we have
\[
(\Lambda_M+\Lambda_N)^\rmt
=(\Lambda_{(M\smallsetminus N)\cup(N\smallsetminus M)})^\rmt
=\Lambda_{(Z_\rmI\smallsetminus(M\cup N))\cup(M\cap N)}
\]
and
\begin{align*}
(\Lambda_M)^\rmt+\Lambda_N
=\Lambda_{Z_\rmI\smallsetminus M}+\Lambda_N
& =\Lambda_{(Z_\rmI\smallsetminus(M\cup N))\cup(N\smallsetminus M)}+\Lambda_{(M\cap N)\cup(N\smallsetminus M)} \\
& =\Lambda_{(Z_\rmI\smallsetminus(M\cup N))\cup(M\cap N)}.
\end{align*} 
\end{proof}

We define a non-singular symplectic form on $\overline\cals_Z$ (over $\bff_2$) by
\begin{equation}\label{0222}
\langle\Lambda_M,\Lambda_{M'}\rangle
=|M\cap M'|\pmod 2.
\end{equation}

\begin{remark}
Our notation ``$\Lambda_M$'' is different from the original one defined in \cite{lg-symplectic} \S 5.7.
Now we check that both are equivalent.
For $M\subset Z_\rmI$,
let $\Lambda_M$ be defined in (\ref{0201}),
and let $\bar\Lambda_M$ denote the ``$\Lambda_M$'' defined in \cite{lg-symplectic} \S 5.7, i.e.,
\[
\bar\Lambda_M=\binom{Z_{\rm II}\sqcup (Z_\rmI\smallsetminus M)}{Z_{\rm II}\sqcup M}
\]
where $Z_\rmI$ is regarded as $(Z_\rmI)_*\cup(Z_\rmI)^*$, and then
$\langle\bar\Lambda_{M_1},\bar\Lambda_{M_2}\rangle=|M_1^\sharp\cap M_2^\sharp|\pmod2$
where $M^\sharp=M\cup(Z_\rmI)_*\smallsetminus(M\cap(Z_\rmI)_*)$.
Then
\[
\Lambda_M
=\binom{Z_{\rm II}\sqcup((Z_\rmI)^*\smallsetminus M^*)\cup M_*}{Z_{\rm II}\sqcup((Z_\rmI)_*\smallsetminus M_*)\cup M^*}
=\bar\Lambda_{(Z_\rmI\smallsetminus M)_*\cup M^*}.
\]
Let $M'=(Z_\rmI\smallsetminus M)_*\cup M^*$.
Then
\[
(M')^\sharp=(M'\cup(Z_\rmI)_*)\smallsetminus (M'\cap(Z_\rmI)_*)
=(M^*\cup(Z_\rmI)_*)\smallsetminus(Z_\rmI\smallsetminus M)_*
=M^*\cup M_*
=M.
\]
Therefore, the definition of $\Lambda_M$ in (\ref{0201}) with the form $\langle,\rangle$ defined in (\ref{0222})
is equivalent to the original definition in \cite{lg-symplectic} \S 5.7.
\end{remark}

Let $Z$ be a special symbol.
If ${\rm def}(Z)=1$, we define
\begin{equation}\label{0224}
\cals_Z = \bigcup_{\beta\equiv 1\pmod 4}\cals_{Z,\beta};
\end{equation}
if ${\rm def}(Z)=0$, we define
\begin{align}\label{0234}
\begin{split}
\cals^+_Z &= \bigcup_{\beta\equiv 0\pmod 4}\cals_{Z,\beta}; \\
\cals^-_Z &= \bigcup_{\beta\equiv 2\pmod 4}\cals_{Z,\beta}.
\end{split}
\end{align}
It is clear that
\[
\begin{cases}
|\cals_Z|=2^{2\delta}, & \text{if ${\rm def}(Z)=1$ and $\deg(Z)=\delta$};\\
|\cals^+_Z|=|\cals^-_Z|=2^{2\delta-1}, & \text{if ${\rm def}(Z)=0$ and $\deg(Z)=\delta$}.
\end{cases}
\]
For simplicity, if ${\rm def}(Z)$ is not specified,
let $\cals_Z$ denote one of $\cals_Z$ (when ${\rm def}(Z)=1$) or $\cals^+_Z,\cals^-_Z$ (when ${\rm def}(Z)=0$).

\subsection{Bi-partitions}\label{0220}
For a partition
\[
\lambda=[\lambda_i]=[\lambda_1,\lambda_2,\ldots,\lambda_k],\quad\text{ with } \lambda_1\geq\lambda_2\geq\cdots\geq\lambda_k\geq 0,
\]
define $|\lambda|=\lambda_1+\lambda_2+\cdots+\lambda_k$.
For a partition $\lambda=[\lambda_i]$, define its \emph{transpose} (or \emph{dual})
$\lambda^\rmT=[\lambda_j^*]$ by $\lambda_j^*=|\{\,i\mid\lambda_i\geq j\,\}|$ for $j\in\bbN$.
It is clear that $|\lambda^\rmT|=|\lambda|$.

\begin{lemma}\label{0212}
Let $\lambda=[\lambda_i]$ be a partition and let $\lambda^\rmT=[\lambda_j^*]$ be its transpose.
Then for $i,j\in\bbN$,
\begin{enumerate}
\item[(i)] $j\leq\lambda_i$ if and only if $\lambda^*_j\geq i$;

\item[(ii)] $j>\lambda_i$ if and only if $\lambda^*_j<i$.
\end{enumerate}
\end{lemma}
\begin{proof}
Part (i) follows from the definition directly,
and part (ii) is equivalent to (i).
\end{proof}

Let $\lambda=[\lambda_1,\ldots,\lambda_k]$ and $\mu=[\mu_1,\ldots,\mu_l]$ be two partitions with
$\lambda_1\geq\cdots\geq\lambda_k$ and $\mu_1\geq\cdots\geq\mu_l$.
We may assume that $k=l$ by adding several $0$'s if necessary.
Then we denote
\[
\lambda\preccurlyeq\mu\quad\text{ if }\mu_i-1\leq\lambda_i\leq\mu_i\text{ for each }i.
\]

\begin{lemma}\label{0202}
Let $\lambda=[\lambda_i]$, $\mu=[\mu_i]$ be two partitions.
Then $\lambda^\rmT\preccurlyeq\mu^\rmT$ if and only if $\mu_{i+1}\leq\lambda_i\leq\mu_i$ for each $i$.
\end{lemma}
\begin{proof}
Write $\lambda^\rmT=[\lambda^*_j]$ and $\mu^\rmT=[\mu^*_j]$.
First suppose that $\mu^*_j\geq\lambda^*_j$ for each $j$.
If $\lambda_i\geq j$, then by Lemma~\ref{0212}, $\mu^*_j\geq\lambda_j^*\geq i$ and hence $\mu_i\geq j$.
Since $j$ is arbitrary, this implies that $\mu_i\geq\lambda_i$ for each $i$.
By the same argument, $\mu_i\geq\lambda_i$ for each $i$ implies that $\mu_j^*\geq\lambda^*_j$ for each $j$.

Next suppose that $\mu_j^*-1\leq\lambda^*_j$ for each $j$.
If $\mu_{i+1}\geq j$,
then $\mu^*_j\geq i+1$ by Lemma~\ref{0212}.
Then $\lambda^*_j\geq i$ and hence $\lambda_i\geq j$.
Since $j$ is arbitrary,
$\mu_j^*-1\leq\lambda^*_j$ for each $j$ implies $\lambda_i\geq\mu_{i+1}$ for each $i$.
Now suppose that $\lambda_i\geq\mu_{i+1}$ for each $i$.
If $\mu_j^*\geq i+1$, then $\mu_{i+1}\geq j$ and $\lambda_i\geq j$, and hence $\lambda^*_j\geq i$.
Since $i$ is arbitrary, $\lambda_i\geq\mu_{i+1}$ for each $i$ implies that $\lambda^*_j\geq\mu^*_j-1$ for each $j$.
\end{proof}

Let $\calp_2(n)$ denote the set of bi-partitions of $n$,
i.e., the set of $\sqbinom{\lambda}{\mu}$ where $\lambda,\mu$ are partitions and $|\lambda|+|\mu|=n$.
To each symbol we can associate a bi-partition by:
\begin{equation}\label{0219}
\Upsilon\colon\binom{a_1,a_2,\ldots,a_{m_1}}{b_1,b_2,\ldots,b_{m_2}}\mapsto
\sqbinom{a_1-(m_1-1),a_2-(m_1-2),\ldots,a_{m_1-1}-1,a_{m_1}}{b_1-(m_2-1),b_2-(m_2-2),\ldots,b_{m_2-1}-1,b_{m_2}}.
\end{equation}
It is easy to check that $\Upsilon$ induces a bijection
\[
\Upsilon\colon\cals_{n,\beta}\longrightarrow
\begin{cases}
\calp_2(n-(\frac{\beta-1}{2})(\frac{\beta+1}{2})), & \text{if $\beta$ is odd};\\
\calp_2(n-(\frac{\beta}{2})^2), & \text{if $\beta$ is even}.
\end{cases}
\]
In particular, $\Upsilon$ induces a bijection from $\cals_{n,1}$ onto $\calp_2(n)$, and
a bijection from $\cals_{n,0}$ onto $\calp_2(n)$.
Therefore, for each $\sqbinom{\lambda}{\mu}\in\calp_2(n-(\frac{\beta}{2})^2)$ with $\beta>0$,
there exist a unique $\Lambda\in\cals_{n,\beta}$ and a unique $\Lambda'\in\cals_{n,-\beta}$ such that
$\Upsilon(\Lambda)=\Upsilon(\Lambda')=\sqbinom{\lambda}{\mu}$.

\subsection{The sets $\calb^\epsilon_{Z,Z'}$ and $\cald_{Z,Z'}$}
Let $Z,Z'$ be special symbols of defect $1,0$ respectively.
For $\Lambda\in\overline\cals_Z$ write $\Upsilon(\Lambda)=\sqbinom{\lambda}{\mu}$;
and for $\Lambda'\in\overline\cals_{Z'}$ write $\Upsilon(\Lambda')=\sqbinom{\xi}{\nu}$.
Define two relations $\overline\calb^+_{Z,Z'},\overline\calb^-_{Z,Z'}$ on $\overline\cals_Z\times\overline\cals_{Z'}$ by
\begin{align}\label{0223}
\begin{split}
\overline\calb^+_{Z,Z'} &=\{\,(\Lambda,\Lambda')\in\overline\cals_Z\times\overline\cals_{Z'}\mid\mu^\rmT\preccurlyeq\xi^\rmT,\ \nu^\rmT\preccurlyeq\lambda^\rmT,
\ {\rm def}(\Lambda')=-{}{\rm def}(\Lambda)+1\,\};\\
\overline\calb^-_{Z,Z'} &=\{\,(\Lambda,\Lambda')\in\overline\cals_Z\times\overline\cals_{Z'}\mid\xi^\rmT\preccurlyeq\mu^\rmT,\ \lambda^\rmT\preccurlyeq\nu^\rmT,
\ {\rm def}(\Lambda')=-{\rm def}(\Lambda)-1\,\}.
\end{split}
\end{align}
We say that symbol $\Lambda\in\overline\cals_Z$ \emph{occurs} in the relation $\overline\calb^\epsilon_{Z,Z'}$ (where $\epsilon=\pm$)
if there exists $\Lambda'\in\overline\cals_{Z'}$ such that $(\Lambda,\Lambda')\in\overline\calb^\epsilon_{Z,Z'}$.
An analogous definition also applies to a symbol $\Lambda'\in\overline\cals_{Z'}$.
Finally we define
\begin{align}\label{0227}
\begin{split}
\calb^\epsilon_{Z,Z'} &=\overline\calb^\epsilon_{Z,Z'}\cap(\cals_Z\times\cals^\epsilon_{Z'});\\
\cald_{Z,Z'} &= \overline\calb^+_{Z,Z'}\cap(\cals_{Z,1}\times\cals_{Z',0});\\
\calb_{\Sp_{2n},\rmO^\epsilon_{2n'}}=\calb^\epsilon_{n,n'} &=\bigsqcup_{Z,Z'}\calb^\epsilon_{Z,Z'};\\
\calb^\epsilon &=\bigsqcup_{n,n'\geq 0}\calb^\epsilon_{n,n'}
\end{split}
\end{align}
where the disjoint union $\bigsqcup_{Z,Z'}$ is taken over all special symbols $Z,Z'$ 
of rank $n,n'$ and defect $1,0$ respectively.
Because $\cals_{Z,1}\subset\cals_Z$ and $\cals_{Z',0}\subset\cals^+_{Z'}$,
it is clear that $\cald_{Z,Z'}\subset\calb^+_{Z,Z'}$.
Note that if $\Lambda\in\cals_Z\subset\overline\cals_Z$ and $(\Lambda,\Lambda')\in\overline\calb^\epsilon_{Z,Z'}$
for some $\Lambda'\in\overline\cals_{Z'}$,
then ${\rm def}(\Lambda)\equiv 1\pmod 4$ and by (\ref{0223}) we know that $\Lambda'\in\cals^\epsilon_{Z'}$,
and consequently $(\Lambda,\Lambda')\in\calb^\epsilon_{Z,Z'}$.
Conversely, if $\Lambda'\in\cals_{Z'}^\epsilon\subset\overline\cals_{Z'}$
and $(\Lambda,\Lambda')\in\overline\calb^\epsilon_{Z,Z'}$
for some $\Lambda\in\overline\cals_Z$,
then $\Lambda\in\cals_Z$ and hence $(\Lambda,\Lambda')\in\calb^\epsilon_{Z,Z'}$.
Define
\begin{align}\label{0235}
\begin{split}
B^\epsilon_\Lambda &=\{\,\Lambda'\in\overline\cals_{Z'}\mid(\Lambda,\Lambda')\in\overline\calb^\epsilon_{Z,Z'}\,\}\quad\text{ for }\Lambda\in\overline\cals_Z,\\
B^\epsilon_{\Lambda'} &=\{\,\Lambda\in\overline\cals_Z\mid(\Lambda,\Lambda')\in\overline\calb^\epsilon_{Z,Z'}\,\}\quad\text{ for }\Lambda'\in\overline\cals_{Z'}.
\end{split}
\end{align}
If $\Sigma\in\cals_{Z,1}$,
then $B^+_\Sigma$ is also denoted by $D_\Sigma$.
Similarly, $B^+_{\Sigma'}$ is denoted by $D_{\Sigma'}$ if $\Sigma'\in\cals_{Z',0}$.

\begin{lemma}\label{0213}
Let $Z,Z'$ be special symbols of size $(m+1,m), (m',m')$ respectively for some non-negative integers $m,m'$.
If\/ $\cald_{Z,Z'}\neq\emptyset$,
then either $m'=m$ or $m'=m+1$.
\end{lemma}
\begin{proof}
Write
\begin{align*}
\Sigma &=\binom{a_1,a_2,\ldots,a_{m+1}}{b_1,b_2,\ldots,b_m}, &
\Sigma' &=\binom{a'_1,a'_2,\ldots,a'_{m'}}{b'_1,b'_2,\ldots,b'_{m'}}; \\
\Upsilon(\Sigma) &=\sqbinom{\lambda_1,\lambda_2,\ldots,\lambda_{m+1}}{\mu_1,\mu_2,\ldots,\mu_m}, &
\Upsilon(\Sigma') &=\sqbinom{\lambda'_1,\lambda'_2,\ldots,\lambda'_{m'}}{\mu'_1,\mu'_2,\ldots,\mu'_{m'}}
\end{align*}
where $\Upsilon$ is given in (\ref{0219}).
Suppose that $(\Sigma,\Sigma')\in\cald_{Z,Z'}$.
If $a_{m+1}\neq 0$, then $\lambda_{m+1}\neq 0$.
This implies that $\mu'_{m}\geq\lambda_{m+1}>0$ by Lemma~\ref{0202}, and hence $m'\geq m$.
If $b_m\neq 0$, then $\mu_m\neq 0$.
This implies that $\lambda_m\geq\mu_m>0$ and hence again $m'\geq m$.
Because we assume that $\Sigma$ is reduced, $a_{m+1},b_m$ can not both be zero,
so we conclude that $m'\geq m$.

Similarly, if $a'_{m'}\neq 0$, then $\lambda'_{m'}\neq 0$.
Then $\mu_{m'-1}\geq\lambda'_{m'}>0$ and hence $m\geq m'-1$.
If $b'_{m'}\neq 0$, then $\mu'_{m'}\neq 0$.
Then $\lambda_{m'}\geq\mu_{m'}>0$ and hence $m+1\geq m'$.
Because $a'_{m'},b'_{m'}$ can not both be zero,
so we conclude that $m\geq m'-1$.
\qed
\end{proof}

\begin{lemma}\label{0210}
Let $Z,Z'$ be two special symbols of sizes $(m+1,m),(m',m')$ respectively
where $m'=m,m+1$.
Let
\[
\Lambda=\binom{a_1,a_2,\ldots,a_{m_1}}{b_1,b_2,\ldots,b_{m_2}}\in\overline\cals_Z,\qquad
\Lambda'=\binom{c_1,c_2,\ldots,c_{m'_1}}{d_1,d_2,\ldots,d_{m'_2}}\in\overline\cals_{Z'}.
\]
\begin{enumerate}
\item[(i)] Then $(\Lambda,\Lambda')\in\overline\calb^+_{Z,Z'}$ if and only if
\[
\begin{cases}
m_1'=m_2,\ \text{ and }a_i>d_i,\ d_i\geq a_{i+1},\ c_i\geq b_i,\ b_i>c_{i+1}\text{ \ for each $i$}, & \text{if $m'=m$};\\
m_1'=m_2+1,\ \text{ and }a_i\geq d_i,\ d_i>a_{i+1},\ c_i>b_i,\ b_i\geq c_{i+1}\text{ \ for each $i$}, & \text{if $m'=m+1$}.
\end{cases}
\]

\item[(ii)] Then $(\Lambda,\Lambda')\in\overline\calb^-_{Z,Z'}$ if and only if
\[
\begin{cases}
m_1'=m_2-1,\ \text{ and }d_i\geq a_i,\ a_i>d_{i+1},\ b_i>c_i,\ c_i\geq b_{i+1}\text{ \ for each $i$}, & \text{if $m'=m$};\\
m_1'=m_2,\ \text{ and }d_i>a_i,\ a_i\geq d_{i+1},\ b_i\geq c_i,\ c_i>b_{i+1}\text{ \ for each $i$}, & \text{if $m'=m+1$}.
\end{cases}
\]
\end{enumerate}
\end{lemma}
\begin{proof}
We know that $m_1+m_2=2m+1$, $m'_1+m'_2=2m'$, ${\rm def}(\Lambda)=m_1-m_2$ and ${\rm def}(\Lambda')=m_1'-m'_2$.
Write $\Upsilon(\Lambda)=\sqbinom{\lambda}{\mu}$,
$\lambda=[\lambda_i]$,
$\lambda^\rmT=[\lambda^*_j]$,
$\mu=[\mu_i]$,
$\mu^\rmT=[\mu^*_j]$;
and $\Upsilon(\Lambda')=\sqbinom{\xi}{\nu}$,
$\xi=[\xi_i]$,
$\xi^\rmT=[\xi^*_j]$,
$\nu=[\nu_i]$,
$\nu^\rmT=[\nu^*_j]$.
Then $\lambda_i=a_i-(m_1-i)$, $\mu_i=b_i-(m_2-i)$,
$\xi_i=c_i-(m_1'-i)$, $\nu_i=d_i-(m_2'-i)$.

\begin{enumerate}
\item[(a)] Suppose that $\epsilon=+$ and $m'=m$.
Suppose that $(\Lambda,\Lambda')\in\overline\calb^+_{Z,Z'}$.
Then $m'_1-m'_2=-(m_1-m_2)+1$ by (\ref{0223}).
Therefore, we have $m'_1=m_2$ and $m'_2=m_1-1$.
Moreover, by Lemma~\ref{0202}, the condition
\begin{equation}
\mu^\rmT\preccurlyeq\xi^\rmT\quad\text{and}\quad\nu^\rmT\preccurlyeq\lambda^\rmT.
\end{equation}
is equivalent to
\begin{equation}\label{0203}
\xi_{i+1}\leq\mu_i\leq\xi_i\quad\text{and}\quad
\lambda_{i+1}\leq\nu_i\leq\lambda_i
\end{equation}
for each $i$.
Moreover, the condition $m'_1=m_2$ and (\ref{0203}) imply that
\begin{equation}\label{0228}
c_{i+1}<b_i\leq c_i\quad\text{and}\quad a_{i+1}\leq d_i<a_i
\end{equation}
for each $i$.

Conversely, the condition $m'_1=m_2$ will imply that
${\rm def}(\Lambda')=-{\rm def}(\Lambda)+1$.
Moreover, the condition $m'_1=m_2$ and (\ref{0228}) will imply (\ref{0203}),
and hence $(\Lambda,\Lambda')\in\overline\calb^+_{Z,Z'}$.

\item[(b)] Suppose that $\epsilon=+$ and $m'=m+1$.
Then the condition ${\rm def}(\Lambda')=-{\rm def}(\Lambda)+1$ is equivalent to $m'_1=m_2+1$.
Moreover, under the condition $m'_1=m_2+1$, (\ref{0203}) is equivalent to
\[
c_{i+1}\leq b_i<c_i\quad\text{and}\quad a_{i+1}<d_i\leq a_i
\]
for each $i$.

\item[(c)] Suppose that $\epsilon=-$ and $m'=m$.
Then the condition ${\rm def}(\Lambda')=-{\rm def}(\Lambda)-1$ is equivalent to $m'_1=m_2-1$.
Moreover, under the condition $m'_1=m_2-1$, the condition
\begin{equation}\label{0229}
\xi^\rmT\preccurlyeq\mu^\rmT\quad\text{and}\quad \lambda^\rmT\preccurlyeq\nu^\rmT
\end{equation}
is equivalent to
\[
b_{i+1}\leq c_i<b_i\quad\text{and}\quad d_{i+1}<a_i\leq d_i
\]
for each $i$.

\item[(d)] Suppose that $\epsilon=-$ and $m'=m+1$.
Then the condition ${\rm def}(\Lambda')=-{\rm def}(\Lambda)-1$ is equivalent to $m'_1=m_2$.
And under the condition $m'_1=m_2$,
the condition (\ref{0229}) is equivalent to
\[
b_{i+1}<c_i\leq b_i\quad\text{and}\quad d_{i+1}\leq a_i<d_i
\]
for each $i$.
\end{enumerate}
Thus the lemma is proved.
\qed
\end{proof}

We will see that the following two examples are basic:

\begin{example}\label{0225}
Let $m$ be a non-negative integer,
and let
\begin{align*}
Z=Z_{(m)} &=\binom{2m,2m-2,\ldots,0}{2m-1,2m-3,\ldots,1}, \\
Z'=Z'_{(m+1)} &=\binom{2m+1,2m-1,\ldots,1}{2m,2m-2,\ldots,0}.
\end{align*}
Then $Z_{(m)}$ is a special symbol of size $(m+1,m)$ and rank $2m(m+1)$;
$Z'_{(m+1)}$ is a special symbol of size $(m+1,m+1)$ and rank $2(m+1)^2$.
Let $\Lambda\in\cals_{Z}$ and $\Lambda'\in\cals^\epsilon_{Z'}$ such that $(\Lambda,\Lambda')\in\calb^\epsilon_{Z,Z'}$.
\begin{enumerate}
\item[(1)] Suppose that $\epsilon=+$.
Then we can write
\[
\Lambda=\binom{a_1,a_2,\ldots,a_{m_1}}{b_1,b_2,\ldots,b_{m_2}},\qquad
\Lambda'=\binom{c_1,c_2,\ldots,c_{m_2+1}}{d_1,d_2,\ldots,d_{m_1}}
\]
for some nonnegative integers $m_1,m_2$ such that $m_1+m_2=2m+1$.
If $0$ is in the first row of $\Lambda$, i.e., $a_{m_1}=0$,
then $0$ must be in the second row of $\Lambda'$, i.e., $d_{m_1}=0$ since we need $a_{m_1}\geq d_{m_1}$ by (1) of Lemma~\ref{0210};
similarly, if $0$ is in the second row of $\Lambda$, i.e., $b_{m_2}=0$,
then $0$ must be in the first row of $\Lambda'$, i.e., $c_{m_2+1}=0$ since we need $b_{m_2}\geq c_{m_2+1}$.
Next, it is easy to see that $i+1$ is in the same row of $i$ in $\Lambda$ if and only if
$i+1$ is in the same row of $i$ in $\Lambda'$ by the rule in Lemma~\ref{0210}.
Therefore the row positions in $\Lambda,\Lambda'$ of the same entry $i$ are always opposite.
Finally, the number of entries in the first row of $\Lambda'$ is one more than the number
of entries in the second row of $\Lambda$,
we see that $2m+1$ is always in the first row of $\Lambda'$.
Thus $\Lambda'=\binom{2m+1}{-}\cup\Lambda^\rmt$.

\item[(2)] Suppose that $\epsilon=-$.
Then we can write
\[
\Lambda=\binom{a_1,a_2,\ldots,a_{m_1}}{b_1,b_2,\ldots,b_{m_2}},\qquad
\Lambda'=\binom{c_1,c_2,\ldots,c_{m_2}}{d_1,d_2,\ldots,d_{m_1+1}}.
\]
If $0$ is in the first row of $\Lambda$, i.e., $a_{m_1}=0$,
then $0$ must be in the second row of $\Lambda'$, i.e., $d_{m_1+1}=0$ since we need $a_{m_1}\geq d_{m_1+1}$;
similarly, if $0$ is in the second row of $\Lambda$, i.e., $b_{m_2}=0$,
then $0$ must be in first row of $\Lambda'$, i.e., $c_{m_2}=0$ since we need $b_{m_2}\geq c_{m_2}$
from Lemma~\ref{0222}.
By the same argument as in (1) we see that the row positions of the entry $i$ in $\Lambda$ and $\Lambda'$ must be opposite for each $i=0,\ldots,2m$,
and the entry $2m+1$ must be in the second row of $\Lambda'$.
Thus $\Lambda'=\binom{-}{2m+1}\cup\Lambda^\rmt$.
\end{enumerate}
\end{example}

\begin{example}\label{0226}
Let $m$ be a non-negative integer,
$Z=Z_{(m)}$ and $Z'=Z'_{(m)}$ where $Z_{(m)},Z'_{(m)}$ are given as in the previous example.
Then $Z$ is a special symbol of size $(m+1,m)$ and rank $2m(m+1)$;
$Z'$ is a special symbol of size $(m,m)$ and rank $2m^2$.
Let $\Lambda\in\cals_Z$ and $\Lambda'\in\cals^\epsilon_{Z'}$.
By the similar argument in the previous example, we see that
$(\Lambda,\Lambda')\in\calb^\epsilon_{Z,Z'}$ if and only if
\[
\Lambda=\begin{cases}
\binom{2m}{-}\cup\Lambda'^\rmt, & \text{if $\epsilon=+$};\\
\binom{-}{2m}\cup\Lambda'^\rmt, & \text{if $\epsilon=-$}.
\end{cases}
\]
\end{example}


\section{Finite Howe Correspondence of Unipotent Characters}\label{0343}
In the first part of this section we review the parametrization of
irreducible unipotent characters by Lusztig in \cite{lg-symplectic} and \cite{lg-orthogonal}.
The comparison of our main results and the conjecture in \cite{amr} is in the final subsection.

\subsection{Weyl groups}\label{0324}
Let $\bff_q$ be a finite field of $q$ element with odd characteristic.
Let $\bfG$ be a classical group defined over $\bff_q$,
$F$ the corresponding Frobenius automorphism,
$G=\bfG^F$ the group of rational points.
If $\bfG$ is $\Sp_{2n}$, $\rmO_n^\epsilon$, or $\rmU_n$,
then $G$ is denoted by $\Sp_{2n}(q)$, $\rmO_n^\epsilon(q)$, or $\rmU_n(q)$ respectively.

Let $\bfT_0$ be a fixed maximally split rational maximal torus of $\bfG$,
$\bfN(\bfT_0)$ the normalizer of $\bfT_0$ in $\bfG$,
$\bfW=\bfW_\bfG=\bfN(\bfT_0)/\bfT_0$ the \emph{Weyl group} of $\bfG$.
For $w\in\bfW$, choose an element $g\in\bfG$ such that $g^{-1}F(g)\in\bfN(\bfT_0)$ whose
image in $\bfW$ is $w$,
and define $\bfT_w=g\bfT_0g^{-1}$,
a rational maximal torus of $\bfG$.
Every rational maximal torus of $\bfG$ is $G$-conjugate to $\bfT_w$ for some $w\in\bfW$.

The Weyl groups of $\Sp_{2n}$ and $\SO_{2n+1}$ can be identified with $W_n$ which consists of permutations on the set
\[
\{1,2,\ldots,n,n^*,(n-1)^*,\ldots,1^*\}
\]
which commutes with the involution $(1,1^*)(2,2^*)\cdots(n,n^*)$
where $(i,j)$ denote the transposition of $i,j$.
Let $s_i=(i,i+1)(i^*,(i+1)^*)$ for $i=1,\ldots,n-1$,
and let $t_n=(n,n^*)$.
It is known that $W_n$ is generated by $\{s_1,\ldots,s_{n-1},t_n\}$.
The kernel $W_n^+$ of the homomorphism $\varepsilon_0\colon W_n\rightarrow\{\pm1\}$
given by $s_i\mapsto 1$ and $t_n\mapsto -1$ is subgroup of index two and is generated
by $\{s_1,\ldots,s_{n-1},t_ns_{n-1}t_n\}$.
Let $W_n^-=W_n\smallsetminus W_n^+$.

\subsection{Deligne-Lusztig virtual characters}\label{0339}
The set of irreducible characters of $G$ is denoted by $\cale(G)$
where an irreducible character means the character of an irreducible representation of $G$.
The space $\calv(G)$ of (complex-valued) class functions on $G$ is an inner product space
and $\cale(G)$ forms an orthonormal basis of $\calv(G)$.

If $\bfG$ is a connected classical group,
let $R_{\bfT,\theta}=R_{\bfT,\theta}^\bfG$ denote the \emph{Deligne-Lusztig virtual character} of $G$
with respect to a rational maximal torus $\bfT$ and an irreducible character $\theta\in\cale(T)$ where $T=\bfT^F$.
If $\bfG=\rmO_n^\epsilon$, then we define
\[
R_{\bfT,\theta}^{\rmO_n^\epsilon}=\Ind_{\SO_n^\epsilon(q)}^{\rmO_n^\epsilon(q)}R_{\bfT,\theta}^{\SO_n^\epsilon}.
\]
A character $\rho$ of $G$ is called \emph{unipotent} if
$\langle\rho,R_{\bfT,1}^\bfG\rangle_\bfG\neq0$ for some $\bfT$.
The set of irreducible unipotent character of $G$ is denoted by $\cale(G)_1$.

Let $\calv(G)^\sharp$ denote the subspace of $\calv(G)$ spanned by all Deligne-Lusztig characters of $G$.
For $f\in\calv(G)$, the orthogonal projection $f^\sharp$ of $f$ over $\calv(G)^\sharp$ is called the \emph{uniform projection} of $f$.
A class function $f$ is called \emph{uniform} if $f^\sharp=f$, i.e., if $f\in\calv(G)^\sharp$.

If $\bfG$ is connected, it is well-known that the regular character ${\rm Reg}_\bfG$ of $\bfG$ is uniform
(\cf.~\cite{carter-finite} corollary~7.5.6).
Because ${\rm Reg}_{\rmO^\epsilon}=\Ind_{\SO^\epsilon}^{\rmO^\epsilon}({\rm Reg}_{\SO^\epsilon})$,
we see that ${\rm Reg}_{\rmO^\epsilon}$ is also uniform.
Therefore,
we have
\begin{equation}\label{0330}
\rho(1)
=\langle\rho,{\rm Reg}_\bfG\rangle
=\langle\rho^\sharp,{\rm Reg}_\bfG\rangle
=\rho^\sharp(1).
\end{equation}
This means that $\rho^\sharp\neq0$ for any $\rho\in\cale(G)$.

Let $\chi$ be an irreducible character of $W_n$,
then we define
\begin{align}\label{0308}
\begin{split}
R^{\Sp_{2n}}_\chi
&=\frac{1}{|W_n|}\sum_{w\in W_n}\chi(w) R_{\bfT_w,\bf1}^{\Sp_{2n}};\\
R_\chi^{\rmO_{2n}^\epsilon}
&=\frac{1}{|W_n^\epsilon|}\sum_{w\in W_n^\epsilon}\chi(w) R_{\bfT_w,\bf1}^{\rmO_{2n}^\epsilon}.
\end{split}
\end{align}
Because $(\varepsilon_0\chi)(w)=\chi(w)$ if $\bfT_w\subset\rmO_{2n}^+$ and
$(\varepsilon_0\chi)(w)=-\chi(w)$ if $\bfT_w\subset\rmO_{2n}^-$ where
$\varepsilon_0\in\cale(W_n)$ is given in Subsection~\ref{0324},
we see that
\begin{equation}\label{0325}
R_{\varepsilon_0\chi}^{\rmO^+}=R_{\chi}^{\rmO^+}\quad\text{and}\quad
R_{\varepsilon_0\chi}^{\rmO^-}=-R_{\chi}^{\rmO^-}.
\end{equation}

\begin{remark}\label{0333}
Definition in (\ref{0308}) need some explanation when $\bfG=\rmO_{2n}^-$.
We identify $\bfW_\bfG=W_n^+$.
The action $F\colon W_n^+\rightarrow W_n^+$ induced by the Frobenius automorphism of $\bfG$
is given by $x\mapsto t_nx t_n$ where $t_n$ is given in Subsection~\ref{0324}.
If $\chi\in\cale(W_n)$,
then the restriction of $\chi$ to $W_n^+$ is fixed by $F$ (but might not be irreducible)
and therefore (3.17.1) in \cite{lg-CBMS} becomes
\begin{equation}\label{0314}
R^{\rmO_{2n}^-}_\chi=\frac{1}{|W_n^+|}\sum_{w\in W_n^+}\chi(w t_n)R_{\bfT_w,\bf1}^{\rmO_{2n}^-}.
\end{equation}
Now the mapping $w\mapsto w t_n$ gives a bijection between $W_n^+$ and $W_n^-$,
and we define
\[
R_{\bfT_{w t_n},\bf1}^{\rmO_{2n}^-}=R_{\bfT_w,\bf1}^{\rmO_{2n}^-}
\]
for $w\in W_n^+$.
So (\ref{0314}) becomes
\[
R^{\rmO_{2n}^-}_\chi=\frac{1}{|W_n^-|}\sum_{w\in W_n^-}\chi(w)R_{\bfT_w,\bf1}^{\rmO_{2n}^-}
\]
for any $\chi\in\cale(W_n)$.
\end{remark}

\begin{remark}\label{0328}
Note that the definition of $R^\bfG_\chi$ in \cite{amr} where $\bfG=\rmO^\epsilon_{2n}$ and $\chi\in\cale(W_n)$
is different from the definition here.
Let $R_\chi^{\rmO^\epsilon_{2n}}$ be defined as above,
and let ``$R_\chi^{\rmO^+_{2n}\oplus\rmO^-_{2n}}$'' be defined as in \cite{amr} p.367.
Then we have
\begin{align*}
R_\chi^{\rmO^+_{2n}\oplus\rmO^-_{2n}}
&=\frac{1}{|W_n|}\sum_{w\in W_n}\chi(w)R^{\rmO^{\epsilon_w}_{2n}}_{\bfT_w,\bf1}\\
&=\frac{1}{2}\biggl(\frac{1}{|W_n^+|}\sum_{w\in W_n^+}\chi(w)R^{\rmO^+_{2n}}_{\bfT_w,\bf1}
+\frac{1}{|W_n^-|}\sum_{w\in W_n^-}\chi(w)R^{\rmO^-_{2n}}_{\bfT_w,\bf1}\biggr)\\
&=\frac{1}{2}\Bigl(R_\chi^{\rmO^+_{2n}}+R_\chi^{\rmO^-_{2n}}\Bigr).
\end{align*}
\end{remark}

\subsection{Unipotent characters of $\Sp_{2n}$}\label{0335}
From \cite{lg} theorem 8.2,
there is a parametrization of unipotent characters $\cals_{\Sp_{2n}}\rightarrow\cale(\Sp_{2n})_1$
denoted by $\Lambda\mapsto\rho_\Lambda$.
Note that our definition of $\cals_{\Sp_{2n}}$ is different from that of $\Phi_n$ in \cite{lg} p.134.
However, we know that $\Lambda\in\cals_{n,\beta}$ if and only if $\Lambda^\rmt\in\cals_{n,-\beta}$,
and $\beta\equiv 3\pmod 4$ if and only if $-\beta\equiv 1\pmod 4$.
Therefore there is a natural bijection
\begin{equation}\label{0326}
\cals_{\Sp_{2n}}\longrightarrow\Phi_n\quad\text{by }
\Lambda\mapsto
\begin{cases}
\Lambda, & \text{if ${\rm def}(\Lambda)\geq 0$};\\
\Lambda^\rmt, & \text{if ${\rm def}(\Lambda)<0$}.
\end{cases}
\end{equation}

Every symbol $\Lambda$ of rank $n$ and defect $\beta\equiv 1\pmod 4$ is in $\cals_Z$ for some unique
special symbol $Z$ of rank $n$ and defect $1$, i.e.,
$\cals_{\Sp_{2n}}=\bigsqcup_Z\cals_Z$ where $Z$ runs over all special symbols of rank $n$ and defect $1$.
Therefore, we have the decomposition
\begin{equation}
\cale(\Sp_{2n})_1=\bigsqcup_{Z}\{\,\rho_\Lambda\mid\Lambda\in\cals_Z\,\}.
\end{equation}

It is know that there is a one-to-one correspondence between the set $\calp_2(n)$ of bipartitions of $n$ and
the set of irreducible characters of the Weyl group $W_n$ of $\Sp_{2n}$ (\cf.~\cite{Geck-Pfeiffer} theorem~5.5.6).
Then for a symbol $\Sigma\in\cals_{n,1}$,
we can associate it a uniform function $R_\Sigma$ on $\Sp_{2n}(q)$ given by $R_\Sigma=R_\chi$
where $R_\chi=R_\chi^\bfG$ is defined in (\ref{0308}) and $\chi$ is the irreducible character of $W_n$
associated to $\Upsilon(\Sigma)$ where $\Upsilon$ is the bijection $\cals_{n,1}\rightarrow\calp_2(n)$
given in (\ref{0219}).
The following proposition is modified from \cite{lg-symplectic} theorem 5.8:

\begin{proposition}[Lusztig]\label{0306}
Let $\bfG=\Sp_{2n}$, $Z$ a special symbol of rank $n$, defect $1$ and degree $\delta$.
For $\Sigma\in\cals_{Z,1}$,
we have
\[
\langle R_\Sigma,\rho_{\Lambda}\rangle_\bfG
=\begin{cases}
(-1)^{\langle\Sigma,\Lambda\rangle}2^{-\delta}, & \text{if $\Lambda\in\cals_Z$};\\
0, & \text{otherwise}
\end{cases}
\]
where
$\langle,\rangle\colon\cals_{Z,1}\times\cals_Z\rightarrow\bff_2$ is given by
$\langle\Lambda_N,\Lambda_M\rangle=|N\cap M|\pmod 2$.
\end{proposition}
\begin{proof}
Note that the unipotent character $\rho_\Lambda$ indexed by symbol $\Lambda$ of rank $n$ and
defect $1\pmod 4$ is indexed by its transpose $\Lambda^\rmt$ in \cite{lg} if ${\rm def}(\Lambda)<0$.
Suppose that $\Lambda=\Lambda_M$ for some $M\subset Z_\rmI$ with $|M^*|\equiv|M^*|\pmod 2$.
By (\ref{0221}), we know that $(\Lambda_M)^\rmt=\Lambda_{Z_\rmI\smallsetminus M}$.
So we need to check that
\begin{equation}\label{0309}
|M\cap N|\equiv|(Z_\rmI\smallsetminus M)\cap N|\pmod 2
\end{equation}
for any $N\subset Z_\rmI$ with $|N^*|=|N_*|$.
The requirement $|N^*|=|N_*|$ implies that $|N|$ is even.
Since $N\subset Z_\rmI$, (\ref{0309}) is obtained by the equality $|N|=|M\cap N|+|(Z_\rmI\smallsetminus M)\cap N|$.
\end{proof}

Because both sets $\{\,\rho_\Lambda\mid\Lambda\in\cals_Z\,\}$ and $\{\,R_\Sigma\mid\Sigma\in\cals_{Z,1}\,\}$ are orthonormal in $\calv(G)$,
the above proposition means that:
\begin{align}\label{0315}
\begin{split}
R_\Sigma &=\frac{1}{2^\delta}\sum_{\Lambda\in\cals_Z}(-1)^{\langle\Sigma,\Lambda\rangle}\rho_{\Lambda}\quad\text{for $\Sigma\in\cals_{Z,1}$};\\
\rho_\Lambda^\sharp &=\frac{1}{2^\delta}\sum_{\Sigma\in\cals_{Z,1}}(-1)^{\langle\Sigma,\Lambda\rangle}R_\Sigma\quad\text{for $\Lambda\in\cals_Z$}.
\end{split}
\end{align}

Let $\calv_Z$ denote the (complex) vector space spanned by $\{\,\rho_\Lambda\mid\Lambda\in\cals_Z\,\}$.
It is known that $\{\,\rho_\Lambda\mid\Lambda\in\cals_Z\,\}$ is an orthonormal basis for $\calv_Z$,
and $\{\,R_\Sigma\mid\Sigma\in\cals_{Z,1}\,\}$ is an orthonormal basis for the uniform projection $\calv_Z^\sharp$ of the space $\calv_Z$.

\subsection{Unipotent characters of $\rmO_{2n}^\epsilon$}\label{0331}
From \cite{lg} theorem 8.2,
we have the following parametrization of unipotent characters:
\begin{enumerate}
\item[(1)]
\begin{enumerate}
\item If $\Lambda\in\cals_{\rmO_{2n}^+}$ is a degenerate symbol,
then $n$ is even, $\Lambda=\Lambda^\rmt$ and
$\Lambda$ corresponds to two irreducible characters $\rho_1,\rho_2$ of $\SO_{2n}^+(q)$.
Now $\Ind_{\SO_{2n}^+(q)}^{\rmO_{2n}^+(q)}\rho_1\simeq\Ind_{\SO_{2n}^+(q)}^{\rmO_{2n}^+(q)}\rho_2$ is irreducible,
hence $\Lambda$ corresponds a unique irreducible unipotent character of $\rmO_{2n}^+(q)$.

\item If $\Lambda\in\cals_{\rmO^+_{2n}}$ is non-degenerate,
then $\Lambda\neq\Lambda^\rmt$ and both $\Lambda,\Lambda^\rmt$ corresponds to the same irreducible
unipotent character $\rho$ of $\SO_{2n}^+(q)$.
We know that $\Ind_{\SO^+_{2n}(q)}^{\rmO_{2n}^+(q)}\rho$ decomposed as a sum $\rho^\rmI\oplus\rho^{\rm II}$
of two irreducible unipotent characters of $\rmO_{2n}^+(q)$.
\end{enumerate}
\noindent Hence there is a one-to-one correspondence between $\cals_{\rmO_{2n}^+}$ and the set of
irreducible unipotent characters of $\rmO_{2n}^+(q)$.

\item[(2)] There is no degenerate symbol in $\cals_{\rmO_{2n}^-}$.
For non-degenerate symbols, the situation is similar to case (1.b).
Again, we have a one-to-one correspondence between $\cals_{\rmO_{2n}^-}$ to the set of
irreducible unipotent characters of $\rmO_{2n}^-(q)$.
\end{enumerate}
For $\Lambda\in\cals_{\rmO_{2n}^\epsilon}$,
the associated irreducible character of $\rmO_{2n}^\epsilon(q)$ is denoted by $\rho_\Lambda$.
It is known that $\rho_{\Lambda^\rmt}=\rho_\Lambda\cdot\sgn$.

Similar to the symplectic case,
for a symbol $\Sigma\in\cals_{n,0}$,
we define a uniform class function $R_\Sigma=R_\chi$ on $\rmO_{2n}^\epsilon(q)$
where $\chi$ is the irreducible character of $W_n$ associated to $\Upsilon(\Sigma)$.
If $\Sigma\in\cals_{n,0}$,
then $\Sigma^\rmt\in\cals_{n,0}$.
Moreover if the character $\chi$ of $W_n$ is associated to $\Upsilon(\Sigma)$,
then the character associated to $\Upsilon(\Sigma^\rmt)$ is $\varepsilon_0\chi$,
and hence
\[
R_{\Sigma^\rmt}^{\rmO^+}=R_\Sigma^{\rmO^+}\quad\text{and}\quad
R_{\Sigma^\rmt}^{\rmO^-}=-R_\Sigma^{\rmO^-}
\]
from (\ref{0325}).

\begin{remark}
If $\bfG=\rmO^-_{2n}$,
$R_\Sigma$ is denoted by $R(\Sigma-\Sigma^\rmt)$ in \cite{lg-orthogonal} p.745.
\end{remark}

Let $Z$ be a special symbol of defect $0$.
As in the symplectic case, we define $\langle,\rangle\colon\cals_{Z,0}\times\cals^\epsilon_Z\rightarrow\bff_2$ by
\begin{equation}\label{0318}
\langle\Lambda_M,\Lambda_N\rangle=|M\cap N|\pmod 2.
\end{equation}

\begin{lemma}\label{0313}
Let $Z$ be a special symbol of defect $0$,
$\Sigma\in\cals_{Z,0}$, $\Lambda\in\cals^\epsilon_Z$.
Then
\[
\langle\Sigma,\Lambda\rangle=\langle\Sigma,\Lambda^\rmt\rangle,\qquad
\langle\Sigma^\rmt,\Lambda\rangle
\begin{cases}
=\langle\Sigma,\Lambda\rangle, & \text{if $\epsilon=+$};\\
\neq\langle\Sigma,\Lambda\rangle, & \text{if $\epsilon=-$}.
\end{cases}
\]
\end{lemma}
\begin{proof}
Write $\Sigma=\Lambda_M$ and $\Lambda=\Lambda_N$ for some $M,N\subset Z_\rmI$ such that $|M^*|=|M_*|$.
In particular, $|M|$ is even.
Now $\Sigma^\rmt=\Lambda_{Z_\rmI\smallsetminus M}$ and $\Lambda^\rmt=\Lambda_{Z_\rmI\smallsetminus N}$
by (\ref{0221}).
Because $|M|=|M\cap N|+|M\cap(Z_\rmI\smallsetminus N)|$ and $|M|$ is even,
we conclude that $|M\cap N|\equiv|M\cap(Z_\rmI\smallsetminus N)|\pmod 2$,
and hence $\langle\Sigma,\Lambda\rangle=\langle\Sigma,\Lambda^\rmt\rangle$.

Note that $\Lambda\in\cals^+_Z$ if $|N|$ is even;
and $\Lambda\in\cals^-_Z$ if $|N|$ is odd.
Therefore
\[
|N\cap M|
\begin{cases}
\equiv|N\cap(Z_\rmI\smallsetminus M)|\pmod 2, & \text{if $\Lambda\in\cals^+_Z$};\\
\not\equiv|N\cap(Z_\rmI\smallsetminus M)|\pmod 2, & \text{if $\Lambda\in\cals^-_Z$},
\end{cases}\]
and the lemma is proved.
\end{proof}

The following proposition is a modification for $\rmO^\epsilon_{2n}(q)$
from \cite{lg-orthogonal} theorem 3.15:

\begin{proposition}\label{0319}
Let $\bfG=\rmO^\epsilon_{2n}$,
$Z$ be a non-degenerate special symbol of defect $0$ and degree $\delta\geq 1$.
For any $\Sigma\in\cals_{Z,0}$,
we have
\[
\langle R^{\rmO^\epsilon}_\Sigma,\rho^{\rmO^\epsilon}_\Lambda\rangle_{\rmO^\epsilon}
=\begin{cases}
(-1)^{\langle\Sigma,\Lambda\rangle}2^{-(\delta-1)}, & \text{if $\Lambda\in\cals^\epsilon_Z$};\\
0, & \text{otherwise}.
\end{cases}
\]
\end{proposition}
\begin{proof}
Recall that
\[
\langle R^{\SO^\epsilon}_\Sigma,\rho^{\SO^\epsilon}_\Lambda\rangle_{\SO^\epsilon}
=\begin{cases}
(-1)^{\langle\Sigma,\Lambda\rangle}2^{-(\delta-1)}, & \text{if $\Lambda\in\cals^\epsilon_Z$};\\
0, & \text{otherwise}.
\end{cases}
\]
from \cite{lg-orthogonal} theorem 3.15.
Moreover, we have $\Ind_{\SO^\epsilon}^{\rmO^\epsilon}\rho_\Lambda^{\SO^\epsilon}
=\rho_\Lambda^{\rmO^\epsilon}+\rho^{\rmO^\epsilon}_{\Lambda^\rmt}$ (since $Z$ is non-degenerate) and
$R_\Sigma^{\rmO^\epsilon}=\Ind_{\SO^\epsilon}^{\rmO^\epsilon}R_\Sigma^{\SO^\epsilon}$.
Hence, $\rho^{\rmO^\epsilon}_\Lambda|_{\SO^\epsilon}=\rho^{\rmO^\epsilon}_{\Lambda^\rmt}|_{\SO^\epsilon}
=\rho^{\SO^\epsilon}_\Lambda$,
and
\[
R^{\rmO^\epsilon}_\Sigma(y)=\begin{cases}
2R^{\SO^\epsilon}_\Sigma(y), & \text{if $y\in\SO^\epsilon(q)$};\\
0, & \text{if $y\in\rmO^\epsilon(q)\smallsetminus\SO^\epsilon(q)$}.
\end{cases}
\]
Note that $\langle\Sigma,\Lambda\rangle=\langle\Sigma,\Lambda^\rmt\rangle$ by Lemma~\ref{0313}.
Therefore
\[
\langle R^{\rmO^\epsilon}_\Sigma,\rho^{\rmO^\epsilon}_\Lambda\rangle_{\rmO^\epsilon}
=\langle R^{\rmO^\epsilon}_\Sigma,\rho^{\rmO^\epsilon}_{\Lambda^\rmt}\rangle_{\rmO^\epsilon}
=\frac{|\SO^\epsilon(q)|}{|\rmO^\epsilon(q)|}\langle 2R^{\SO^\epsilon}_\Sigma,\rho^{\SO^\epsilon}_\Lambda\rangle_{\SO^\epsilon}
=\langle R^{\SO^\epsilon}_\Sigma,\rho^{\SO^\epsilon}_\Lambda\rangle_{\SO^\epsilon},
\]
and hence the lemma is proved.
\end{proof}

Because $\{\,\rho_\Lambda\mid\Lambda\in\cals^\epsilon_Z\,\}$ is an orthonormal set in $\calv(G)$,
the proposition implies that
\begin{equation}\label{0316}
R^{\rmO^\epsilon}_\Sigma
=\frac{1}{2^{\delta-1}}\sum_{\Lambda\in\cals^\epsilon_Z}
(-1)^{\langle\Sigma,\Lambda\rangle}\rho^{\rmO^\epsilon}_\Lambda\text{ \ for $\Sigma\in\cals_{Z,0}$},
\end{equation}

\begin{corollary}\label{0317}
If $G$ is a nontrivial orthogonal group, then
\[
\langle R_\Sigma, R_\Sigma\rangle_\bfG=\begin{cases}
1, & \text{if\/ $\Sigma$ is degenerate and $\bfG=\rmO^+$;}\\
2, & \text{if\/ $\Sigma$ is non-degenerate}.
\end{cases}
\]
\end{corollary}
\begin{proof}
If $\Sigma$ is degenerate,
then $R_\Sigma^{\rmO^+}=\Ind_{\SO^+}^{\rmO^+}R_\Sigma^{\SO^+}$ is an irreducible
character of $\rmO^+$,
and hence $\langle R^{\rmO^+}_\Sigma,R^{\rmO^+}_\Sigma\rangle_\bfG=1$.

Next, suppose that $\Sigma$ is non-degenerate.
Then
\begin{align*}
\langle R^{\rmO^\epsilon}_\Sigma,R^{\rmO^\epsilon}_\Sigma\rangle_\bfG
& =\frac{1}{2^{2(\delta-1)}}\sum_{\Lambda\in\cals^\epsilon_Z}
\sum_{\Lambda'\in\cals^\epsilon_Z}
(-1)^{\langle\Sigma,\Lambda\rangle+\langle\Sigma,\Lambda'\rangle}
\langle\rho_\Lambda^{\rmO^\epsilon},\rho_{\Lambda'}^{\rmO^\epsilon}\rangle_\bfG  =\frac{1}{2^{2(\delta-1)}}\sum_{\Lambda\in\cals^\epsilon_Z}1
=2.
\end{align*}
\end{proof}

As in the symplectic case,
let $\calv^\epsilon_Z$ denote the (complex) vector space spanned by $\{\,\rho_\Lambda\mid\Lambda\in\cals^\epsilon_Z\,\}$.
Then $\{\,\rho_\Lambda\mid\Lambda\in\cals^\epsilon_Z\,\}$ is an orthonormal basis for $\calv^\epsilon_Z$.
If $\Sigma$ is degenerate,
then $Z=\Sigma$, $\cals^-_Z=\emptyset$, $\cals^+_Z=\{Z\}$ and $\rho_\Sigma^{\rmO^+}=R_\Sigma^{\rmO^+}$.
If $\Sigma$ is non-degenerate,
let $\bar\cals_{Z,0}$ denote a complete set of representatives of cosets $\{\Sigma,\Sigma^\rmt\}$ in $\cals_{Z,0}$,
then $\{\,\frac{1}{\sqrt 2}R^{\rmO^\epsilon}_\Sigma\mid\Sigma\in\bar\cals_{Z,0}\,\}$ is an orthonormal basis for $(\calv_Z^\epsilon)^\sharp$.
Note that $R^{\rmO^\epsilon}_{\Sigma^\rmt}=\epsilon R^{\rmO^\epsilon}_\Sigma$ and
$(-1)^{\langle\Sigma^\rmt,\Lambda\rangle}=\epsilon(-1)^{\langle\Sigma,\Lambda\rangle}$
for $\Lambda\in\cals_Z^\epsilon$ by Lemma~\ref{0313},
so we have
\begin{equation}\label{0327}
(\rho^{\rmO^\epsilon}_\Lambda)^\sharp
=(\rho^{\rmO^\epsilon}_{\Lambda^\rmt})^\sharp
=\frac{1}{2^\delta}\sum_{\Sigma\in\bar\cals_{Z,0}}(-1)^{\langle\Sigma,\Lambda\rangle}R^{\rmO^\epsilon}_\Sigma
=\frac{1}{2^{\delta+1}}\sum_{\Sigma\in\cals_{Z,0}}(-1)^{\langle\Sigma,\Lambda\rangle}R^{\rmO^\epsilon}_\Sigma
\end{equation}
for $\Lambda\in\cals^\epsilon_Z$.

\subsection{The decomposition of the Weil character}\label{0340}
In this subsection, let $(\bfG,\bfG')$ be a reductive dual pair.
By the homomorphism
$G\times G'\rightarrow GG'\hookrightarrow\Sp_{2N}(q)$,
the Weil character $\omega^\psi_{\Sp_{2N}}$ of $\Sp_{2N}(q)$ with respect to a nontrivial character $\psi$
can be regarded as an character $\omega^\psi_{\bfG,\bfG'}$ of $G\times G'$.
Define the uniform class function $\omega^\natural_{\bfG,\bfG'}$ of $G\times G'$ as follows:
\begin{enumerate}
\item[(I)] for $(\GL_n,\GL_{n'})$
\begin{multline*}
\omega^\natural_{\GL_n,\GL_{n'}}
=\sum_{k=0}^{\min(n,n')} \frac{1}{|\bfW_{\GL_k}|}\frac{1}{|\bfW_{\GL_{n-k}}|}\frac{1}{|\bfW_{\GL_{n'-k}}|}\sum_{v\in\bfW_{\GL_k}} \\
\sum_{\theta\in\cale(T_v)}\sum_{w\in\bfW_{\GL_{n-k}}}\sum_{w'\in\bfW_{\GL_{n'-k}}}
R^{\GL_n}_{\bfT_v\times\bfT_w,\theta\otimes\bf 1}\otimes R^{\GL_{n'}}_{\bfT_v\times\bfT_{w'},\theta\otimes\bf 1};
\end{multline*}

\item[(II)] for $(\rmU_n,\rmU_{n'})$
\begin{multline*}
\omega^\natural_{\rmU_n,\rmU_{n'}}
=\sum_{k=0}^{\min(n,n')} \frac{(-1)^k}{|\bfW_{\rmU_k}|}\frac{1}{|\bfW_{\rmU_{n-k}}|}\frac{1}{|\bfW_{\rmU_{n'-k}}|}\sum_{v\in\bfW_{\rmU_k}} \\
\sum_{\theta\in\cale(T_v)}\sum_{w\in\bfW_{\rmU_{n-k}}}\sum_{w'\in\bfW_{\rmU_{n'-k}}}
R^{\rmU_n}_{\bfT_v\times\bfT_w,\theta\otimes\bf 1}\otimes R^{\rmU_{n'}}_{\bfT_v\times\bfT_{w'},\theta\otimes\bf 1};
\end{multline*}

\item[(III)] for $(\Sp_{2n},\SO_{2n'}^\epsilon)$
\begin{itemize}
\item if $n'>n$,
\begin{multline*}
\omega^\natural_{\Sp_{2n},\SO^\epsilon_{2n'}}
=\sum_{k=0}^n \frac{1}{|\bfW_{\Sp_{2k}}|}\frac{1}{|\bfW_{\Sp_{2(n-k)}}|}\frac{1}{|\bfW_{\SO_{2(n'-k)}}|}\sum_{v\in\bfW_k}
\sum_{\theta\in\cale(T_v)} \\
\sum_{w\in\bfW_{\Sp_{2(n-k)}}}\sum_{w'\in\bfW_{\SO_{2(n'-k)}^{\epsilon_v\epsilon}}}
\epsilon_v R^{\Sp_{2n}}_{\bfT_v\times\bfT_w,\theta\otimes\bf 1}\otimes
R^{\SO^{\epsilon}_{2n'}}_{\bfT_v\times\bfT_{w'},\theta\otimes\bf 1};
\end{multline*}
\item if $n'\leq n$,
\begin{multline*}
\omega^\natural_{\Sp_{2n},\SO^\epsilon_{2n'}}
=\sum_{k=0}^{n'-1}\frac{1}{|\bfW_{\Sp_{2k}}|}\frac{1}{|\bfW_{\Sp_{2(n-k)}}|}\frac{1}{|\bfW_{\SO_{2(n'-k)}}|}\sum_{v\in W_k}
\sum_{\theta\in\cale(T_v)} \\
\sum_{w\in\bfW_{\Sp_{2(n-k)}}}\sum_{w'\in\bfW_{\SO_{2(n'-k)}^{\epsilon_v\epsilon}}}
\epsilon_v R^{\Sp_{2n}}_{\bfT_v\times\bfT_w,\theta\otimes\bf 1}\otimes R^{\SO^{\epsilon}_{2n'}}_{\bfT_v\times\bfT_{w'},\theta\times\bf 1}\\
+\frac{1}{|W_{n'}^\epsilon|}\frac{1}{|\bfW_{\Sp_{2(n-n')}}|}
\sum_{v\in W_{n'}^\epsilon}\sum_{\theta\in\cale(T_v)}\sum_{w\in\bfW_{\Sp_{2(n-n')}}}
\epsilon R^{\Sp_{2n}}_{\bfT_v\times\bfT_w,\theta\otimes\bf 1}\otimes R^{\SO^\epsilon_{2n'}}_{\bfT_v,\theta};
\end{multline*}
\end{itemize}

\item[(IV)] for $(\Sp_{2n},\SO_{2n'+1})$
\begin{multline*}
\omega^\natural_{\Sp_{2n},\SO_{2n'+1}}
= \sum_{k=0}^{\min(n,n')}\frac{1}{|W_k|}\frac{1}{|\bfW_{\Sp_{2(n-k)}}|}\frac{1}{|\bfW_{\SO_{2(n'-k)+1}}|}
\sum_{v\in W_k}\sum_{\theta\in\cale(T_v)}\\
\sum_{w\in\bfW_{\Sp_{2(n-k)}}}\sum_{w'\in\bfW_{\SO_{2(n'-k)+1}}}
\epsilon_w R^{\Sp_{2n}}_{\bfT_v\times\bfT_w,\theta\otimes\theta_w}\otimes R^{\SO_{2n'+1}}_{\bfT_v\times\bfT_{w'},\theta\otimes\theta_{w'}}.
\end{multline*}
\end{enumerate}

The following result is the main theorem of \cite{srinivasan}.

\begin{proposition}[Srinivasan]\label{0336}
Assume that $q$ is sufficiently large, such that every torus $T$ in $G$ or $G'$
satisfies the condition that $T/Z$ has at least two $W(T)$-orbits of characters in general position,
where $Z$ is the center of the group.
Then
\begin{enumerate}
\item[(i)] $\omega^\psi_{\GL_n,\GL_{n'}}=\omega^\natural_{\GL_n,\GL_{n'}}$;

\item[(ii)] $(-1)^{nn'}\omega^\psi_{\rmU_n,\rmU_{n'}}=\omega^\natural_{\rmU_n,\rmU_{n'}}$;

\item[(iii)] $\omega^\sharp_{\Sp_{2n},\SO^\epsilon_{2n'}}=\omega^\natural_{\Sp_{2n},\SO^\epsilon_{2n'}}$.
\end{enumerate}
\end{proposition}

Based on Proposition~\ref{0336},
the following proposition is proved in \cite{pan-odd}:

\begin{proposition}\label{0337}
Suppose that $q$ is sufficiently large so that Proposition~\ref{0336} holds.
Then
\[
\omega^\sharp_{\Sp_{2n},\SO_{2n'+1}}\cdot({\bf 1}\otimes\chi_{\SO_{2n'+1}})
=\omega^\natural_{\Sp_{2n},\SO_{2n'+1}}.
\]
\end{proposition}

Now we show that the above two propositions are still true without assuming that $q$ is sufficiently large.
First we need some result on the Green functions of classical groups.
Let $\bfG$ be a connected classical group, and let $u$ be a unipotent element of $G$.
It is known that the value $R_{\bfT,\theta}^\bfG(u)$ is independent of the character $\theta$,
and the restriction of $R_{\bfT,\bf1}^\bfG$ to the set of unipotent elements of $G$ is called the
\emph{Green function} of $\bfG$ associated to the rational maximal torus $\bfT$ (\cf.~\cite{carter-finite}).
It is known that the value $R^\bfG_{\bfT,\bf1}(u)$ is a polynomial in $q$ for any classical group $\bfG$:
\begin{itemize}
\item If $\bfG$ is a general linear group,
the result is proved in \cite{green} p.420.

\item If $\bfG$ is a unitary group,
the result is proved in \cite{kawanaka} theorem 4.1.2.

\item If $\bfG$ is a symplectic group or a special orthogonal group,
the result is proved in \cite{srinivasan-polynomial} p.1242 by using the result in \cite{Kazhdan} and assuming that
$q$ and $p$ are large enough.
The restriction on $p$ is removed in \cite{lusztig-green-function} p.355.
Finally the restriction on $q$ is removed in \cite{shoji-polynomial} theorem 4.2 when the center of $\bfG$ is connected.
Note that in \cite{dl} p.123, it is known that the map $x\mapsto\bar x$ from $\bfG\rightarrow\bfG^{\rm ad}$
gives a bijection on the sets of unipotent elements and we have
\[
R^\bfG_{\bfT,\bf1}(u)
=R^{\bfG^{\rm ad}}_{\bar\bfT,\bf1}(\bar u).
\]
So the result is true for any symplectic group or special orthogonal group.

\item If $\bfG$ is an orthogonal group,
note that $R_{\bfT,\bf1}^{\rmO_n^\epsilon}=\Ind_{\SO_n^\epsilon(q)}^{\rmO_n^\epsilon(q)}R_{\bfT,\bf1}^{\SO_n^\epsilon}$,
we see that $R_{\bfT,\bf1}^{\rmO_n^\epsilon}(u)$ is again a polynomial in $q$.
\end{itemize}

\begin{theorem}\label{0338}
Suppose that $p\neq 2$.
Then
\begin{enumerate}
\item[(i)] $\omega^\psi_{\GL_n,\GL_{n'}}=\omega^\natural_{\GL_n,\GL_{n'}}$;

\item[(ii)] $(-1)^{nn'}\omega^\psi_{\rmU_n,\rmU_{n'}}=\omega^\natural_{\rmU_n,\rmU_{n'}}$;

\item[(iii)] $\omega^\sharp_{\Sp_{2n},\SO^\epsilon_{2n'}}=\omega^\natural_{\Sp_{2n},\SO^\epsilon_{2n'}}$;

\item[(iv)] $\omega^\sharp_{\Sp_{2n},\SO_{2n'+1}}\cdot({\bf 1}\otimes\chi_{\SO_{2n'+1}})
=\omega^\natural_{\Sp_{2n},\SO_{2n'+1}}$.
\end{enumerate}
\end{theorem}
\begin{proof}
Let $(\bfG,\bfG')$ be a finite reductive dual pair.
Define the class function $\Phi$ on $G\times G'$ by
\[
\Phi=\begin{cases}
\omega^\psi_{\GL_n,\GL_{n'}}-\omega^\natural_{\GL_n,\GL_{n'}}, & \text{if $(\bfG,\bfG')=(\GL_n,\GL_{n'})$};\\
(-1)^{nn'}\omega^\psi_{\rmU_n,\rmU_{n'}}-\omega^\natural_{\rmU_n,\rmU_{n'}}, & \text{if $(\bfG,\bfG')=(\rmU_n,\rmU_{n'})$};\\
\omega^\sharp_{\Sp_{2n},\SO^\epsilon_{2n'}}-\omega^\natural_{\Sp_{2n},\SO^\epsilon_{2n'}}, & \text{if $(\bfG,\bfG')=(\Sp_{2n},\SO^\epsilon_{2n'})$};\\
\omega^\sharp_{\Sp_{2n},\SO_{2n'+1}}\cdot({\bf 1}\otimes\chi_{\SO_{2n'+1}})-\omega^\natural_{\Sp_{2n},\SO_{2n'+1}},
& \text{if $(\bfG,\bfG')=(\Sp_{2n},\SO_{2n'+1})$}.
\end{cases}
\]
So we need to show that $\Phi$ is the zero function.
We know that $\Phi$ is a uniform class function on $G\times G'$,
so we can write
\begin{equation}
\Phi=\sum_{(\bfT,\theta)}\sum_{(\bfT',\theta')}a_{\bfT,\theta,\bfT',\theta'}R^\bfG_{\bfT,\theta}\otimes R^{\bfG'}_{\bfT',\theta'}
\end{equation}
where $a_{\bfT,\theta,\bfT',\theta'}\in\bbC$.
Now let $y\in G$ and $y'\in G'$ with respective Jordan decompositions $y=su$ and $y'=s'u'$.
\begin{enumerate}
\item[(1)] Suppose $s,s'$ are not both central in $G,G'$ respectively.
It is proved in \cite{srinivasan} p.150--p.151 that $\Phi(y,y')=0$ for Cases (I),(II), and (III) by the induction hypothesis.
Similarly, we can also show that $\Phi(y,y')=0$ for Case (IV) (\cf.~\cite{pan-odd} subsection 6.4).

\item[(2)] Suppose that $s,s'$ are in the centers of $G,G'$ respectively.
Then we know that (\cf.~\cite{carter-finite} proposition 7.2.8)
\[
R^\bfG_{\bfT,\theta}(su)=\theta(s)R^\bfG_{\bfT,\bf1}(u)=\theta(s)Q^\bfG_\bfT(u)
\]
where $Q^\bfG_\bfT$ is the \emph{Green function} defined on the set of unipotent elements in $G$.
It is known that the value $Q^\bfG_\bfT(u)$ is a complex polynomial in $q$ (\cf.~the paragraph before the theorem).
Therefore $\Phi(y,y')$ is a complex polynomial in $q$.
But from Proposition~\ref{0336} and Proposition~\ref{0337} we know that
$\Phi(y,y')$ equals to $0$ if $q$ is sufficiently large.
This means that $\Phi(y,y')$ is in fact the zero polynomial in $q$,
and hence $\Phi(y,y')=0$ for any $\bff_q$ of odd characteristic.
\end{enumerate}
\end{proof}

\subsection{The main results}\label{0322}
Let $(\bfG,\bfG')=(\Sp_{2n},\rmO^\epsilon_{2n'})$.
It is known that unipotent characters are preserved in the Howe correspondence for the dual pair $(G,G')$.
Let $\omega_{\bfG,\bfG',1}$ denote the unipotent part of the Weil character $\omega^\psi_{\bfG,\bfG'}$.
Then \cite{amr} proposition 4.1 can be written as follows (\cf.~\cite{Kable-Sanat} proposition 5.3):

\begin{proposition}\label{0329}
Let $(\bfG,\bfG')=(\Sp_{2n},\rmO^\epsilon_{2n'})$.
Then we have the decomposition
\[
\omega_{\bfG,\bfG',1}^\sharp
=\frac{1}{2}\sum_{k=0}^{\min(n,n')}\sum_{\chi\in\cale(W_k)}
R_{I_{n,k}(\chi)}^\bfG
\otimes R_{I_{n',k}(\varepsilon_0\chi)}^{\bfG'}
\]
where $I_{n,k}(\chi)=\Ind^{W_n}_{W_k\times W_{n-k}}(\chi\otimes{\bf 1})$,
$I_{n',k}(\varepsilon_0\chi)
=\Ind^{W_{n'}}_{W_k\times W_{n'-k}}(\varepsilon_0\chi\otimes{\bf 1})$,
and $\varepsilon_0$ the irreducible character of $W_k$ given in Subsection~\ref{0324}.
\end{proposition}

\begin{remark}
\begin{enumerate}
\item[(1)] In \cite{amr} p.369, they consider the decomposition of
\[
\omega^\sharp_{\Sp_{2n},\rmO_{2n'}^+,1}+\omega^\sharp_{\Sp_{2n},\rmO_{2n'}^-,1}
\]
at the same time but here we consider the decompositions of $\omega^\sharp_{\Sp_{2n},\rmO_{2n'}^+,1}$
and $\omega^\sharp_{\Sp_{2n},\rmO_{2n'}^-,1}$ separately.

\item[(2)] The appearance of the factor ``$\frac{1}{2}$'' in Proposition~\ref{0329} is due to the different
definitions of $R_\chi^{\rmO^\epsilon_{2n'}}$ (\cf.~Remark~\ref{0328}).
\end{enumerate}
\end{remark}

The irreducible character of $W_n$ associated to the bi-partition $\sqbinom{\lambda}{\mu}$ is denoted by $\chi_{\sqbinom{\lambda}{\mu}}$.
The following special case of the Littlewood-Richardson rule is well-known (\cf.~\cite{Geck-Pfeiffer} lemma 6.1.3)
\begin{equation}
\Ind^{W_n}_{W_k\times W_{n-k}}(\chi_{\sqbinom{\lambda'}{\mu'}}\otimes{\bf 1})
=\sum_{|\lambda|=|\lambda'|+(n-k),\ \lambda'^\rmT\preccurlyeq\lambda^\rmT}\chi_{\sqbinom{\lambda}{\mu'}}
\end{equation}
where $\lambda'^\rmT\preccurlyeq\lambda^\rmT$ is defined in Subsection~\ref{0220}.

For a special symbol $Z$ of rank $n$ and defect $1$,
and a special symbol $Z'$ of rank $n'$ and defect $0$,
let $\omega_{Z,Z'}$ denote the orthogonal projection of $\omega_{\bfG,\bfG',1}$ over $\calv_Z\otimes\calv^\epsilon_{Z'}$, i.e.,
$\omega_{Z,Z'}$ is the sum of irreducible characters in the set
\[
\{\,\rho_\Lambda\otimes\rho_{\Lambda'}\in\omega_{\bfG,\bfG',1}\mid\Lambda\in\cals_Z,\ \Lambda'\in\cals^\epsilon_{Z'}\,\}
\]
Then by Proposition~\ref{0306} and Proposition~\ref{0319} we have
\begin{equation}\label{0341}
\omega_{\bfG,\bfG',1}=\sum_{Z,Z'}\omega_{Z,Z'}\quad\text{ and }\quad
\omega_{\bfG,\bfG',1}^\sharp=\sum_{Z,Z'}\omega^\sharp_{Z,Z'}
\end{equation}
where $Z,Z'$ run over all special symbols of rank $n,n'$ and defect $1,0$ respectively.
Then \cite{amr} th\'eor\`eme 4.4 can be rephrased as follows:
\begin{equation}\label{0312}
\omega_{Z,Z'}^\sharp
=\frac{1}{2}\sum_{(\Sigma,\Sigma')\in\cald_{Z,Z'}} R_\Sigma\otimes R_{\Sigma'}
\end{equation}
where $\cald_{Z,Z'}$ is defined in (\ref{0227}).
In particular,
if $\cald_{Z,Z'}=\emptyset$,
then $\omega_{Z,Z'}^\sharp=0$ and hence $\omega_{Z,Z'}=0$ (\cf.~(\ref{0330})).

The following theorem is the main result of this work:
\begin{theorem}\label{0310}
Let $(\bfG,\bfG')=(\Sp_{2n},\rmO^\epsilon_{2n'})$,
$Z,Z'$ special symbols of rank $n,n'$ and defect $1,0$ respectively, and $\epsilon=\pm$.
Then
\begin{equation}\label{0342}
\frac{1}{2}\sum_{(\Sigma,\Sigma')\in\cald_{Z,Z'}} R^\bfG_\Sigma\otimes R^{\bfG'}_{\Sigma'}
=\sum_{(\Lambda,\Lambda')\in\calb^\epsilon_{Z,Z'}}\rho_\Lambda^\sharp\otimes\rho^\sharp_{\Lambda'}.
\end{equation}
\end{theorem}
The proof of the theorem are divided into several stages:
\begin{itemize}
\item We first prove the theorem for the two basic cases that 
$(Z,Z')=(Z_{(m)},Z'_{(m+1)})$, or $(Z_{(m)},Z'_{(m)})$ in 
Proposition~\ref{0512} and Proposition~\ref{0513} respectively where
$Z_{(m)}$ and $Z'_{(m)}$ are given in Example~\ref{0225}.
For these two cases, the result is proved by applying directly Lusztig's theory on 
unipotent characters.

\item Then the theorem is proved in Proposition~\ref{0806} for the case that $\epsilon=+$,
both $Z,Z'$ are regular and $\cald_{Z,Z'}$ is one-to-one.
For this case, we reduce the situation to the basic two cases considered above via
an isometry of inner product spaces.

\item And the theorem is proved in Proposition~\ref{0813} for the case $\epsilon=+$.
We will define subsets of pairs $\Psi_0,\Psi'_0$ of $Z_\rmI,Z'_\rmI$ respectively to measure the 
complexity of the relation $\cald_{Z,Z'}$. 
In particular, $\cald_{Z,Z'}$ is one-to-one when both $\Psi_0,\Psi'_0$ are empty.
Then we will prove the theorem iteratively by reducing the sizes of $\Psi_0,\Psi'_0$
via inner product space isometries.

\item Finally we prove the theorem in Proposition~\ref{0921} for case $\epsilon=-$.
\end{itemize}

\begin{corollary}
Let $(\bfG,\bfG')=(\Sp_{2n},\rmO^\epsilon_{2n'})$.
Then
\[
\omega^\sharp_{\bfG,\bfG',1}
=\sum_{(\Lambda,\Lambda')\in\calb_{\bfG,\bfG'}}\rho_\Lambda^\sharp\otimes\rho_{\Lambda'}^\sharp
\]
\end{corollary}
\begin{proof}
This follows from (\ref{0341}), (\ref{0312}), Theorem~\ref{0310}, and (\ref{0227}) directly.
\end{proof}


\section{Correspondence for Cuspidal Symbols}\label{0403}

\subsection{Basic case I: from $\Sp_{2m(m+1)}$ to $\rmO^\epsilon_{2(m+1)^2}$}\label{0501}
In this subsection,
we assume that $\bfG=\Sp_{2m(m+1)}$ and $\bfG'=\rmO^\epsilon_{2(m+1)^2}$ where $m$ is a non-negative integer,
and we fix $Z=Z_{(m)}$ and $Z'=Z'_{(m+1)}$ (\cf.~Example~\ref{0225}).
Now $Z$ is a regular special symbol of rank $2m(m+1)$ and defect $1$,
and $Z'=\binom{2m+1}{-}\cup Z^\rmt$ is a regular special symbol of rank $2(m+1)^2$ and defect $0$.

It is easy to check that $|\cals_{Z,1}|=\binom{2m+1}{m}$ and $|\cals_{Z',0}|=\binom{2m+2}{m+1}=2\times\binom{2m+1}{m}$
where $\binom{k}{l}=\frac {k!}{l!(k-l)!}$ denotes the binomial coefficient.
Define
\begin{align}\label{0502}
\begin{split}
\theta=\theta_{Z,Z'}\colon\cals_{Z,1} &\longrightarrow\cals_{Z',0}\quad\text{ by }
\Sigma\mapsto\textstyle\binom{2m+1}{-}\cup\Sigma^\rmt, \\
\theta^\epsilon\colon\cals_Z &\longrightarrow\cals^\epsilon_{Z'}
\quad\text{ by }\Lambda\mapsto\begin{cases}
\binom{2m+1}{-}\cup\Lambda^\rmt, & \text{if $\epsilon=+$};\\
\binom{-}{2m+1}\cup\Lambda^\rmt, & \text{if $\epsilon=-$}.
\end{cases}
\end{split}
\end{align}
It is clear that the mapping $\theta\colon\cals_{Z,1}\rightarrow\cals_{Z',0}$ is injective and
$\theta=\theta^+|_{\cals_{Z,1}}$.
Moreover, for each $\Sigma'\in\cals_{Z',0}$
exactly one of $\Sigma',\Sigma'^\rmt$ is in the image of $\theta$,
i.e., $\theta(\cals_{Z,1})$ forms a complete set of representatives of cosets $\{\Sigma',\Sigma'^\rmt\}$ in $\cals_{Z',0}$.

\begin{lemma}\label{0507}
Let $Z=Z_{(m)}$ and $Z'=Z'_{(m+1)}$.
Then
\[
\omega_{Z,Z'}^\sharp
=\frac{1}{2}\sum_{\Sigma\in\cals_{Z,1}} R_\Sigma\otimes R_{\theta(\Sigma)}.
\]
\end{lemma}
\begin{proof}
From (\ref{0312}), we need to prove that
$\cald_{Z,Z'}=\{\,(\Sigma,\theta(\Sigma))\mid\Sigma\in\cals_{Z,1}\,\}$,
which is just Example~\ref{0225}.
\end{proof}

\begin{lemma}\label{0510}
Let $M$ be a subset of $Z_\rmI$.
Then
\[
\theta^\epsilon(\Lambda_M)=\begin{cases}
\Lambda_{M^\rmt}, & \text{if $\epsilon=+$};\\
\Lambda_{\binom{-}{2m+1}\cup M^\rmt}, & \text{if $\epsilon=-$}.
\end{cases}
\]
\end{lemma}
\begin{proof}
From (\ref{0201}) and (\ref{0502}),
we have
\begin{align*}
\theta^+(\Lambda_M)
=\theta^+((Z\smallsetminus M)\cup M^\rmt)
=\textstyle\binom{2m+1}{-}\cup(Z^\rmt\smallsetminus M^\rmt)\cup M
=(Z'\smallsetminus M^\rmt)\cup M 
=\Lambda_{M^\rmt},
\end{align*}
and
\begin{align*}
\textstyle\theta^-(\Lambda_M)
=\theta^-((Z\smallsetminus M)\cup M^\rmt)
&=\textstyle\binom{-}{2m+1}\cup(Z^\rmt\smallsetminus M^\rmt)\cup M \\
&=\textstyle(Z'\smallsetminus (\binom{2m+1}{-}\cup M)^\rmt)\cup(\binom{2m+1}{-}\cup M)
=\textstyle\Lambda_{\binom{-}{2m+1}\cup M^\rmt}.
\end{align*}
\end{proof}

\begin{lemma}\label{0505}
We have $\langle\Sigma,\Lambda\rangle=\langle\theta(\Sigma),\theta^\epsilon(\Lambda)\rangle$
for any $\Sigma\in\cals_{Z,1}$ and $\Lambda\in\cals_Z$
where $\langle,\rangle$ is defined in Proposition~\ref{0306}.
\end{lemma}
\begin{proof}
Write $\Lambda=\Lambda_M$, $\Sigma=\Lambda_N$, $\theta^\epsilon(\Lambda)=\Lambda_{M'}$,
$\theta(\Sigma)=\Lambda_{N'}$ for some $M,N\subset Z_\rmI$,
and some $M',N'\subset Z'_\rmI$.
Suppose that $\epsilon=+$.
By Lemma~\ref{0510}, we have $M'=M^\rmt$ and $N'=N^\rmt$.
Thus from the definition in Proposition~\ref{0306}, we have
\begin{align*}
\langle\Sigma,\Lambda\rangle
& =|N\cap M|\pmod 2 \\
& \equiv |N^\rmt\cap M^\rmt|\pmod 2
=\langle\theta(\Sigma),\theta^\epsilon(\Lambda)\rangle.
\end{align*}

Next, suppose that $\epsilon=-$.
Then we see that $N'=N^\rmt$ and $M'=\binom{-}{2m+1}\cup M^\rmt$.
Note that $2m+1\not\in N$,
so
\begin{align*}
\langle\Sigma,\Lambda\rangle
& = |N\cap M|\pmod 2 \\
& \equiv |N^\rmt\cap(\textstyle\binom{-}{2m+1}\cup M^\rmt)|\pmod 2
=\langle\theta(\Sigma),\theta^\epsilon(\Lambda)\rangle.
\end{align*}
\end{proof}

\begin{lemma}\label{0506}
Let $\Sigma,\Sigma'\in\cals_{Z,1}$.
Then
\[
\sum_{\Lambda\in\cals_Z}(-1)^{\langle\Sigma,\Lambda\rangle+\langle\Sigma',\Lambda\rangle}
=\begin{cases}
2^{2m}, & \text{if\/ $\Sigma=\Sigma'$};\\
0, & \text{otherwise}.
\end{cases}
\]
\end{lemma}
\begin{proof}
Now $\bfG=\Sp_{2m(m+1)}$,
we know that
$R_\Sigma=2^{-m}\sum_{\Lambda\in\cals_Z}(-1)^{\langle\Sigma,\Lambda\rangle}\rho_\Lambda$
from (\ref{0315}).
Hence
\begin{align*}
\langle R_\Sigma,R_{\Sigma'}\rangle_\bfG
& =\frac{1}{2^{2m}}\sum_{\Lambda\in\cals_Z}\sum_{\Lambda'\in\cals_Z}
(-1)^{\langle\Sigma,\Lambda\rangle+\langle\Sigma',\Lambda'\rangle}\langle\rho_\Lambda,\rho_{\Lambda'}\rangle_\bfG 
=\frac{1}{2^{2m}}\sum_{\Lambda\in\cals_Z}(-1)^{\langle\Sigma,\Lambda\rangle+\langle\Sigma',\Lambda\rangle}.
\end{align*}
Now the set $\{\,R_\Sigma\mid\Sigma\in\cals_{Z,1}\,\}$ is orthonormal,
so the lemma follows.
\end{proof}

\begin{proposition}\label{0512}
Let $(\bfG,\bfG')=(\Sp_{2m(m+1)},\rmO^\epsilon_{2(m+1)^2})$
where $\epsilon=(-1)^{m+1}$.
Suppose that $Z=Z_{(m)}$ and $Z'=Z'_{(m+1)}$.
Then
\[
\omega^\sharp_{Z,Z'}
=\sum_{\Lambda\in\cals_Z}\rho_\Lambda^\sharp\otimes\rho^\sharp_{\theta^\epsilon(\Lambda)}.
\]
\end{proposition}
\begin{proof}
Now $\deg(Z)=m$ and $\deg(Z')=m+1$.
We know that $\{\,\theta(\Sigma')\mid\Sigma'\in\cals_{Z,1}\,\}$ for a complete set
of representatives of coset $\{\Sigma'',\Sigma''^\rmt\}$ in $\cals_{Z',0}$.
From (\ref{0315}) and (\ref{0327}),
for $\Lambda\in\cals_Z$, we know that
\[
\rho_\Lambda^\sharp
=\frac{1}{2^m}\sum_{\Sigma\in\cals_{Z,1}}(-1)^{\langle\Sigma,\Lambda\rangle}R_\Sigma,\quad\text{and}\quad
\rho_{\theta^\epsilon(\Lambda)}^\sharp
=\frac{1}{2^{m+1}}\sum_{\Sigma'\in\cals_{Z,1}}(-1)^{\langle\theta(\Sigma'),\theta^\epsilon(\Lambda)\rangle}R_{\theta(\Sigma')}.
\]
Therefore, we can write
\[
\sum_{\Lambda\in\cals_Z}\rho^\sharp_\Lambda\otimes\rho^\sharp_{\theta^\epsilon(\Lambda)}
=\sum_{\Sigma\in\cals_{Z,1}}\sum_{\Sigma'\in\cals_{Z,1}}c_{\Sigma,\Sigma'}R_{\Sigma}\otimes R_{\theta(\Sigma')}
\]
where
\[
c_{\Sigma,\Sigma'}
=\frac{1}{2^{2m+1}}\sum_{\Lambda\in\cals_Z}(-1)^{\langle\Sigma,\Lambda\rangle+\langle\theta(\Sigma'),\theta^\epsilon(\Lambda)\rangle}
=\frac{1}{2^{2m+1}}\sum_{\Lambda\in\cals_Z}(-1)^{\langle\Sigma,\Lambda\rangle+\langle\Sigma',\Lambda\rangle}
\]
by Lemma~\ref{0505}.
Hence by Lemma~\ref{0506},
we have
\[
c_{\Sigma,\Sigma'}=\begin{cases}
\frac{1}{2}, & \text{if $\Sigma'=\Sigma$};\\
0, & \text{otherwise}.
\end{cases}
\]
Then the proposition follows from Lemma~\ref{0507}.
\end{proof}

\subsection{Basic case II: from $\rmO^\epsilon_{2m^2}$ to $\Sp_{2m(m+1)}$}\label{0503}
In this subsection,
we assume that $\bfG'=\rmO_{2m^2}^\epsilon$ and $\bfG=\Sp_{2m(m+1)}$.
Let $Z'=Z'_{(m)}$ and $Z=Z_{(m)}$.
Then $Z'$ is a special symbol of rank $2m^2$ and defect $0$,
and $Z$ is a special symbol of rank $2m(m+1)$ and defect $1$.

Define a mapping
\begin{align}\label{0504}
\begin{split}
\theta=\theta_{Z',Z}\colon\cals_{Z',0} &\longrightarrow\cals_{Z,1}\quad\text{ by }
\Sigma'\mapsto\textstyle\binom{2m}{-}\cup\Sigma'^\rmt,\\
\theta^\epsilon\colon\cals^\epsilon_{Z'} &\longrightarrow\cals_Z\quad\text{ by }
\Lambda'\mapsto\begin{cases}
\binom{2m}{-}\cup\Lambda'^\rmt, & \text{if $\epsilon=+$};\\
\binom{-}{2m}\cup\Lambda'^\rmt, & \text{if $\epsilon=-$}.
\end{cases}
\end{split}
\end{align}
It is clear that the mapping $\theta_{Z,Z'}$ is injective and an element
$\Lambda$ is in the image of $\theta_{Z',Z}$ if and only the entry $2m$ is in the
first row of $\Lambda$.
Moreover, it is also clear that
\[
\theta^+(\cals^+_{Z'})\cup\theta^-(\cals^-_{Z'})=\cals_Z.
\]

\begin{lemma}\label{0401}
Suppose that $Z'=Z'_{(m)}$ and $Z=Z_{(m)}$.
Then
\[
\omega_{Z,Z'}^\sharp
=\frac{1}{2}\sum_{\Sigma'\in\cals_{Z',0}}R_{\theta(\Sigma')}\otimes R_{\Sigma'}.
\]
\end{lemma}
\begin{proof}
Similar to the proof of Lemma~\ref{0507},
we need to prove that
\[
\cald_{Z,Z'}=\{\,(\theta(\Sigma'),\Sigma')\mid\Sigma'\in\cals_{Z',0}\,\},
\]
which is just Example~\ref{0226}.
\end{proof}

\begin{lemma}\label{0511}
Let $N$ be a subset of $Z'_\rmI$.
Then
\[
\theta^\epsilon(\Lambda_N)=\begin{cases}
\Lambda_{N^\rmt}, & \text{if $\epsilon=+$};\\
\Lambda_{\binom{-}{2m}\cup N^\rmt}, & \text{if $\epsilon=-$}.
\end{cases}
\]
\end{lemma}
\begin{proof}
The proof is similar to that of Lemma~\ref{0510}.
\end{proof}

\begin{lemma}\label{0508}
We have $\langle\Sigma',\Lambda'\rangle=\langle\theta(\Sigma'),\theta^\epsilon(\Lambda')\rangle$
for each $\Lambda'\in\cals^\epsilon_{Z'}$ and $\Sigma'\in\cals_{Z',0}$.
\end{lemma}
\begin{proof}
The proof is similar to that of Lemma~\ref{0505}.
\end{proof}

\begin{lemma}\label{0402}
Let $\Sigma,\Sigma'\in\cals_{Z,1}$.
Suppose that $\Sigma'\in\theta(\cals_{Z',0})$ and $\Sigma\not\in\theta(\cals_{Z',0})$.
Then
\[
\sum_{\Lambda'\in\cals^\epsilon_{Z'}}(-1)^{\langle\Sigma,\theta^\epsilon(\Lambda')\rangle+\langle\Sigma',\theta^\epsilon(\Lambda')\rangle}
=0.
\]
\end{lemma}
\begin{proof}
By Lemma~\ref{0511}, $\Sigma'\in\theta(\cals_{Z',0})$ means that $\Sigma'=\Lambda_N$ for some $N\subset Z'_\rmI$ such that
$|N^*|=|N_*|$, in particular, $|N|$ is even.
Moreover, $\Sigma\not\in\theta(\cals_{Z',0})$ means that $\Sigma=\Lambda_{\binom{2m}{-}\cup M}$
for some subset $M\subset Z'_\rmI$ such that $|M^*|+1=|M_*|$, in particular, $|M|$ is odd.
Suppose that $\Lambda'=\Lambda_K$ for some subset $K\subset Z'_\rmI$.
Then by Lemma~\ref{0511}, we know that
$\langle\Sigma,\theta^\epsilon(\Lambda')\rangle=|M\cap K|\pmod 2$ and $\langle\Sigma',\theta^\epsilon(\Lambda')\rangle=|N\cap K|\pmod 2$.
Now for any subset $K\subset Z'_\rmI$, it is clearly that
$|M\cap K|\equiv|N\cap K|\pmod 2$ if and only if $|M\cap(Z'_\rmI\smallsetminus K)|\not\equiv|N\cap(Z'_\rmI\smallsetminus K)|\pmod 2$.
Moreover, for $K\subset Z'_\rmI$,
we know that $\Lambda_K\in\cals^\epsilon_{Z'}$ if and only if $\Lambda_{Z'_\rmI\smallsetminus K}=(\Lambda_K)^\rmt\in\cals^\epsilon_{Z'}$.
Then the lemma is proved.
\end{proof}

\begin{lemma}\label{0509}
Let $\Sigma,\Sigma'\in\cals_{Z',0}$.
Then
\[
\sum_{\Lambda\in\cals^\epsilon_{Z'}}(-1)^{\langle\Sigma,\Lambda\rangle+\langle\Sigma',\Lambda\rangle}
=\begin{cases}
2^{2m-1}, & \text{if\/ $\Sigma=\Sigma'$};\\
2^{2m-1}\epsilon, & \text{if\/ $\Sigma=\Sigma'^\rmt$};\\
0, & \text{otherwise}.
\end{cases}
\]
\end{lemma}
\begin{proof}
Recall that $R_\Sigma=2^{-(m-1)}\sum_{\Lambda\in\cals^\epsilon_{Z'}}(-1)^{\langle\Sigma,\Lambda\rangle}\rho_\Lambda$
from (\ref{0316}).
Hence
\begin{align*}
\langle R_\Sigma,R_{\Sigma'}\rangle_{\bfG'}
& =2^{-2(m-1)}\sum_{\Lambda\in\cals^\epsilon_{Z'}}\sum_{\Lambda'\in\cals^\epsilon_{Z'}}
(-1)^{\langle\Sigma,\Lambda\rangle+\langle\Sigma',\Lambda'\rangle}\langle\rho_\Lambda,\rho_{\Lambda'}\rangle_{\bfG'} \\
&=2^{-2(m-1)}\sum_{\Lambda\in\cals^\epsilon_{Z'}}(-1)^{\langle\Sigma,\Lambda\rangle+\langle\Sigma',\Lambda\rangle}.
\end{align*}
Note that $\Sigma$ is always non-degenerate for our $Z$ right now.
Then we know that
\[
\langle R_\Sigma,R_{\Sigma'}\rangle_{\bfG'}=
\begin{cases}
2, & \text{if $\Sigma=\Sigma'$};\\
2\epsilon, & \text{if $\Sigma=\Sigma'^\rmt$};\\
0, & \text{if $\Sigma\neq\Sigma',\Sigma'^\rmt$}.
\end{cases}
\]
from (\ref{0317}) and the fact $R^{\rmO^\epsilon}_{\Sigma^\rmt}=\epsilon R^{\rmO^\epsilon}_\Sigma$,
so the lemma follows.
\end{proof}

\begin{proposition}\label{0513}
Let $(\bfG,\bfG')=(\Sp_{2m(m+1)},\rmO^\epsilon_{2m^2})$ where $\epsilon=(-1)^m$.
Suppose that $Z=Z_{(m)}$ and $Z'=Z'_{(m)}$ for some positive integer $m$.
Then
\[
\omega_{Z,Z'}^\sharp
=\sum_{\Lambda'\in\cals^\epsilon_{Z'}}\rho^\sharp_{\theta^\epsilon(\Lambda')}\otimes\rho^\sharp_{\Lambda'}.
\]
\end{proposition}
\begin{proof}
Now $\deg(Z)=\deg(Z')=m$.
From (\ref{0316}) and (\ref{0327}),
for $\Lambda'\in\cals^\epsilon_{Z'}$, we know that
\[
\rho_{\theta^\epsilon(\Lambda')}^\sharp
=\frac{1}{2^m}\sum_{\Sigma\in\cals_{Z,1}}(-1)^{\langle\Sigma,\theta^\epsilon(\Lambda')\rangle}R_\Sigma\quad\text{and}\quad
\rho_{\Lambda'}^\sharp
=\frac{1}{2^{m+1}}\sum_{\Sigma'\in\cals_{Z',0}}(-1)^{\langle\Sigma',\Lambda'\rangle}R_{\Sigma'}.
\]
Therefore, we can write
\[
\sum_{\Lambda'\in\cals^\epsilon_{Z'}}\rho^\sharp_{\theta^\epsilon(\Lambda')}\otimes\rho^\sharp_{\Lambda'}
=\sum_{\Sigma\in\cals_{Z,1}}\sum_{\Sigma'\in\cals_{Z',0}}c_{\Sigma,\Sigma'}R_\Sigma\otimes R_{\Sigma'}
\]
where
\begin{align*}
c_{\Sigma,\Sigma'}
&=\frac{1}{2^{2m+1}}\sum_{\Lambda'\in\cals^\epsilon_{Z'}}
(-1)^{\langle\Sigma,\theta^\epsilon(\Lambda')\rangle+\langle\Sigma',\Lambda'\rangle} \\
&=\frac{1}{2^{2m+1}}\sum_{\Lambda'\in\cals^\epsilon_{Z'}}
(-1)^{\langle\Sigma,\theta^\epsilon(\Lambda')\rangle+\langle\theta(\Sigma'),\theta^\epsilon(\Lambda')\rangle}
\end{align*}
by Lemma~\ref{0508}.
By Lemma~\ref{0402}, we see that $c_{\Sigma,\Sigma'}=0$ if $\Sigma\not\in\theta(\cals_{Z',0})$.
Now suppose that $\Sigma=\theta(\Sigma'')$ for some $\Sigma''\in\cals_{Z',0}$.
Then by Lemma~\ref{0508} and Lemma~\ref{0509},
we have
\begin{align*}
c_{\theta(\Sigma''),\Sigma'}
=\frac{1}{2^{2m+1}}\sum_{\Lambda'\in\cals^\epsilon_{Z'}}(-1)^{\langle\theta(\Sigma''),\theta^\epsilon(\Lambda')\rangle+\langle\Sigma',\Lambda'\rangle}
& =\frac{1}{2^{2m+1}}\sum_{\Lambda'\in\cals^\epsilon_{Z'}}(-1)^{\langle\Sigma'',\Lambda'\rangle+\langle\Sigma',\Lambda'\rangle} \\
&=\begin{cases}
\frac{1}{4}, & \text{if $\Sigma''=\Sigma'$};\\
\frac{1}{4}\epsilon, & \text{if $\Sigma''=\Sigma'^\rmt$};\\
0, & \text{otherwise}.
\end{cases}
\end{align*}
Therefore, we have
\begin{align*}
\sum_{\Lambda'\in\cals^\epsilon_{Z'}}\rho^\sharp_{\theta^\epsilon(\Lambda')}\otimes\rho^\sharp_{\Lambda'}
& =\frac{1}{4}\sum_{\Sigma'\in\cals_{Z',0}}R_{\theta(\Sigma')}\otimes R_{\Sigma'}
 +\epsilon\frac{1}{4}\sum_{\Sigma'\in\cals_{Z',0}}R_{\theta(\Sigma')}\otimes R_{\Sigma'^\rmt} \\
& =\frac{1}{2}\sum_{\Sigma'\in\cals_{Z',0}}R_{\theta(\Sigma')}\otimes R_{\Sigma'}
\end{align*}
since $R_{\Sigma'^\rmt}=\epsilon R_{\Sigma'}$.
Then the proposition follows from Lemma~\ref{0401} immediately.
\end{proof}


\section{The Relation $\cald_{Z,Z'}$}\label{0514}
In this section, we fix special symbols $Z,Z'$ of sizes $(m+1,m),(m',m')$ respectively such that $m'=m,m+1$.
The material in this section is really technical and our purpose is Proposition~\ref{0216}.

\subsection{Consecutive pairs}\label{0609}
Let $Z$ be a special symbol of degree $\delta$.
Write
\[
Z_\rmI=\binom{s_1,s_2,\ldots,s_{\delta'}}{t_1,t_2,\ldots,t_\delta}
\]
where $\delta'=\delta,\delta+1$ depending on ${\rm def}(Z)=0,1$ respectively.
A pair $\binom{s}{t}$ in $Z_\rmI$ of the form $\binom{s_i}{t_i}$ or $\binom{s_i}{t_{i+1}}$ is called \emph{semi-consecutive}.
A semi-consecutive pair $\binom{s}{t}$ in $Z_\rmI$ is called \emph{consecutive} if there is no other entry in $Z$ lying between $s$ and $t$ i.e., there is no entry $x$ in $Z$ such that $s<x<t$ or $t<x<s$.
Two consecutive pairs $\binom{s}{t},\binom{s'}{t'}$ in $Z_\rmI$ are called \emph{disjoint} if
$\{s,t\}\cap\{s',t'\}=\emptyset$.

\begin{example}
Suppose that $Z=\binom{7,6,4,3}{6,5,4,1}$, a special symbol of defect $0$ and degree $2$.
Now $Z_\rmI=\binom{7,3}{5,1}$.
Then the pair $\binom{3}{1}$ in $Z_\rmI$ is consecutive,
and $\binom{7}{5},\binom{3}{5}$ are only semi-consecutive pairs which are not consecutive.
\end{example}

For a subset of consecutive pairs $\Psi_0$ in $Z_\rmI$,
we define several subsets of $\cals_Z$ or $\overline\cals_Z$ as follows.
\begin{enumerate}
\item[(1)] Suppose $Z$ is of defect $1$ and degree $\delta$.
Define
\begin{align*}
\cals_{Z,\Psi_0} &=\{\,\Lambda_M\in\cals_Z\mid M\leq\Psi_0\,\}\subset\cals_{Z,1};\\
\overline\cals_Z^{\Psi_0} &=\{\,\Lambda_M\in\overline\cals_Z\mid M\subset Z_\rmI\smallsetminus\Psi_0\,\};\\
\cals_Z^{\Psi_0} &=\overline\cals_Z^{\Psi_0}\cap\cals_Z.
\end{align*}
Clearly, $\cals_{Z,\Psi_0}=\{Z\}$ and $\cals_Z^{\Psi_0}=\cals_Z$ if $\Psi_0=\emptyset$.
If $\delta_0$ is the number of pairs in $\Psi_0$,
then it is clear that $|\cals_{Z,\Psi_0}|=2^{\delta_0}$ and $|\cals_Z^{\Psi_0}|=2^{2(\delta-\delta_0)}$.

\item[(2)] Suppose $Z$ is of defect $0$.
Define
\begin{align*}
\cals_{Z,\Psi_0} &=\{\,\Lambda_M\in\cals^+_Z\mid M\leq\Psi_0\,\}\subset\cals_{Z,0};\\
\overline\cals_Z^{\Psi_0} &=\{\,\Lambda_M\in\overline\cals_Z\mid M\subset Z_\rmI\smallsetminus\Psi_0\,\};\\
\cals_Z^{\epsilon,\Psi_0} &=\overline\cals_Z^{\Psi_0}\cap\cals^\epsilon_Z.
\end{align*}
\end{enumerate}
Note that if ${\rm def}(Z)$ is not specified,
we will just use the notation $\cals_Z^{\Psi_0}$ to denote $\cals_Z^{\Psi_0}$ (when ${\rm def}(Z)=1$) or
$\cals_Z^{+,\Psi_0},\cals_Z^{-,\Psi_0}$ (when ${\rm def}(Z)=0$).

\subsection{The set $\calb^+_{Z,Z'}$ I}
Let
\begin{equation}\label{0204}
\Lambda_M=\binom{a_1,a_2,\ldots,a_{m_1}}{b_1,b_2,\ldots,b_{m_2}}\in\overline\cals_Z,\qquad
\Lambda_N=\binom{c_1,c_2,\ldots,c_{m'_1}}{d_1,d_2,\ldots,d_{m'_2}}\in\overline\cals_{Z'}
\end{equation}
for some $M\subset Z_\rmI$ and $N\subset Z'_\rmI$ such that $m_1'=m_2$ if $m'=m$; and $m_2'=m_1$ if $m'=m+1$.

\begin{lemma}\label{0615}
Let $\Lambda_M$ be given as in (\ref{0204}) and suppose that $M\neq\emptyset$.
Let $x=\max(M)$ be the largest element in $M$.
\begin{enumerate}
\item[(i)] If $x=a_k$,
then $k\geq 2$ and $b_{k-2}>a_k>b_{k-1}$.

\item[(ii)] If $x=b_k$,
then $a_{k-1}>b_k>a_k$.
\end{enumerate}
\end{lemma}
\begin{proof}
Recall that $\Lambda_M=(Z\smallsetminus M)\cup M^\rmt$ and write
$Z=\binom{s_1,\ldots,s_{m+1}}{t_1,\ldots,t_m}$.

First suppose that $x=a_k$.
Then the natural position of $x$ is in the second row of $Z$,
so $s_1\not\in M$, and hence $s_1$ is still in the first row of $\Lambda_M$.
Therefore $k\neq 1$.
Suppose that $b_{k'-1}>a_k>b_{k'}$ for some $k'$.
Now $a_1,\ldots,a_{k-1},b_1,\ldots,b_{k'-1}$ are all the entries in $\Lambda_M$ which are greater than $a_k$,
and hence the following entries
\[
\binom{a_1,a_2,\ldots,a_{k-2},a_{k-1}}{b_1,b_2,\ldots,b_{k'-1},a_k}
\]
are all in their natural positions.
This implies that $a_{k-1}=s_{k-1}$, $a_k=t_{k-1}$, and hence $b_{k'-1}=t_{k-2}$, i.e., $k'=k-1$.

The proof for (ii) is similar.
\end{proof}

\begin{lemma}\label{0214}
Keep the notations of (\ref{0204}).
Suppose that $(\Lambda_M,\Lambda_N)\in\overline\calb^+_{Z,Z'}$ and $\max(M)=a_k$.
\begin{enumerate}
\item[(i)] If either
\begin{enumerate}
\item $m'=m$ and $a_k>c_{k-1}$, or

\item $m'=m+1$ and $a_k\geq c_{k-1}$,
\end{enumerate}
then $c_{k-2}>d_{k-1}>c_{k-1}$, in particular, $d_{k-1}\in Z'_\rmI$.

\item[(ii)] If either
\begin{enumerate}
\item $m'=m$, $c_{k-1}\geq a_k$ and $d_k\geq b_{k-1}$, or

\item $m'=m+1$, $c_{k-1}>a_k$ and $d_k>b_{k-1}$,
\end{enumerate}
then $c_{k-1}>d_k>c_k$, in particular, $d_k\in Z'_\rmI$.

\item[(iii)] If either
\begin{enumerate}
\item $m'=m$, $c_{k-1}\geq a_k$ and $b_{k-1}>d_k$, or

\item $m'=m+1$, $c_{k-1}>a_k$ and $b_{k-1}\geq d_k$,
\end{enumerate}
then $a_k>b_{k-1}>a_{k+1}$, in particular, $b_{k-1}\in Z_\rmI$.
\end{enumerate}
\end{lemma}
\begin{proof}
Now $\max(M)=a_k$,
from the proof of the previous lemma
we know that the entries $a_1,\ldots,a_{k-1}$ and $b_1,\ldots,b_{k-2}$ in $\Lambda_M$
are in their natural positions.
In particular, we have $b_{k-2}\geq a_{k-1}$.

Suppose that $a_k>c_{k-1}$ if $m'=m$; or suppose that $a_k\geq c_{k-1}$ if $m'=m+1$.
Then by Lemma~\ref{0210},
we have
\[
\begin{cases}
c_{k-2}\geq b_{k-2}\geq a_{k-1}>d_{k-1}\geq a_k>c_{k-1}, & \text{if $m'=m$};\\
c_{k-2}> b_{k-2}\geq a_{k-1}\geq d_{k-1}> a_k\geq c_{k-1}, & \text{if $m'=m+1$}.
\end{cases}
\]
We have $c_{k-2}>d_{k-1}>c_{k-1}$ for both cases,
and hence $d_{k-1}\in Z_\rmI'$.

The proofs for other cases are similar.
\end{proof}

\begin{lemma}\label{0215}
Keep the notations of (\ref{0204}).
Suppose that $(\Lambda_M,\Lambda_N)\in\overline\calb^+_{Z,Z'}$ and $\max(M)=b_k$.
\begin{enumerate}
\item[(i)] If either
\begin{enumerate}
\item $m'=m$ and $b_k>d_{k-1}$, or

\item $m'=m+1$ and $b_k\geq d_{k-1}$,
\end{enumerate}
then $d_{k-2}>c_k>d_{k-1}$, in particular, $c_k\in Z'_\rmI$.

\item[(ii)] If either
\begin{enumerate}
\item $m'=m$, $d_{k-1}\geq b_k$ and $c_{k+1}\geq a_k$, or

\item $m'=m+1$, $d_{k-1}>b_k$ and $c_{k+1}>a_k$,
\end{enumerate}
then $d_{k-1}>c_{k+1}>d_k$, in particular, $c_{k+1}\in Z'_\rmI$.

\item[(iii)] If either
\begin{enumerate}
\item $m'=m$, $d_{k-1}\geq b_k$ and $a_k>c_{k+1}$, or

\item $m'=m+1$, $d_{k-1}>b_k$ and $a_k\geq c_{k+1}$,
\end{enumerate}
then $b_k>a_k>b_{k+1}$, in particular, $a_k\in Z_\rmI$.
\end{enumerate}
\end{lemma}
\begin{proof}
Now $\max(M)=b_k$, by Lemma~\ref{0615},
the entries $a_1,\ldots,a_{k-1}$ and $b_1,\ldots,b_{k-1}$ in $\Lambda_M$ are in their natural positions.
Then we have $a_{k-1}\geq b_{k-1}$.

Suppose that $b_k>d_{k-1}$ if $m'=m$; or suppose that $b_k\geq d_{k-1}$ if $m'=m+1$.
Then by Lemma~\ref{0210}, we have
\[
\begin{cases}
d_{k-2}\geq a_{k-1}\geq b_{k-1}>c_k\geq b_k>d_{k-1}, & \text{if $m'=m$};\\
d_{k-2}>a_{k-1}\geq b_{k-1}\geq c_k>b_k\geq d_{k-1}, & \text{if $m'=m+1$}.
\end{cases}
\]
We have $d_{k-2}>c_k>d_{k-1}$ for both cases, in particular, $c_k\in Z'_\rmI$.

The proofs for other cases are similar.
\end{proof}

\begin{lemma}
Keep the notations of (\ref{0204}).
Suppose that $(\Lambda_M,\Lambda_N)\in\overline\calb^+_{Z,Z'}$ and $M\neq\emptyset$.
\begin{enumerate}
\item[(i)] If\/ $\max(M)=a_k$,
then $b_{k-1}$ or $d_k$ exist.

\item[(ii)]  If\/ $\max(M)=b_k$ and $d_{k-1}$ exists,
then $a_k$ or $c_{k+1}$ exist.
\end{enumerate}
\end{lemma}
\begin{proof}
First suppose that $\max(M)=a_k$ and $b_{k-1}$ does not exist.
From Lemma~\ref{0615}, we know that either $k=2$ or $b_{k-2}$ exists.
Because ${\rm def}(Z)=1$,
$\Lambda_M$ has odd number of entries and hence $a_{k+1}$ exists.
By Lemma~\ref{0210} we need $d_k\geq a_{k+1}$ (resp.~$d_k>a_{k+1}$) if $m'=m$ (resp.~$m'=m+1$),
so $d_k$ must exist.

The proof for (ii) is similar.
\end{proof}

\subsection{The set $\calb^+_{Z,Z'}$ II}
Keep the notations in (\ref{0204}) and suppose that $(\Lambda_M,\Lambda_N)\in\overline\calb^+_{Z,Z'}$
and $M\neq\emptyset$.
Now we want to define a new symbol
$\Lambda_{M'}=\binom{a'_1,\ldots,a'_{m_3}}{b'_1,\ldots,b'_{m_4}}\in\overline\cals_Z$ by moving (at least) the maximal element in
$M$ back to its natural position and also define a new symbol $\Lambda_{N'}=\binom{c'_1,\ldots,c'_{m'_3}}{d'_1,\ldots,d'_{m'_4}}\in\overline\cals_{Z'}$
such that $(\Lambda_{M'},\Lambda_{N'})\in\overline\calb^+_{Z,Z'}$.

\begin{enumerate}
\item[(1)] Suppose that $\max(M)=a_k$.
Then $a_1,\ldots,a_{k-1}$ and $b_1,\ldots,b_{k-2}$ in $\Lambda_M$ are in their natural positions
by Lemma~\ref{0615}.

\begin{enumerate}
\item
Suppose that
\[
\begin{cases}
a_k>c_{k-1}\text{ (or $c_{k-1}$ does not exist)}, & \text{if $m'=m$};\\
a_k\geq c_{k-1}\text{ (or $c_{k-1}$ does not exist)}, & \text{if $m'=m+1$}.
\end{cases}
\]
By Lemma~\ref{0214}, we know that $d_{k-1}\in Z'_\rmI$.
Let $\Lambda_{M'}$ be obtained from $\Lambda_M$ by moving $a_k$ to the second row,
i.e., $M'=M\smallsetminus\{a_k\}$,
and let $\Lambda_{N'}$ be obtained from $\Lambda_N$ by moving $d_{k-1}$ to the first row.
By Lemma~\ref{0615} and Lemma~\ref{0214}, we have the following table for $\Lambda_{M'}$ and $\Lambda_{N'}$:
\[
\begin{tabular}{c|c|c|c|c|c|c|c}
\toprule
$i$ & $1$ & $\ldots$ & $k-2$ & $k-1$ & $k$ & $k+1$ & $\ldots$ \\
\midrule
$a'_i$ & $a_1$ & $\ldots$ & $a_{k-2}$ & $a_{k-1}$ & $a_{k+1}$ & $a_{k+2}$ & $\ldots$ \\
\midrule
$b'_i$ & $b_1$ & $\ldots$ & $b_{k-2}$ & $a_k$ & $b_{k-1}$ & $b_k$ & $\ldots$ \\
\midrule
$c'_i$ & $c_1$ & $\ldots$ & $c_{k-2}$ & $d_{k-1}$ & $c_{k-1}$ & $c_k$ & $\ldots$ \\
\midrule
$d'_i$ & $d_1$ & $\ldots$ & $d_{k-2}$ & $d_k$ & $d_{k+1}$ & $d_{k+2}$ & $\ldots$ \\
\bottomrule
\end{tabular}
\]

\item Suppose that
\[
\begin{cases}
c_{k-1}\geq a_k\text{ and }d_k\geq b_{k-1}\text{ (or $b_{k-1}$ does not exist)}, & \text{if $m'=m$};\\
c_{k-1}>a_k\text{ and }d_k>b_{k-1}\text{ (or $b_{k-1}$ does not exist)}, & \text{if $m'=m+1$}.
\end{cases}
\]
Then by Lemma~\ref{0214}, we know that $d_k\in Z'_\rmI$.
Let $\Lambda_{M'}$ be given as in (a);
and let $\Lambda_{N'}$ be obtained from $\Lambda_N$ by moving $d_k$ to the first row.
By Lemma~\ref{0615} and Lemma~\ref{0214}, we have the following table for $\Lambda_{M'}$ and $\Lambda_{N'}$:
\[
\begin{tabular}{c|c|c|c|c|c|c|c}
\toprule
$i$ & $1$ & $\ldots$ & $k-2$ & $k-1$ & $k$ & $k+1$ & $\ldots$ \\
\midrule
$a'_i$ & $a_1$ & $\ldots$ & $a_{k-2}$ & $a_{k-1}$ & $a_{k+1}$ & $a_{k+2}$ & $\ldots$ \\
\midrule
$b'_i$ & $b_1$ & $\ldots$ & $b_{k-2}$ & $a_k$ & $b_{k-1}$ & $b_k$ & $\ldots$ \\
\midrule
$c'_i$ & $c_1$ & $\ldots$ & $c_{k-2}$ & $c_{k-1}$ & $d_k$ & $c_k$ & $\ldots$ \\
\midrule
$d'_i$ & $d_1$ & $\ldots$ & $d_{k-2}$ & $d_{k-1}$ & $d_{k+1}$ & $d_{k+2}$ & $\ldots$ \\
\bottomrule
\end{tabular}
\]

\item
Suppose that
\[
\begin{cases}
c_{k-1}\geq a_k\text{ and }b_{k-1}>d_k\text{ (or $d_k$ does not exist)}, & \text{if $m'=m$};\\
c_{k-1}>a_k\text{ and }b_{k-1}\geq d_k\text{ (or $d_k$ does not exist)}, & \text{if $m'=m+1$}.
\end{cases}
\]
Then by Lemma~\ref{0214}, we know that $b_{k-1}\in Z_\rmI$.
Let $\Lambda_{M'}$ be obtained from $\Lambda_M$ by moving $a_k$ to the second row and
$b_{k-1}$ to the first row,
and let $\Lambda_{N'}=\Lambda_N$.
By Lemma~\ref{0615} and Lemma~\ref{0214}, we have the following table:
\[
\begin{tabular}{c|c|c|c|c|c|c|c}
\toprule
$i$ & $1$ & $\ldots$ & $k-2$ & $k-1$ & $k$ & $k+1$ & $\ldots$ \\
\midrule
$a'_i$ & $a_1$ & $\ldots$ & $a_{k-2}$ & $a_{k-1}$ & $b_{k-1}$ & $a_{k+1}$ & $\ldots$ \\
\midrule
$b'_i$ & $b_1$ & $\ldots$ & $b_{k-2}$ & $a_k$ & $b_k$ & $b_{k+1}$ & $\ldots$ \\
\midrule
$c'_i$ & $c_1$ & $\ldots$ & $c_{k-2}$ & $c_{k-1}$ & $c_k$ & $c_{k+1}$ & $\ldots$ \\
\midrule
$d'_i$ & $d_1$ & $\ldots$ & $d_{k-2}$ & $d_{k-1}$ & $d_k$ & $d_{k+1}$ & $\ldots$ \\
\bottomrule
\end{tabular}
\]
\end{enumerate}

\item[(2)] Suppose that $\max(M)=b_k$.
Then $a_1,\ldots,a_{k-1}$ and $b_1,\ldots,b_{k-1}$ in $\Lambda_M$ are in their natural positions
by Lemma~\ref{0615}.

\begin{enumerate}
\item[(d)]
Suppose that
\[
\begin{cases}
b_k>d_{k-1}\text{ (or $d_{k-1}$ does not exist)}, & \text{if $m'=m$};\\
b_k\geq d_{k-1}\text{ (or $d_{k-1}$ does not exist)}, & \text{if $m'=m+1$}.
\end{cases}
\]
Then by Lemma~\ref{0215}, we know that $c_k\in Z'_\rmI$.
Let $\Lambda_{M'}$ be obtained from $\Lambda_M$ by moving $b_k$ to the first row,
i.e., $M'=M\smallsetminus\{b_k\}$,
and let $\Lambda_{N'}$ be obtained from $\Lambda_N$ by moving $c_k$ to the second row.
By Lemma~\ref{0615} and Lemma~\ref{0215} we have the following table for $\Lambda_{M'}$ and $\Lambda_{N'}$:
\[
\begin{tabular}{c|c|c|c|c|c|c|c}
\toprule
$i$ & $1$ & $\ldots$ & $k-2$ & $k-1$ & $k$ & $k+1$ & $\ldots$ \\
\midrule
$a'_i$ & $a_1$ & $\ldots$ & $a_{k-2}$ & $a_{k-1}$ & $b_k$ & $a_k$ & $\ldots$ \\
\midrule
$b'_i$ & $b_1$ & $\ldots$ & $b_{k-2}$ & $b_{k-1}$ & $b_{k+1}$ & $b_{k+2}$ & $\ldots$ \\
\midrule
$c'_i$ & $c_1$ & $\ldots$ & $c_{k-2}$ & $c_{k-1}$ & $c_{k+1}$ & $c_{k+2}$ & $\ldots$ \\
\midrule
$d'_i$ & $d_1$ & $\ldots$ & $d_{k-2}$ & $c_k$ & $d_{k-1}$ & $d_k$ & $\ldots$ \\
\bottomrule
\end{tabular}
\]

\item[(e)]
Suppose that
\[
\begin{cases}
d_{k-1}\geq b_k\text{ and }c_{k+1}\geq a_k\text{ (or $a_k$ does not exist)}, & \text{if $m'=m$};\\
d_{k-1}>b_k\text{ and }c_{k+1}>a_k\text{ (or $a_k$ does not exist)}, & \text{if $m'=m+1$}.
\end{cases}
\]
Then by Lemma~\ref{0215}, we know that $c_{k+1}\in Z'_\rmI$.
Let $\Lambda_{M'}$ be as in (d),
and let $\Lambda_{N'}$ be obtained from $\Lambda_N$ by moving $c_{k+1}$ to the second row.
By Lemma~\ref{0615} and Lemma~\ref{0215} we have the following table for $\Lambda_{M'}$ and $\Lambda_{N'}$:
\[
\begin{tabular}{c|c|c|c|c|c|c|c}
\toprule
$i$ & $1$ & $\ldots$ & $k-2$ & $k-1$ & $k$ & $k+1$ & $\ldots$ \\
\midrule
$a'_i$ & $a_1$ & $\ldots$ & $a_{k-2}$ & $a_{k-1}$ & $b_k$ & $a_k$ & $\ldots$ \\
\midrule
$b'_i$ & $b_1$ & $\ldots$ & $b_{k-2}$ & $b_{k-1}$ & $b_{k+1}$ & $b_{k+2}$ & $\ldots$ \\
\midrule
$c'_i$ & $c_1$ & $\ldots$ & $c_{k-2}$ & $c_{k-1}$ & $c_k$ & $c_{k+2}$ & $\ldots$ \\
\midrule
$d'_i$ & $d_1$ & $\ldots$ & $d_{k-2}$ & $d_{k-1}$ & $c_{k+1}$ & $d_k$ & $\ldots$ \\
\bottomrule
\end{tabular}
\]

\item[(f)]
Suppose that
\[
\begin{cases}
d_{k-1}\geq b_k\text{ and }a_k>c_{k+1}\text{ (or $c_{k+1}$ does not exist)}, & \text{if $m'=m$};\\
d_{k-1}>b_k\text{ and }a_k\geq c_{k+1}\text{ (or $c_{k+1}$ does not exist)}, & \text{if $m'=m+1$}.
\end{cases}
\]
Then by Lemma~\ref{0215}, we know that $a_k\in Z_\rmI$.
Let $\Lambda_{M'}$ be obtained from $\Lambda_M$ by moving $b_k$ to the first row
and $a_k$ to the second row,
and let $\Lambda_{N'}=\Lambda_N$.
By Lemma~\ref{0615} and Lemma~\ref{0215} we have the following table for $\Lambda_{M'}$ and $\Lambda_{N'}$:
\[
\begin{tabular}{c|c|c|c|c|c|c|c}
\toprule
$i$ & $1$ & $\ldots$ & $k-2$ & $k-1$ & $k$ & $k+1$ & $\ldots$ \\
\midrule
$a'_i$ & $a_1$ & $\ldots$ & $a_{k-2}$ & $a_{k-1}$ & $b_k$ & $a_{k+1}$ & $\ldots$ \\
\midrule
$b'_i$ & $b_1$ & $\ldots$ & $b_{k-2}$ & $b_{k-1}$ & $a_k$ & $b_{k+1}$ & $\ldots$ \\
\midrule
$c'_i$ & $c_1$ & $\ldots$ & $c_{k-2}$ & $c_{k-1}$ & $c_k$ & $c_{k+1}$ & $\ldots$ \\
\midrule
$d'_i$ & $d_1$ & $\ldots$ & $d_{k-2}$ & $d_{k-1}$ & $d_k$ & $d_{k+1}$ & $\ldots$ \\
\bottomrule
\end{tabular}
\]
\end{enumerate}
\end{enumerate}

\begin{example}
Suppose that
\[
\Lambda_M=\binom{8,6,5}{3,1}=\binom{a_1,a_2,a_3}{b_1,b_2},\quad
\Lambda_N=\binom{6,3,0}{8,6,2}=\binom{c_1,c_2,c_3}{d_1,d_2,d_3}.
\]
Now $m=2$ and $m'=3$, $Z=\binom{8,5,1}{6,3}$, $M=\binom{1}{6}$,
$Z'=\binom{8,6,2}{6,3,0}$, $N=\binom{8,2}{3,0}$, and $(\Lambda_M,\Lambda_N)\in\overline\calb^+_{Z,Z'}$.
Then $\max(M)=a_2=6$, $k=2$, and $a_2\geq c_1=6$.
So we are in Case (a).
Then $a_2$ in $\Lambda_M$ is moved the second row, $d_1$ in $\Lambda_N$ is moved to the first row,
and we get
\[
\Lambda_{M'}=\binom{8,5}{6,3,1}=\binom{a'_1,a'_2}{b'_1,b'_2,b'_3},\quad
\Lambda_{N'}=\binom{8,6,3,0}{6,2}=\binom{c'_1,c'_2,c'_3,c'_4}{d'_1,d'_2}.
\]
For the next stage,
now $M'=\binom{1}{-}$ with $\max(M')=b'_3=1$ and $k=3$.
Now $b'_3<d'_2=2$ and $a'_3$ does not exist,
so we are in Case (e) and then $b'_3$ in $\Lambda_{M'}$ is moved to the first row,
and $c'_4$ in $\Lambda_{N'}$ is moved to the second row:
\[
\Lambda_{M''}=Z=\binom{8,5,1}{6,3},\quad
\Lambda_{N''}=\binom{8,6,3}{6,2,0}.
\]
It is easy to see that both $(\Lambda_{M'},\Lambda_{N'})$ and $(\Lambda_{M''},\Lambda_{N''})$
occur in $\overline\calb^+_{Z,Z'}$ by Lemma~\ref{0210}.
\end{example}

\begin{lemma}\label{0617}
For all Cases (a)--(f) above,
if $M'\neq\emptyset$, then $\max(M')<\max(M)$.
\end{lemma}
\begin{proof}
For Cases (a) and (b), we have $\max(M)=a_k$ and $M'=M\smallsetminus\{a_k\}$,
so the lemma is clearly true.

For Case (c), we have $\max(M)=a_k$ and
$M'=M\smallsetminus\{a_k,b_{k-1}\}$ if $b_{k-1}\in M$;
$M'=(M\smallsetminus\{a_k\})\cup\{b_{k-1}\}$ if $b_{k-1}\not\in M$.
Since $a_k>b_{k-1}$ by (i) of Lemma~\ref{0615}, the lemma is true, again.

For Cases (d) and (e), we have $\max(M)=b_k$ and $M'=M\smallsetminus\{b_k\}$,
so the lemma is true.

For Case (f), we have $\max(M)=b_k$ and $M'=M\smallsetminus\{a_k,b_k\}$ if $a_k\in M$;
$M'=(M\smallsetminus\{b_k\})\cup\{a_k\}$ if $a_k\not\in M$.
Since now $b_k>a_k$ by (ii) of Lemma~\ref{0615}, the lemma is true, again.
\qed
\end{proof}

\begin{lemma}\label{0205}
Keep notation in (\ref{0204}) and suppose that $(\Lambda_M,\Lambda_N)\in\overline\calb^+_{Z,Z'}$.
Let $\Lambda_{M'},\Lambda_{N'}$ be defined as above.
Then $(\Lambda_{M'},\Lambda_{N'})\in\overline\calb^+_{Z,Z'}$.
\end{lemma}
\begin{proof}
First we consider case (a).
\begin{itemize}
\item Suppose that $m'=m$.
Then
\begin{itemize}
\item $a'_k=a_{k+1}>d_{k+1}=d'_k$;

\item $d'_{k-1}=d_k\geq a_{k+1}=a'_k$;

\item $c'_{k-1}=d_{k-1}\geq a_k=b'_{k-1}$;

\item $b'_{k-2}=b_{k-2}\geq a_{k-1}>d_{k-1}=c'_{k-1}$; $b'_{k-1}=a_k>c_{k-1}=c'_k$.
\end{itemize}

\item Suppose that $m'=m+1$.
Then
\begin{itemize}
\item $a'_k=a_{k+1}\geq d_{k+1}=d'_k$;

\item $d'_{k-1}=d_k>a_{k+1}=a'_k$;

\item $c'_{k-1}=d_{k-1}>a_k=b'_{k-1}$;

\item $b'_{k-2}=b_{k-2}\geq a_{k-1}\geq d_{k-1}=c'_{k-1}$; $b'_{k-1}=a_k\geq c_{k-1}=c'_k$.
\end{itemize}
\end{itemize}
Because we just move $a_k,d_{k-1}$ from their respective rows,
for the requirement for $a'_i$, $b'_i$, $c'_i$, $d'_i$ in Lemma~\ref{0210} the above are all
the conditions which need to be checked.

The proofs for other cases are similar.
\end{proof}

\begin{proposition}\label{0208}
Let $Z,Z'$ be two special symbols of sizes $(m+1,m),(m',m')$ respectively
where $m'=m,m+1$.
If\/ $\overline\calb^+_{Z,Z'}\neq\emptyset$,
then $Z$ occurs in the relation $\cald_{Z,Z'}$.
\end{proposition}
\begin{proof}
Suppose that $\overline\calb^+_{Z,Z'}\neq\emptyset$.
Then there exists $\Lambda=\Lambda_M$ occurring in the relation $\overline\calb^+_{Z,Z'}$ for some $M\subset Z_\rmI$.
If $M=\emptyset$, then $\Lambda_M=Z$ and the lemma is proved.
Now we suppose that $M\neq\emptyset$.
So we are in one of Cases (a)--(f) above.
From Lemma~\ref{0205} and Lemma~\ref{0617},
we see that there exists $M'\subset Z_\rmI$ such that $\Lambda_{M'}$
also occurs in the relation $\overline\calb^+_{Z,Z'}$ and either $M'=\emptyset$ or $\max(M')<\max(M)$.
Repeat the same process,
we will finally conclude that $Z=\Lambda_\emptyset$ also occurs in the relation $\overline\calb^+_{Z,Z'}$.
Because $Z\in\cals_{Z,1}$,
then $Z$ occurs in $\calb^+_{Z,Z'}$, and in fact $Z$ occurs in $\cald_{Z,Z'}$.
\end{proof}

\begin{corollary}
$\cald_{Z,Z'}\neq\emptyset$ if and only if $\calb^+_{Z,Z'}\neq\emptyset$.
\end{corollary}
\begin{proof}
Because $\cald_{Z,Z'}$ is a subset of $\calb^+_{Z,Z'}$,
it is obvious that $\cald_{Z,Z'}\neq\emptyset$ implies $\calb^+_{Z,Z'}\neq\emptyset$.
Conversely, suppose that $\calb^+_{Z,Z'}\neq\emptyset$.
Then $\overline\calb^+_{Z,Z'}\neq\emptyset$,
and hence $\cald_{Z,Z'}\neq\emptyset$ by Proposition~\ref{0208}.
\end{proof}

Recall that $B^+_\Lambda,D_Z$ are defined in (\ref{0235}).

\begin{lemma}\label{0616}
Let $\Lambda\in\overline\cals_Z$.
Then $|B^+_\Lambda|\leq|D_Z|$.
\end{lemma}
\begin{proof}
If $\Lambda$ does not occur in $\calb^+_{Z,Z'}$, i.e., $B^+_\Lambda=\emptyset$,
then the lemma is clearly true.
Now suppose that $(\Lambda,\Lambda')\in\overline\calb^+_{Z,Z'}$ where $\Lambda=\Lambda_M$, $\Lambda'=\Lambda_N$
for some $M\subset Z_\rmI$ and $N\subset Z'_\rmI$.
Write $M=M^{(0)}$, $N=N^{(0)}$ and $(M^{(i)})'=M^{(i+1)}$, $(N^{(i)})'=N^{(i+1)}$ where
$(M^{(i)})'$ and $(N^{(i)})'$ are constructed as in the beginning of this subsection.
From the construction, we know that $\max(M^{(i+1)})<\max(M^{(i)})$
and $(\Lambda_{M^{(i+1)}},\Lambda_{N^{(i+1)}})\in\overline\calb^+_{Z,Z'}$.
Therefore, $M^{(t)}=\emptyset$ for some $t$ and then $\Lambda_{N^{(t)}}\in D_Z$.
So now it suffices to show that the function $f\colon B^+_\Lambda\rightarrow D_Z$ given by
$\Lambda'\mapsto\Lambda_{N^{(t)}}$ is one-to-one.

Suppose that
\[
\Lambda_M=\Lambda_{\bar M}=\binom{a_1,\ldots,a_{m_1}}{b_1,\ldots,b_{m_2}},\qquad
\Lambda_N=\binom{c_1,\ldots,c_{m'_1}}{d_1,\ldots,d_{m'_2}},\qquad
\Lambda_{\bar N}=\binom{\bar c_1,\ldots,\bar c_{m'_1}}{\bar d_1,\ldots,\bar d_{m'_2}}
\]
where $M=\bar M\neq\emptyset$ and $N\neq\bar N$,
and suppose that both $(\Lambda_M,\Lambda_N)$ and $(\Lambda_{\bar M},\Lambda_{\bar N})$ are in $\overline\calb^+_{Z,Z'}$.
Write
\[
f(\Lambda_N)=\binom{c_1^{(t)},\ldots,c^{(t)}_{m'}}{d_1^{(t)},\ldots,d^{(t)}_{m'}},\qquad
f(\Lambda_{\bar N})=\binom{\bar c^{(t)}_1,\ldots,\bar c^{(t)}_{m'}}{\bar d^{(t)}_1,\ldots,\bar d^{(t)}_{m'}}.
\]
Now we consider the following cases when $\max(M)=a_k$ for some $k$:
\begin{enumerate}
\item[(1)] Suppose that both $(\Lambda_M,\Lambda_N)$ and $(\Lambda_{\bar M},\Lambda_{\bar N})$ are in the same case
((a), (b) or (c)).
Then it is clear that $M'=\bar M'$ and $N'\neq\bar N'$,
hence $\Lambda_{M'}=\Lambda_{\bar M'}$ and $\Lambda_{N'}\neq\Lambda_{\bar N'}$.

\item[(2)] Suppose $(\Lambda_M,\Lambda_N)$ is in Case (a) and $(\Lambda_{\bar M},\Lambda_{\bar N})$ is in Case (b).
Note that now $\Lambda_{M'}=\Lambda_{\bar M'}$
and we have $\Lambda_{M'}$ and $\Lambda_{N'}$ by the following table:
\[
\begin{tabular}{c|c|c|c|c|c|c|c}
\toprule
$i$ & $1$ & $\ldots$ & $k-2$ & $k-1$ & $k$ & $k+1$ & $\ldots$ \\
\midrule
$a'_i$ & $a_1$ & $\ldots$ & $a_{k-2}$ & $a_{k-1}$ & $a_{k+1}$ & $a_{k+2}$ & $\ldots$ \\
\midrule
$b'_i$ & $b_1$ & $\ldots$ & $b_{k-2}$ & $a_k$ & $b_{k-1}$ & $b_k$ & $\ldots$ \\
\midrule
$c'_i$ & $c_1$ & $\ldots$ & $c_{k-2}$ & $d_{k-1}$ & $c_{k-1}$ & $c_k$ & $\ldots$ \\
\midrule
$d'_i$ & $d_1$ & $\ldots$ & $d_{k-2}$ & $d_k$ & $d_{k+1}$ & $d_{k+2}$ & $\ldots$ \\
\bottomrule
\end{tabular}
\]
Similarly, we have $\Lambda_{\bar M'}$ and $\Lambda_{\bar N'}$ by the following table:
\[
\begin{tabular}{c|c|c|c|c|c|c|c}
\toprule
$i$ & $1$ & $\ldots$ & $k-2$ & $k-1$ & $k$ & $k+1$ & $\ldots$ \\
\midrule
$\bar a'_i$ & $a_1$ & $\ldots$ & $a_{k-2}$ & $a_{k-1}$ & $a_{k+1}$ & $a_{k+2}$ & $\ldots$ \\
\midrule
$\bar b'_i$ & $b_1$ & $\ldots$ & $b_{k-2}$ & $a_k$ & $b_{k-1}$ & $b_k$ & $\ldots$ \\
\midrule
$\bar c'_i$ & $\bar c_1$ & $\ldots$ & $\bar c_{k-2}$ & $\bar c_{k-1}$ & $\bar d_k$ & $\bar c_k$ & $\ldots$ \\
\midrule
$\bar d'_i$ & $\bar d_1$ & $\ldots$ & $\bar d_{k-2}$ & $\bar d_{k-1}$ & $\bar d_{k+1}$ & $\bar d_{k+2}$ & $\ldots$ \\
\bottomrule
\end{tabular}
\]
Note that
\[
\begin{cases}
\bar d'_{k-1}=\bar d_{k-1}\geq a_k>d_k=d'_{k-1}, & \text{if $m'=m$};\\
\bar d'_{k-1}=\bar d_{k-1}>a_k\geq d_k=d'_{k-1}, & \text{if $m'=m+1$},
\end{cases}
\]
by Lemma~\ref{0210},
so we have $\Lambda_{N'}\neq\Lambda_{\bar N'}$.

\item[(3)] Suppose $(\Lambda_M,\Lambda_N)$ is in Case (b) and $(\Lambda_{\bar M},\Lambda_{\bar N})$ is in  Case (c).
Then we have $\Lambda_{M'}$ and $\Lambda_{N'}$ by the following table:
\[
\begin{tabular}{c|c|c|c|c|c|c|c}
\toprule
$i$ & $1$ & $\ldots$ & $k-2$ & $k-1$ & $k$ & $k+1$ & $\ldots$ \\
\midrule
$a'_i$ & $a_1$ & $\ldots$ & $a_{k-2}$ & $a_{k-1}$ & $a_{k+1}$ & $a_{k+2}$ & $\ldots$ \\
\midrule
$b'_i$ & $b_1$ & $\ldots$ & $b_{k-2}$ & $a_k$ & $b_{k-1}$ & $b_k$ & $\ldots$ \\
\midrule
$c'_i$ & $c_1$ & $\ldots$ & $c_{k-2}$ & $c_{k-1}$ & $d_k$ & $c_k$ & $\ldots$ \\
\midrule
$d'_i$ & $d_1$ & $\ldots$ & $d_{k-2}$ & $d_{k-1}$ & $d_{k+1}$ & $d_{k+2}$ & $\ldots$ \\
\bottomrule
\end{tabular}
\]
Now because $a'_k=a_{k+1}<a_k=b'_{k-1}$,
by Lemma~\ref{0615}, we see that $\max(M')\neq a'_k$.
Because now $(\Lambda_M,\Lambda_N)$ is in Case (b),
we have $d_k\geq b_{k-1}$ (resp.~$d_k>b_{k-1}$) if $m'=m$ (resp.~$m'=m+1$).
So if $\max(M')=b'_k=b_{k-1}$,
we must be in the Cases (e) or (f).
Then the entry $c'_k$ is not changed during the process $N'\mapsto N''$.
For all other possibilities, i.e., $\max(M')$ is $a'_l$ or $b'_l$ for
$l>k$, the entry $c'_k$ is also clearly unchanged during the process $N'\mapsto N''$.
Therefore we can conclude that $c^{(t)}_k=c'_k=d_k$.

Similarly, we have $\Lambda_{\bar M'}$ and $\Lambda_{\bar N'}$ by the following table:
\[
\begin{tabular}{c|c|c|c|c|c|c|c}
\toprule
$i$ & $1$ & $\ldots$ & $k-2$ & $k-1$ & $k$ & $k+1$ & $\ldots$ \\
\midrule
$\bar a'_i$ & $a_1$ & $\ldots$ & $a_{k-2}$ & $a_{k-1}$ & $b_{k-1}$ & $a_{k+1}$ & $\ldots$ \\
\midrule
$\bar b'_i$ & $b_1$ & $\ldots$ & $b_{k-2}$ & $a_k$ & $b_k$ & $b_{k+1}$ & $\ldots$ \\
\midrule
$\bar c'_i$ & $\bar c_1$ & $\ldots$ & $\bar c_{k-2}$ & $\bar c_{k-1}$ & $\bar c_k$ & $\bar c_{k+1}$ & $\ldots$ \\
\midrule
$\bar d'_i$ & $\bar d_1$ & $\ldots$ & $\bar d_{k-2}$ & $\bar d_{k-1}$ & $\bar d_k$ & $\bar d_{k+1}$ & $\ldots$ \\
\bottomrule
\end{tabular}
\]
Now because $\bar a'_k=b_{k-1}<a_k=\bar b'_{k-1}$ by Lemma~\ref{0615},
we see that $\max(\bar M')\neq\bar a'_k$.
Because now $(\Lambda_{\bar M},\Lambda_{\bar N})$ is in Case (c),
we have $\bar d_{k-1}\geq b_{k-1}>b_k=\bar b_k$ from the proof of Lemma~\ref{0205}.
So if $\max(\bar M')=\bar b'_k=b_k$,
then $(\Lambda_{\bar M'},\Lambda_{\bar N'})$ is in Case (e) or (f).
By the similar argument, we have $\bar c^{(t)}_k=\bar c'_k=\bar c_k$.

Now $d_k\geq b_{k-1}>\bar c_k$ if $m'=m$,
and $d_k>b_{k-1}\geq \bar c_k$ if $m'=m+1$ by Lemma~\ref{0210}.
This means that $c^{(t)}_k\neq\bar c^{(t)}_k$, hence
$f(\Lambda_N)\neq f(\Lambda_{\bar N})$.

\item[(4)] Suppose $(\Lambda_M,\Lambda_N)$ is in Case (a) and $(\Lambda_{\bar M},\Lambda_{\bar N})$ is in  Case (c).
Then we have $\Lambda_{\bar M'}$ and $\Lambda_{\bar N'}$ by the following table:
\[
\begin{tabular}{c|c|c|c|c|c|c|c}
\toprule
$i$ & $1$ & $\ldots$ & $k-2$ & $k-1$ & $k$ & $k+1$ & $\ldots$ \\
\midrule
$\bar a'_i$ & $a_1$ & $\ldots$ & $a_{k-2}$ & $a_{k-1}$ & $b_{k-1}$ & $a_{k+1}$ & $\ldots$ \\
\midrule
$\bar b'_i$ & $b_1$ & $\ldots$ & $b_{k-2}$ & $a_k$ & $b_k$ & $b_{k+1}$ & $\ldots$ \\
\midrule
$\bar c'_i$ & $\bar c_1$ & $\ldots$ & $\bar c_{k-2}$ & $\bar c_{k-1}$ & $\bar c_k$ & $\bar c_{k+1}$ & $\ldots$ \\
\midrule
$\bar d'_i$ & $\bar d_1$ & $\ldots$ & $\bar d_{k-2}$ & $\bar d_{k-1}$ & $\bar d_k$ & $\bar d_{k+1}$ & $\ldots$ \\
\bottomrule
\end{tabular}
\]
By the same argument in (3), we have $\bar c^{(t)}_k=\bar c_k$.

Similarly, we have $\Lambda_{M'}$ and $\Lambda_{N'}$ by the following table:
\[
\begin{tabular}{c|c|c|c|c|c|c|c}
\toprule
$i$ & $1$ & $\ldots$ & $k-2$ & $k-1$ & $k$ & $k+1$ & $\ldots$ \\
\midrule
$a'_i$ & $a_1$ & $\ldots$ & $a_{k-2}$ & $a_{k-1}$ & $a_{k+1}$ & $a_{k+2}$ & $\ldots$ \\
\midrule
$b'_i$ & $b_1$ & $\ldots$ & $b_{k-2}$ & $a_k$ & $b_{k-1}$ & $b_k$ & $\ldots$ \\
\midrule
$c'_i$ & $c_1$ & $\ldots$ & $c_{k-2}$ & $d_{k-1}$ & $c_{k-1}$ & $c_k$ & $\ldots$ \\
\midrule
$d'_i$ & $d_1$ & $\ldots$ & $d_{k-2}$ & $d_k$ & $d_{k+1}$ & $d_{k+2}$ & $\ldots$ \\
\bottomrule
\end{tabular}
\]
As in (3), we know that $\max(M')\neq a'_k$.
Suppose that $\max(M')=b'_k$.
We have the following two situations:
\begin{enumerate}
\item Suppose that $d'_{k-1}\geq b'_k$ (resp.~$d'_{k-1}>b'_k$) if $m'=m$ (resp.~$m'=m+1$),
i.e., $(\Lambda_{M'},\Lambda_{N'})$ is in Case (e) or (f).
By the same argument in (3), we have $c^{(t)}_k=c_{k-1}$.
Now $c_{k-1}\geq b_{k-1}>\bar c_k$ if $m'=m$,
and $c_{k-1}>b_{k-1}\geq \bar c_k$ if $m'=m+1$ by Lemma~\ref{0210}.
This means that $c^{(t)}_k\neq\bar c^{(t)}_k$, hence
$f(\Lambda_N)\neq f(\Lambda_{\bar N})$.

\item Suppose that $b'_k>d'_{k-1}$ (resp.~$b'_k\geq d'_{k-1}$) if $m'=m$ (resp.~$m'=m+1$),
i.e., $(\Lambda_{M'},\Lambda_{N'})$ is in Case (d).
Then we have the following table for $\Lambda_{M''}$ and $\Lambda_{N''}$:
\[
\begin{tabular}{c|c|c|c|c|c|c|c}
\toprule
$i$ & $1$ & $\ldots$ & $k-2$ & $k-1$ & $k$ & $k+1$ & $\ldots$ \\
\midrule
$a''_i$ & $a_1$ & $\ldots$ & $a_{k-2}$ & $a_{k-1}$ & $b_{k-1}$ & $a_{k+1}$ & $\ldots$ \\
\midrule
$b''_i$ & $b_1$ & $\ldots$ & $b_{k-2}$ & $a_k$ & $b_k$ & $b_{k+1}$ & $\ldots$ \\
\midrule
$c''_i$ & $c_1$ & $\ldots$ & $c_{k-2}$ & $d_{k-1}$ & $c_k$ & $c_{k+1}$ & $\ldots$ \\
\midrule
$d''_i$ & $d_1$ & $\ldots$ & $d_{k-2}$ & $c_{k-1}$ & $d_k$ & $d_{k+1}$ & $\ldots$ \\
\bottomrule
\end{tabular}
\]
Now $(\Lambda_M,\Lambda_N)$ is in Case (a), we have
$a_k>c_{k-1}$ (resp.~$a_k\geq c_{k-1}$) if $m'=m$ (resp.~$m'=m+1$).
Then by Lemma~\ref{0210},
\[
\begin{cases}
\bar d'_{k-1}=d_{k-1}\geq a_k>c_{k-1}=d''_{k-1}, & \text{if $m'=m$};\\
\bar d'_{k-1}=d_{k-1}>a_k\geq c_{k-1}=d''_{k-1}, & \text{if $m'=m+1$},
\end{cases}
\]
in particular, $\bar d'_{k-1}\neq d''_{k-1}$.
Then we have $\Lambda_{\bar M'}=\Lambda_{M''}$ and $\Lambda_{\bar N'}\neq\Lambda_{N''}$.
\end{enumerate}
\end{enumerate}

Next we consider the case that $\max(M)=b_k$ for some $k$:
\begin{enumerate}
\item[(5)] Suppose that both $(\Lambda_M,\Lambda_N)$ and $(\Lambda_{\bar M},\Lambda_{\bar N})$ are in the same case
((d), (e) or (f)).
The proof is similar to (1).

\item[(6)] Suppose $(\Lambda_M,\Lambda_N)$ is in Case (d) and $(\Lambda_{\bar M},\Lambda_{\bar N})$ is in Case (e).
The proof is similar to (2).

\item[(7)] Suppose $(\Lambda_M,\Lambda_N)$ is in Case (e) and $(\Lambda_{\bar M},\Lambda_{\bar N})$ is in  Case (f).
The proof is similar to (3).

\item[(8)] Suppose $(\Lambda_M,\Lambda_N)$ is in Case (d) and $(\Lambda_{\bar M},\Lambda_{\bar N})$ is in  Case (f).
The proof is similar to (4).
\end{enumerate}
For all above 8 cases,
either we can show directly $f(\Lambda_N)\neq f(\Lambda_{\bar N})$
or we show that $\Lambda_{M'}=\Lambda_{\bar M'}$ and $\Lambda_{N'}\neq\Lambda_{\bar N'}$
(or $\Lambda_{M''}=\Lambda_{\bar M'}$ and $\Lambda_{N''}\neq\Lambda_{\bar N'}$).
For the second situation, we can repeat the same process and reach the conclusion $f(\Lambda_N)\neq f(\Lambda_{\bar N})$ finally.
Therefore $f$ is one--to-one.
\end{proof}

\begin{remark}
We shall see in the next section that in fact $|B^+_\Lambda|=|D_Z|$ if $B^+_\Lambda\neq\emptyset$.
\end{remark}

\subsection{The set $\cald_{Z,Z'}$}
Suppose that $\Psi_1,\ldots,\Psi_k$ is a set of pairwise disjoint consecutive pairs in $Z_\rmI$.
Then decomposition of $\Psi=\Psi_1\cup\Psi_2\cup\cdots\cup\Psi_k$ as a disjoint union
of consecutive pairs is clearly unique.
The set $\Psi$ is called a \emph{set of consecutive pairs}.
We denote $\Psi'\leq\Psi$ if $\Psi'$ is a union of a subset of $\{\Psi_1,\ldots,\Psi_k\}$.
Clearly $|\{\,\Psi'\mid\Psi'\leq\Psi\,\}|=2^k$.

Write
\begin{equation}\label{0801}
Z'=\binom{s'_1,s'_2,\ldots,s'_{m'}}{t'_1,t'_2,\ldots, t'_{m'}}.
\end{equation}

\begin{lemma}\label{0206}
Let $\Psi$ be a subset of pairs in $Z_\rmI'$, and let $x=\max(\Psi)$.
\begin{enumerate}
\item[(i)] If $x=s'_i$ and $t'_i\not\in\Psi$,
then $(Z,\Lambda_\Psi)\not\in\cald_{Z,Z'}$.

\item[(ii)] If $x=t'_i\in\Psi$ and $s'_{i+1}\not\in\Psi$,
then $(Z,\Lambda_\Psi)\not\in\cald_{Z,Z'}$.
\end{enumerate}
\end{lemma}
\begin{proof}
Write
\[
Z=\binom{a_1,a_2,\ldots,a_{m+1}}{b_1,b_2,\ldots,b_m},\qquad
\Lambda_\Psi=\binom{c_1,c_2,\ldots,c_{m'}}{d_1,d_2,\ldots,d_{m'}}.
\]
First suppose that $x=s'_i$ and $t'_i\not\in\Psi$.
This means that $t'_{i-1},s'_i,t'_i$ are all in the second row of $\Lambda_\Psi$,
so we can write $t'_{i-1}=d_{l-1}$, $s'_i=d_l$, $t'_i=d_{l+1}$ for some $l$.
Now $x=s'_i$ is the largest entry in $\Psi$,
this means that $s'_1,\ldots,s'_{i-1},t'_1,\ldots,t'_{i-1}$ are all in their natural positions in $\Lambda_\Psi$.
Hence we have $l=i$, $c_{l-1}=s'_{i-1}$, and $c_l\leq t'_i=d_{l+1}$.
If $(Z,\Lambda_\Psi)\in\cald_{Z,Z'}$,
then by Lemma~\ref{0210} we have
\[
\begin{cases}
c_l\geq b_l\geq a_{l+1}>d_{l+1}, &  \text{if $m'=m$}; \\
c_l>b_l\geq a_{l+1}\geq d_{l+1}, &  \text{if $m'=m+1$}.
\end{cases}
\]
So we have $c_l>d_{l+1}$ for both cases and this inequality contradicts the inequality $c_l\leq d_{l+1}$
we just obtained.

The proof for (ii) is similar.
\end{proof}

\begin{lemma}\label{0207}
Suppose that $(Z,\Lambda_\Psi)\in\cald_{Z,Z'}$ for some nonempty subset of pairs $\Psi\subset Z'_\rmI$.
Let $s'$ (resp.~$t'$) be the maximal entry in the first (resp.~second) row of\/ $\Psi$.
Then $\binom{s'}{t'}$ is a consecutive pair,
and $(Z,\Lambda_{\Psi'})\in\cald_{Z,Z'}$ where $\Psi'=\Psi\smallsetminus\binom{s'}{t'}$.
\end{lemma}
\begin{proof}
Write
\[
Z=\binom{a_1,a_2,\ldots,a_{m+1}}{b_1,b_2,\ldots,b_m},\quad
\Lambda_\Psi=\binom{c_1,c_2,\ldots,c_{m'}}{d_1,d_2,\ldots,d_{m'}},\quad
\Lambda_{\Psi'}=\binom{c'_1,c'_2,\ldots,c'_{m'}}{d'_1,d'_2,\ldots,d'_{m'}}
\]
where $m'=m$ or $m+1$.
Because $Z$ is special,
we have
\begin{equation}\label{0211}
a_i\geq b_i\geq a_{i+1}\qquad\text{for each $i$}.
\end{equation}
Because $(Z,\Lambda_\Psi)\in\cald_{Z,Z'}$, by Lemma~\ref{0206}, we see that if $s'>t'$,
then $s'=s'_i$ and $'t=t'_i$ for some $i$;
if $t'>s'$, then $s'=s'_{i+1}$ and $t'=t'_i$ for some $i$.
Then $\binom{s'}{t'}$ is a consecutive pair.
Write $s'=d_l$ and $t'=c_k$ for some $k,l$.
Note that $\Lambda_{\Psi'}$ is obtained from $\Lambda_\Psi$ by exchanging the positions of $c_k,d_l$
since $\binom{c_k}{d_l}$ is a consecutive pair.
\begin{enumerate}
\item[(1)] Suppose that $c_k>d_l$.
In this case, the entries $c_1,\ldots,c_{k-1}$ in $\Lambda_\Psi$ are in their natural position,
we must have $d_{k-2}>c_k>d_{k-1}$ and hence $l=k-1$.
Then we have the following table for $\Lambda_{\Psi'}$:
\[
\begin{tabular}{c|c|c|c|c|c|c|c}
\toprule
$i$ & $1$ & $\ldots$ & $k-2$ & $k-1$ & $k$ & $k+1$ & $\ldots$ \\
\midrule
$c'_i$ & $c_1$ & $\ldots$ & $c_{k-2}$ & $c_{k-1}$ & $d_{k-1}$ & $c_{k+1}$ & $\ldots$ \\
\midrule
$d'_i$ & $d_1$ & $\ldots$ & $d_{k-2}$ & $c_k$ & $d_k$ & $d_{k+1}$ & $\ldots$ \\
\bottomrule
\end{tabular}
\]
\begin{enumerate}
\item Suppose that $m'=m$.
Then by Lemma~\ref{0210} and (\ref{0211}),
the assumption $(Z,\Lambda_\Psi)\in\cald_{Z,Z'}$ implies that
\begin{itemize}
\item $a_{k-1}\geq b_{k-1}>c_k=d'_{k-1}$;

\item $d'_{k-1}=c_k>d_{k-1}\geq a_k$;

\item $c'_k=d_{k-1}\geq a_k\geq b_k$;

\item $b_{k-1}>c_k>d_{k-1}=c'_k$.
\end{itemize}

\item Suppose that $m'=m+1$.
Then by the same argument as in (a), we have
\begin{itemize}
\item $a_{k-1}\geq b_{k-1}\geq c_k=d'_{k-1}$;

\item $d'_{k-1}=c_k>d_{k-1}>a_k$;

\item $c'_k=d_{k-1}>a_k\geq b_k$;

\item $b_{k-1}\geq c_k>d_{k-1}=c'_k$.
\end{itemize}
\end{enumerate}
Because $c'_i=c_i$ for $i\neq k$ and $d'_i=d_i$ for $i\neq k-1$,
the above are all the inequalities that we need to check.
Then by Lemma~\ref{0210}, we conclude that
$(Z,\Lambda_{\Psi'})\in\cald_{Z,Z'}$ for both Cases (a) and (b).

\item[(2)] Suppose that $d_l>c_k$.
In this case, the entries $d_1,\ldots,d_{l-1}$ in $\Lambda_\Psi$ are in their natural position,
so we have $c_{l-1}>d_l>c_l$ and hence $l=k$.
The remaining proof is similar to that for Case (1).
\end{enumerate}
\end{proof}

\begin{lemma}\label{0209}
If $(Z,\Lambda_\Psi)\in\cald_{Z,Z'}$ for some subset of pairs $\Psi$ of $Z'_\rmI$,
then $\Psi$ consists of consecutive pairs.
\end{lemma}
\begin{proof}
Write $Z'$ as in (\ref{0801}).
Suppose $(Z,\Lambda_\Psi)\in\cald_{Z,Z'}$ such that $\Psi$ can not expressed as a subset of consecutive pairs.
This means that there are (1) $s'_i\in\Psi$ and $t'_{i-1},t'_i\not\in\Psi$ for some $i$;
or (2) $t'_i\in\Psi$ and $s'_i,s'_{i+1}\not\in\Psi$ for some $i$.
By Lemma~\ref{0207}, after deleting some consecutive pairs if necessary,
we may assume the entry $s'_i$ in (1) (resp.~$t'_i$ in (2)) is the maximal entry of $\Psi$.
Then we get a contradiction from Lemma~\ref{0206}.
\end{proof}

\begin{proposition}\label{0216}
Let $Z,Z'$ be special symbols of defect $1,0$ respectively.
If $\cald_{Z,Z'}\neq\emptyset$, then $(Z,Z')\in\cald_{Z,Z'}$.
\end{proposition}
\begin{proof}
Suppose that the sizes of $Z,Z'$ are $(m+1,m),(m',m')$ respectively.
We know that $m'=m,m+1$ from Lemma~\ref{0213}.
Suppose that $\cald_{Z,Z'}\neq\emptyset$.
Then $\overline\calb^+_{Z,Z'}\neq\emptyset$ and by Proposition\ \ref{0208},
we see that $(Z,\Lambda_\Psi)\in\cald_{Z,Z'}$ for some subset $\Psi\subset Z'_\rmI$.
Because now ${\rm def}(\Lambda_\Psi)=0$,
we have $|\Psi^*|=|\Psi_*|$, i.e., $\Psi$ is a subset of pairs in $Z'_\rmI$.
Then by Lemma\ \ref{0209}, we see that $\Psi$ consists of consecutive pairs.
By Lemma\ \ref{0207}, we have $(Z,\Lambda_{\Psi'})\in\cald_{Z,Z'}$ where $\Psi'$ is obtained from $\Psi$
by removing a consecutive pair.
Repeat the same process, we conclude that $(Z,Z')\in\cald_{Z,Z'}$.
\end{proof}


\section{The Relation $\calb^+_{Z,Z'}$ I: Regular and One-to-one Case}\label{0706}
The main purpose (Proposition~\ref{0806}) is to show that the statement in Theorem~\ref{0310}
holds when $\epsilon=+$, both $Z,Z'$ are regular, and $\cald_{Z,Z'}$ is one-to-one.
The strategy is to reduced the case via an isometry of inner product spaces to the case
considered in Proposition~\ref{0512} and Proposition~\ref{0513}.
In this section,
we always assume that $Z$ is a special symbol of size $(m+1,m)$ (hence of defect $1$),
and $Z'$ is a special symbol of size $(m',m')$ (hence of defect $0$) such that $m'=m$ or $m'=m+1$.

\subsection{Core of the relation $\cald_{Z,Z'}$}
Let
\[
\Lambda=\binom{a_1,a_2,\ldots,a_{m_1}}{b_1,b_2,\ldots,b_{m_2}}\in\overline\cals_Z,\qquad
\Lambda'=\binom{c_1,c_2,\ldots,c_{m_1'}}{d_1,d_2,\ldots,d_{m'_2}}\in\overline\cals_{Z'}
\]
such that $(\Lambda,\Lambda')\in\overline\calb^+_{Z,Z'}$.
Suppose that $\Lambda'=\Lambda_N$ for some $N\subset Z'_\rmI$ and $\binom{c_k}{d_l}$ is a consecutive pair such that
either $\{c_k,d_l\}\subset N$ or $\{c_k,d_l\}\cap N=\emptyset$.
Let $N'=(N\cup\binom{c_k}{d_l})\smallsetminus(N\cap\binom{c_k}{d_l})$.
Then we know that $\Lambda_{N'}=\binom{c'_1,c'_2,\ldots,c'_{m_1'}}{d'_1,d'_2,\ldots,d'_{m'_2}}$
is obtained from $\Lambda_N$ by moving $c_k$ to the second row and moving $d_l$ to the first row.
Because now $c_k,d_l$ are consecutive,
we see that $d_l$ is moved to the original position of $c_k$ and vice versa, i.e.,
we have the following tables:
\[
\begin{tabular}{c|c|c|c|c|c|c}
\toprule
$i$ & $1$ & $\ldots$ & $k-1$ & $k$ & $k+1$ & $\ldots$ \\
\midrule
$c'_i$ & $c_1$ & $\ldots$ & $c_{k-1}$ & $d_l$ & $c_{k+1}$ & $\ldots$ \\
\bottomrule
\end{tabular}
\]
\[
\begin{tabular}{c|c|c|c|c|c|c}
\toprule
$i$ & $1$ & $\ldots$ & $l-1$ & $l$ & $l+1$ & $\ldots$ \\
\midrule
$d'_i$ & $d_1$ & $\ldots$ & $d_{l-1}$ & $c_k$ & $d_{l+1}$ & $\ldots$ \\
\bottomrule
\end{tabular}
\]
Then $(\Lambda,\Lambda_{N'})\in\overline\calb^+_{Z,Z'}$ if and only if
\[
\begin{cases}
a_l>d'_l=c_k\geq a_{l+1}\text{ and }b_{k-1}>c'_k=d_l\geq b_k, & \text{if $m'=m$}; \\
a_l\geq d'_l=c_k>a_{l+1}\text{ and }b_{k-1}\geq c'_k=d_l>b_k, & \text{if $m'=m+1$}
\end{cases}
\]
by Lemma~\ref{0210}.
Now consider the special case that
\begin{equation}\label{0311}
Z=\binom{a_1,a_2,\ldots,a_{m+1}}{b_1,b_2,\ldots,b_m},\qquad
Z'=\binom{c_1,c_2,\ldots,c_{m'}}{d_1,d_2,\ldots,d_{m'}}.
\end{equation}
If $\binom{c_k}{d_l}$ is a consecutive pair in $Z'_\rmI$,
then $l=k-1,k$.

\begin{lemma}\label{0301}
Let $Z,Z'$ be given as in (\ref{0311}), and let $\Psi'=\binom{c_k}{d_l}$ be a consecutive pair in $Z'_\rmI$.
Suppose that $\cald_{Z,Z'}\neq\emptyset$.
Then $(Z,\Lambda_{\Psi'})\in\cald_{Z,Z'}$ if and only if
\[
\begin{cases}
a_k>c_k\text{ and }d_k\geq b_k, & \text{if $m'=m$ and $l=k$};\\
c_k\geq a_k\text{ and }b_{k-1}>d_{k-1}, & \text{ if $m'=m$ and $l=k-1$};\\
a_k\geq c_k\text{ and }d_k>b_k, & \text{ if $m'=m+1$ and $l=k$};\\
c_k>a_k\text{ and }b_{k-1}\geq d_{k-1}, & \text{ if $m'=m+1$ and $l=k-1$}.
\end{cases}
\]
\end{lemma}
\begin{proof}
The assumption $\cald_{Z,Z'}\neq\emptyset$ implies that $(Z,Z')\in\cald_{Z,Z'}$ by Proposition~\ref{0216}.
By Lemma~\ref{0210},
we have
\[
\begin{cases}
a_l>d_l\geq a_{l+1}\text{ and }b_{k-1}>c_k\geq b_k, & \text{if $m'=m$};\\
a_l\geq d_l>a_{l+1}\text{ and }b_{k-1}\geq c_k>b_k. & \text{if $m'=m+1$}.
\end{cases}
\]
First suppose that $m'=m$ and $l=k$.
From the discussion before the lemma, we see that
the condition $(Z,\Lambda_{\Psi'})\in\cald_{Z,Z'}$ is equivalent to
\[
a_k>c_k\geq a_{k+1}\qquad\text{and}\qquad
b_{k-1}>d_k\geq b_k.
\]
However, the inequalities $c_k\geq a_{k+1}$ and $b_{k-1}>d_k$ are implied
by the assumptions that $Z,Z'$ are special and $(Z,Z')\in\cald_{Z,Z'}$:
$c_k\geq b_k\geq a_{k+1}$ and $b_{k-1}\geq a_k>d_k$.
Hence the lemma is proved for the first case.
The proof for other cases are similar.
\qed
\end{proof}

\begin{lemma}\label{0701}
Let $Z,Z'$ be given as in (\ref{0311}), and let $\Psi=\binom{a_k}{b_l}$ be a consecutive pair in $Z_\rmI$.
Suppose that $\cald_{Z,Z'}\neq\emptyset$.
Then $(\Lambda_\Psi,Z')\in\cald_{Z,Z'}$ if and only if
\[
\begin{cases}
c_k\geq a_k\text{ and }b_k>d_k, & \text{ if $m'=m$ and $l=k$};\\
a_k>c_k\text{ and }d_{k-1}\geq b_{k-1}, & \text{ if $m'=m$ and $l=k-1$};\\
c_k>a_k\text{ and }b_k\geq d_k, & \text{if $m'=m+1$ and $l=k$};\\
a_k\geq c_k\text{ and }d_{k-1}>b_{k-1}, & \text{ if $m'=m+1$ and $l=k-1$}.
\end{cases}
\]
\end{lemma}
\begin{proof}
The proof is similar to that of Lemma~\ref{0301}.
\end{proof}

\begin{lemma}\label{0305}
Suppose that $\Psi=\binom{c_k}{d_l}$ and $\Psi'=\binom{c_{k'}}{d_{l'}}$ are two consecutive pairs in $Z'_\rmI$.
If $(Z,\Lambda_\Psi)$ and $(Z,\Lambda_{\Psi'})$ are both in $\cald_{Z,Z'}$,
then $\Psi,\Psi'$ are disjoint.
\end{lemma}
\begin{proof}
Write $Z,Z'$ as in (\ref{0311}).
Suppose that $(Z,\Lambda_\Psi),(Z,\Lambda_{\Psi'})\in\cald_{Z,Z'}$ and $\Psi,\Psi'$ are not disjoint.
First assume that $\Psi=\binom{c_k}{d_{k-1}}$ and $\Psi'=\binom{c_k}{d_k}$ for some $k$.
Then by Lemma~\ref{0301} we have both $a_k>c_k$ and $c_k\geq a_k$ if $m'=m$;
$a_k\geq c_k$ and $c_k>a_k$ if $m'=m+1$.
Hence we get a contradiction.

The other possibility is $\Psi=\binom{c_l}{d_l},\Psi'=\binom{c_{l+1}}{d_l}$ for some $l$,
the proof is similar.
\end{proof}

\begin{proposition}\label{0307}
Let $(\bfG,\bfG')=(\Sp_{2n},\rmO^\epsilon_{2n'})$,
and let $Z,Z'$ be special symbols of ranks $n,n'$ and defects $1,0$ respectively.
Suppose that $\cald_{Z,Z'}\neq\emptyset$.
Then
\begin{enumerate}
\item[(i)] $D_Z=\{\,\Lambda_{\Psi'}\mid\Psi'\leq\Psi_0'\,\}$ for some subset of consecutive pairs 
$\Psi_0'$ of $Z'_\rmI$;

\item[(ii)] $D_{Z'}=\{\,\Lambda_\Psi\mid\Psi\leq\Psi_0\,\}$ for some subset of consecutive pairs 
$\Psi_0$ of $Z_\rmI$.
\end{enumerate}
\end{proposition}
\begin{proof}
Let $\Psi_1,\ldots,\Psi_k$ be the set of all consecutive pairs in $Z'_\rmI$ such that
$(Z,\Lambda_{\Psi_i})\in\cald_{Z,Z'}$.
By Lemma~\ref{0305}, we know that $\Psi_i,\Psi_j$ are disjoint if $i\neq j$.
Define $\Psi'_0=\Psi_1\cup\Psi_2\cup\cdots\cup\Psi_k$.

Suppose that $\Psi'$ is a set of consecutive pairs in $Z'_\rmI$.
Then we have a unique decomposition $\Psi'=\Psi'_1\cup\cdots\cup\Psi'_{k'}$
as a union of (disjoint) consecutive pairs.
Then by Lemma~\ref{0301}, we see that
$(Z,\Lambda_{\Psi'})\in\cald_{Z,Z'}$ if and only if $(Z,\Lambda_{\Psi'_i})\in\cald_{Z,Z'}$
for each $i=1,\ldots,k'$.
Therefore, by Lemma~\ref{0305},
$(Z,\Lambda_{\Psi'})\in\cald_{Z,Z'}$ implies that $\{\Psi'_1,\ldots,\Psi'_{k'}\}$ is a subset of
$\{\Psi_1,\ldots,\Psi_k\}$ and then $\Psi'\leq \Psi'_0$.
From Lemma~\ref{0209},
we know that each $\Lambda'\in D_Z$ is of the form $\Lambda_{\Psi'}$ for some subset of consecutive pairs $\Psi'$ of $Z'_\rmI$.
Hence, part (i) of the lemma follows.
The proof of (ii) is similar.
\end{proof}

The set $\Psi_0$ (resp.~$\Psi'_0$) in Proposition~\ref{0307} will be called the \emph{core} in $Z_\rmI$ (resp.~$Z'_\rmI$) of
the relation $\cald_{Z,Z'}$.
The relation $\cald_{Z,Z'}$ is said to be \emph{one-to-one} if
$\cald_{Z,Z'}\neq\emptyset$ and $\Psi_0=\Psi'_0=\emptyset$, i.e.,
if $D_Z=\{Z'\}$ and $D_{Z'}=\{Z\}$.

\subsection{Decomposition of $\calb^+_{Z,Z'}$}\label{0712}

\begin{lemma}\label{0809}
Suppose that $\cald_{Z,Z'}\neq\emptyset$.
Let $\Psi'_0$ be the core in $Z'_\rmI$ of $\cald_{Z,Z'}$,
and let $N$ be a subset of $Z'_\rmI$.
Suppose that there exists a consecutive pair $\Psi'$ in $\Psi'_0$ such that $|N\cap\Psi'|=1$.
Then $\Lambda_N$ does not occur in the relation $\overline\calb^+_{Z,Z'}$.
\end{lemma}
\begin{proof}
Write $Z$ and $Z'$ as in (\ref{0311}), and $\Psi'=\binom{c_k}{d_l}$.
Because $N$ contains exactly one of $c_k,d_l$,
these two entries will be in the same row of $\Lambda_N$.
Because $\binom{c_k}{d_l}$ is consecutive,
we know that $l=k,k-1$.
Note that $Z$ (resp.~$Z'$) is of size $(m+1,m)$ (resp.~$(m',m')$) such that $m'=m,m+1$.
Suppose that $(\Lambda,\Lambda_N)\in\overline\calb^+_{Z,Z'}$ for some $\Lambda\in\overline\cals_Z$,
and write
\[
\Lambda=\binom{a'_1,a'_2,\ldots,a'_{m_1}}{b'_1,b'_2,\ldots,b'_{m_2}}
\qquad\Lambda_N=\binom{c'_1,c'_2,\ldots,c'_{m'_1}}{d'_1,d'_2,\ldots,d'_{m'_2}}.
\]

\begin{enumerate}
\item[(1)] Suppose that $m'=m$ and $l=k$.
Because $\cald_{Z,Z'}\neq\emptyset$, we must have $(Z,Z')\in\cald_{Z,Z'}$ by Proposition~\ref{0216},
and hence $b_{k-1}>c_k\geq b_k$ and $a_k>d_k\geq a_{k+1}$ by Lemma~\ref{0210}.
Because now $(Z,\Lambda_{\Psi'})\in\cald_{Z,Z'}$, by Lemma~\ref{0301} and Lemma~\ref{0307},
we have $a_k>c_k\geq a_{k+1}$ and $b_{k-1}>d_k\geq b_k$.
This means that any entry of $Z$ is either (a) greater than $c_k$, or (b) less than or equal to $d_k$,
i.e., there is no entry of $Z$ which lies between $c_k$ and $d_k$.

Suppose $c_k,d_k$ are both in the second row of $\Lambda_N$,
so $c_k=d'_{k'}$ and $d_k=d'_{k'+1}$ for some $k'$.
By Lemma~\ref{0210}, the assumption $(\Lambda,\Lambda_N)\in\overline\calb^+_{Z,Z'}$ implies that
$a'_{k'}>d'_{k'}=c_k\geq a'_{k'+1}$,
similarly, $a'_{k'+1}>d'_{k'+1}=d_k\geq a'_{k'+2}$.
Therefore $c_k\geq a'_{k'+1}>d_k$ and we get a contradiction.

If both $c_k,d_k$ are in the first row of $\Lambda_N$,
then $c_k=c'_{k'}$ and $d_k=c'_{k'+1}$ for some $k'$.
By Lemma~\ref{0210}, the assumption $(\Lambda,\Lambda_N)\in\overline\calb^+_{Z,Z'}$ implies that
$b'_{k'-1}>c'_{k'}=c_k\geq b'_{k'}$, similarly, $b'_{k'}>c'_{k'+1}=d_k\geq b'_{k'+1}$.
Then $c_k\geq b'_{k'}>d_k$ and we get a contradiction again.

\item[(2)] Suppose that $m'=m$ and $l=k-1$.
Because $\cald_{Z,Z'}\neq\emptyset$, we must have $(Z,Z')\in\cald_{Z,Z'}$ by Proposition~\ref{0216},
and hence $b_{k-1}>c_k\geq b_k$ and $a_{k-1}>d_{k-1}\geq a_k$ by Lemma~\ref{0210}.
Because now $(Z,\Lambda_{\Psi'})\in\cald_{Z,Z'}$, by Lemma~\ref{0301} and Lemma~\ref{0307},
we have $a_{k-1}>c_k\geq a_k$ and $b_{k-1}>d_{k-1}\geq b_k$.
This means that every entry of $Z$ is either (a) greater than $d_{k-1}$, or (b) less than or equal to $c_k$.
The remaining argument is similar to that in case (1).

\item[(3)] Suppose that $m'=m+1$ and $l=k$.
The argument is similar to that in case (1).

\item[(4)] Suppose that $m'=m+1$ and $l=k-1$.
The argument is similar to that in case (2).
\end{enumerate}
For the four cases,
we conclude that $(\Lambda,\Lambda_N)\not\in\overline\calb^+_{Z,Z'}$ for any $\Lambda\in\overline\cals_Z$,
so $\Lambda_N$ does not occur in the relation $\overline\calb^+_{Z,Z'}$.
\end{proof}

\begin{lemma}\label{0808}
Suppose that $\cald_{Z,Z'}\neq\emptyset$.
Let $\Psi_0,\Psi'_0$ be the cores in $Z_\rmI,Z'_\rmI$ of $\cald_{Z,Z'}$ respectively.
\begin{enumerate}
\item[(i)] Suppose that $\binom{c_k}{d_l}\leq\Psi'_0$.
Then $(\Lambda,\Lambda')\in\overline\calb^+_{Z,Z'}$ if and only if
$(\Lambda,\Lambda'+\Lambda_{\binom{c_k}{d_l}})\in\overline\calb^+_{Z,Z'}$.

\item[(ii)] Suppose that $\binom{a_k}{b_l}\leq\Psi_0$.
Then $(\Lambda,\Lambda')\in\overline\calb^+_{Z,Z'}$ if and only if
$(\Lambda+\Lambda_{\binom{a_k}{b_l}},\Lambda')\in\overline\calb^+_{Z,Z'}$.
\end{enumerate}
\end{lemma}
\begin{proof}
Write
\begin{align*}
Z &=\binom{a_1,a_2,\ldots,a_{m+1}}{b_1,b_2,\ldots,b_m}, &
\Lambda &= \binom{a'_1,a'_2,\ldots,a'_{m_1}}{b'_1,b'_2,\ldots,b'_{m_2}};\\
Z' &=\binom{c_1,c_2,\ldots,c_{m'}}{d_1,d_2,\ldots,d_{m'}}, &
\Lambda' &= \binom{c_1',c'_2,\ldots,c'_{m_1'}}{d'_1,d'_2,\ldots,d'_{m'_2}}, \\
\Lambda'+\Lambda_{\binom{c_k}{d_l}} &=\binom{c''_1,c''_2,\ldots,c''_{m_1'}}{d''_1,d''_2,\ldots,d''_{m'_2}}.
\end{align*}
Suppose that $\Lambda'=\Lambda_N$ for some $N\subset Z'_\rmI$.
By Lemma~\ref{0809}, we know that either $\{c_k,d_l\}\subset N$ or $\{c_k,d_l\}\cap N=\emptyset$.
If $\{c_k,d_l\}\subset N$, then $\Lambda'+\Lambda_{\binom{c_k}{d_l}}=\Lambda_{N\smallsetminus\binom{c_k}{d_l}}$.
Therefore, without loss of generality, we may assume that $\{c_k,d_l\}\cap N=\emptyset$.
Then $c_k=c'_{k'}$, $d_l=d'_{l'}$ for some $k',l'$.

\begin{enumerate}
\item[(1)] Suppose that $m'=m$ and $l=k$.
Now $\binom{c_k}{d_l}\leq\Psi'_0$ means that $(Z,\Lambda_{\binom{c_k}{d_l}})\in\cald_{Z,Z'}$.
By Lemma~\ref{0210} and Lemma~\ref{0301}, we have
\begin{equation}\label{0807}
a_k>c_k\geq b_k\quad\text{ and }\quad a_k>d_k\geq b_k.
\end{equation}
Because $\binom{c_k}{d_l}=\binom{c'_{k'}}{d'_{l'}}$ is a consecutive pair,
the symbol $\Lambda'+\Lambda_{\binom{c_k}{d_l}}$
is obtained from $\Lambda'$ by switching the positions of entries $c'_{k'},d'_{l'}$, i.e.,
$c''_{k'}=d'_{l'}$, $c''_i=c'_i$ if $i\neq k'$; and $d''_{l'}=c'_{k'}$, $d''_i=d'_i$ if $i\neq l'$.
Note that there is no entry in $Z$ lying strictly between $a_k,b_k$,
so (\ref{0807}) implies that the condition:
\[
a'_{l'}>d'_{l'}=d_k\geq a'_{l'+1}\quad\text{ and }\quad b'_{k'-1}>c'_{k'}=c_k\geq b'_{k'}
\]
is equivalent to the condition:
\[
a'_{l'}>d''_{l'}=c_k\geq a'_{l'+1}\quad\text{ and }\quad b'_{k'-1}>c''_{k'}=d_k\geq b'_{k'}.
\]
Because any entry other than $c_k,d_l$ has the same position in $\Lambda'$ and in $\Lambda'+\Lambda_{\binom{c_k}{d_l}}$,
we conclude that $(\Lambda,\Lambda')\in\overline\calb^+_{Z,Z'}$ if and only if
$(\Lambda,\Lambda'+\Lambda_{\binom{c_k}{d_l}})\in\overline\calb^+_{Z,Z'}$.

\item[(2)] Suppose that $m'=m$ and $l=k-1$.
Then (\ref{0807}) is modified to
\[
b_{k-1}>c_k\geq a_k\quad\text{ and }\quad b_{k-1}>d_{k-1}\geq a_k.
\]
The remaining proof is similar to that in (1).

\item[(3)] Suppose that $m'=m+1$ and $l=k$.
Then (\ref{0807}) is modified to
\[
a_k\geq c_k>b_k\quad\text{ and }\quad a_k\geq d_k>b_k.
\]
The remaining proof is similar to that in (1).

\item[(4)] Suppose that $m'=m+1$ and $l=k-1$.
Then (\ref{0807}) is modified to
\[
b_{k-1}\geq c_k>a_k\quad\text{ and }\quad b_{k-1}\geq d_{k-1}>a_k.
\]
The remaining proof is similar to that in (1).
\end{enumerate}
The proof of (ii) is similar to that of (i).
\end{proof}

\begin{proposition}\label{0810}
Let $(\bfG,\bfG')=(\Sp_{2n},\rmO^\epsilon_{2n'})$,
and let $Z,Z'$ be special symbols of ranks $n,n'$ and defects $1,0$ respectively.
Suppose that $\cald_{Z,Z'}\neq\emptyset$.
\begin{enumerate}
\item[(i)] If $\Lambda'_0\in B^+_\Lambda$ for $\Lambda\in\overline\cals_Z$,
then $B^+_\Lambda=\{\,\Lambda'_0+\Lambda'\mid\Lambda'\in D_Z\,\}$.

\item[(ii)] If $\Lambda_0\in B^+_{\Lambda'}$ for $\Lambda'\in\overline\cals_{Z'}$,
then $B^+_{\Lambda'}=\{\,\Lambda_0+\Lambda\mid\Lambda\in D_{Z'}\,\}$.
\end{enumerate}
\end{proposition}
\begin{proof}
Recall that there exists a subset of (disjoint) consecutive pairs $\Psi'_0$ in $Z'_\rmI$ such that
$D_Z=\{\,\Lambda_{\Psi'}\mid\Psi'\leq\Psi'_0\,\}$ by Proposition~\ref{0307}.
Suppose that $\Lambda'_0\in B^+_\Lambda$.
Then Lemma~\ref{0808} implies that
$\{\,\Lambda'_0+\Lambda'\mid\Lambda'\in D_Z\,\}\subset B^+_\Lambda$.
Note that $\Lambda'_0+\Lambda'=\Lambda'_0+\Lambda''$ implies that $\Lambda'=\Lambda''$.
Then Lemma~\ref{0616} implies that in fact $\{\,\Lambda'_0+\Lambda'\mid\Lambda'\in D_Z\,\}=B^+_\Lambda$.

The proof of (ii) is similar.
\end{proof}

Now let $\Psi_0,\Psi'_0$ be the cores of $\cald_{Z,Z'}$ in $Z_\rmI,Z'_\rmI$ respectively.
Suppose that $\Lambda_{N}\in B^+_\Lambda$ for some $N\subset Z'_\rmI$.
Let $\binom{c_k}{d_l}$ be a consecutive pair in $\Psi'_0$.
By Lemma~\ref{0809}, we know that either $\{c_k,d_l\}\subset N$ or $N\cap\{c_k,d_l\}=\emptyset$.
Therefore by Proposition~\ref{0810},
there is a unique $\Lambda'\in\cals_{Z'}^{+,\Psi_0'}$ such that $\Lambda'\in B^+_\Lambda$, i.e.,
$|B^+_\Lambda\cap\cals_{Z'}^{+,\Psi'_0}|=1$ if $B^+_\Lambda\neq\emptyset$
where $\cals_{Z'}^{+,\Psi_0'}$ is defined in Subsection~\ref{0609}.
For $\Psi\leq\Psi_0$ and $\Psi'\leq\Psi'_0$,
we define
\begin{align}\label{0812}
\begin{split}
\calb^{+,\Psi,\Psi'}_{Z,Z'} &=\calb^+_{Z,Z'}\cap(\cals_Z^{\Psi}\times\cals_{Z'}^{+,\Psi'}),\\
\cald^{\Psi,\Psi'}_{Z,Z'} &=\cald_{Z,Z'}\cap(\cals_Z^{\Psi}\times\cals_{Z'}^{+,\Psi'}).
\end{split}
\end{align}
Then we have
\begin{equation}\label{0811}
\calb^+_{Z,Z'}=\{\,(\Lambda+\Lambda_M,\Lambda'+\Lambda_N)\mid(\Lambda,\Lambda')\in\calb^{+,\Psi,\Psi'}_{Z,Z'},
\ M\leq\Psi,\ N\leq\Psi'\,\}
\end{equation}
If $\Psi=\Psi_0$ and $\Psi'=\Psi'_0$,
the $\calb^{+,\Psi,\Psi'}_{Z,Z'}$ (resp.~$\cald^{\Psi,\Psi'}_{Z,Z'}$) is denoted by $\calb^{+,\natural}_{Z,Z'}$ (resp.~$\cald^\natural_{Z,Z'}$).
Moreover, from above explanation we see that $\calb^{+,\natural}_{Z,Z'}$ and $\cald^\natural_{Z,Z'}$,
as sub-relations of $\calb^+_{Z,Z'}$, are one-to-one.

\subsection{The relation for regular special symbols}

\begin{lemma}\label{0702}
Let $Z,Z'$ be given as in (\ref{0311}).
\begin{enumerate}
\item[(i)] If $m'=m+1$, $Z'$ regular and $D_{Z}=\{Z'\}$,
then $c_i>a_i$ and $d_i>b_i$ for each $i$.

\item[(ii)] If $m'=m$, $Z$ regular and $D_{Z'}=\{Z\}$,
then $a_i>c_i$ and $b_i>d_i$ for each $i$.
\end{enumerate}
\end{lemma}
\begin{proof}
First suppose that $m'=m+1$, $Z'$ regular and $D_{Z}=\{Z'\}$.
Because $Z'$ is regular, i.e., $Z'=Z'_\rmI$,
any $\binom{c_i}{d_i}$ or $\binom{c_i}{d_{i-1}}$ is a consecutive pair for each $i$.
If $a_{m+1}\geq c_{m+1}$, then $(Z,\Lambda_{\binom{c_{m+1}}{d_{m+1}}})\in\cald_{Z,Z'}$ by Lemma~\ref{0301}.
Then we get a contradiction because now $D_Z=\{Z'\}$.
Hence we conclude that $c_{m+1}>a_{m+1}$.
Then if $b_m\geq d_m$, then $(Z,\Lambda_{\binom{c_{m+1}}{d_m}})\in\cald_{Z,Z'}$ by Lemma~\ref{0301}, again.
Hence we conclude that $d_m>b_m$.
Then (i) is proved by induction.

The proof for (ii) is similar.
\end{proof}

\begin{lemma}\label{0804}
Let $Z,Z'$ be given as in (\ref{0311}).
\begin{enumerate}
\item[(i)] If $m'=m+1$, $Z'$ is regular and $D_Z=\{Z'\}$, then $Z$ is regular.

\item[(ii)] If $m'=m$, $Z$ is regular and $D_{Z'}=\{Z\}$, then $Z'$ is regular.
\end{enumerate}
\end{lemma}
\begin{proof}
Suppose $m'=m+1$, $Z'$ is regular and $D_Z=\{Z'\}$.
By Lemma~\ref{0702}, we have $c_i>a_i$ and $d_i>b_i$ for each $i$.
Now if $a_k=b_k$ for some $k$, then $d_k>a_k$;
if $a_{k+1}=b_k$ for some $k$, then $c_{k+1}>b_k$,
and both lead to a contradiction by Lemma~\ref{0210}
because now $(Z,Z')\in\cald_{Z,Z'}$.

The proof of (ii) is similar.
\end{proof}

\begin{lemma}\label{0805}
Let $Z,Z'$ be given as in (\ref{0311}).
\begin{enumerate}
\item[(i)] If $m'=m+1$, $Z'$ is regular and $D_Z=\{Z'\}$, then $D_{Z'}=\{Z\}$.

\item[(ii)] If $m'=m$, $Z$ is regular and $D_{Z'}=\{Z\}$, then $D_Z=\{Z'\}$.
\end{enumerate}
\end{lemma}
\begin{proof}
Suppose $m'=m+1$, $Z'$ is regular and $D_Z=\{Z'\}$.
By Lemma~\ref{0702}, we have $c_i>a_i$ and $d_i>b_i$ for each $i$.
Suppose that $\Lambda_{\binom{a_k}{b_l}}\in D_{Z'}$ for some consecutive pair $\binom{a_k}{b_l}$ in $Z_\rmI$.
Hence, we know that $l=k,k-1$.
By Lemma~\ref{0701}, we have $b_k\geq d_k$ if $l=k$; and $a_k\geq c_k$ if $l=k-1$,
and both lead to contradictions.

The proof of (ii) is similar.
\end{proof}

\begin{lemma}\label{0901}
Suppose that $m'=m+1$, $Z'$ are regular, and $D_{Z'}=\{Z\}$.
Let $M\subsetneq Z_\rmI$ such that $\Lambda_M$ occurs in the relation $\overline\calb^+_{Z,Z'}$.
Suppose $x$ is an entry in $Z_\rmI\smallsetminus M$ and $x$ is greater than any entry of $M$.
Then $\Lambda_{M\cup\{x\}}$ occurs in the relation $\overline\calb^+_{Z,Z'}$.
\end{lemma}
\begin{proof}
Suppose that $(\Lambda_M,\Lambda_N)\in\overline\calb^+_{Z,Z'}$ for some $N\subset Z'_\rmI$ and write
\[
\Lambda_M=\binom{a_1,a_2,\ldots,a_{m_1}}{b_1,b_2,\ldots,b_{m_2}},\qquad
\Lambda_N=\binom{c_1,c_2,\ldots,c_{m_1'}}{d_1,d_2,\ldots,d_{m'_2}}.
\]
Note that any entry in $\Lambda_M$ that is greater than all the entries in $M$
is in its natural position.
Note that now $m'=m+1$.
\begin{enumerate}
\item[(1)] Suppose that $x=a_k$ for some $k$.
So now every entry that is greater than $a_k$ is in its natural position in $\Lambda_M$.
Let $M'=M\cup\binom{a_k}{-}$ and write
\begin{equation}\label{0710}
\Lambda_{M'}=\binom{a'_1,\ldots,a'_{m_1-1}}{b'_1,\ldots,b'_{m_2+1}},\qquad
\Lambda_{N'}=\binom{c'_1,\ldots,c'_{m'_1+1}}{d'_1,\ldots,d'_{m_2'-1}}
\end{equation}
where $N'$ will be specified as follows:
\begin{enumerate}
\item
Suppose that $a_k\geq c_k$.
If $d_{k-1}>b_{k-1}$, we see that $\Lambda_{\binom{a_k}{b_{k-1}}}\in D_{Z'}$ by Lemma~\ref{0701}
and get a contradiction,
so we must have $b_{k-1}\geq d_{k-1}$.
Since $(\Lambda_M,\Lambda_N)\in\overline\calb^+_{Z,Z'}$,
by Lemma~\ref{0210}, we have
$c_{k-1}>b_{k-1}\geq d_{k-1}>a_k\geq c_k$, in particular, $d_{k-1}\in Z'_\rmI$.
Define $N'=(N\cup\binom{-}{d_{k-1}})\smallsetminus(N\cap\binom{-}{d_{k-1}})$.
Then $\Lambda_{N'}$ is obtained from $\Lambda_N$ by moving $d_{k-1}$ to the first row, i.e.,
\[
\begin{tabular}{c|c|c|c|c|c|c|c}
\toprule
$i$ & $1$ & $\ldots$ & $k-2$ & $k-1$ & $k$ & $k+1$ & $\ldots$ \\
\midrule
$a'_i$ & $a_1$ & $\ldots$ & $a_{k-2}$ & $a_{k-1}$ & $a_{k+1}$ & $a_{k+2}$ & $\ldots$ \\
\midrule
$b'_i$ & $b_1$ & $\ldots$ & $b_{k-2}$ & $b_{k-1}$ & $a_k$ & $b_k$ & $\ldots$ \\
\midrule
$c'_i$ & $c_1$ & $\ldots$ & $c_{k-2}$ & $c_{k-1}$ & $d_{k-1}$ & $c_k$ & $\ldots$ \\
\midrule
$d'_i$ & $d_1$ & $\ldots$ & $d_{k-2}$ & $d_k$ & $d_{k+1}$ & $d_{k+2}$ & $\ldots$ \\
\bottomrule
\end{tabular}
\]
and then
\begin{itemize}
\item $c'_k=d_{k-1}>a_k=b'_k$;

\item $b'_{k-1}=b_{k-1}\geq d_{k-1}=c'_k$ and $b'_k=a_k\geq c_k=c'_{k+1}$.
\end{itemize}

\item
Suppose that $c_k>a_k$.
If $b_k\geq d_k$ or $d_k$ does not exist,
then $\Lambda_{\binom{a_k}{b_k}}\in D_{Z'}$ and we get a contradiction.
Thus we must have $d_k>b_k$.
By Lemma~\ref{0210}, we have
$c_k>a_k\geq d_k>b_k\geq c_{k+1}$, in particular, $d_k\in Z'_\rmI$.
Define $N'=(N\cup\binom{-}{d_k})\smallsetminus(N\cap\binom{-}{d_k})$.
Then we have
\[
\begin{tabular}{c|c|c|c|c|c|c|c}
\toprule
$i$ & $1$ & $\ldots$ & $k-1$ & $k$ & $k+1$ & $k+2$ & $\ldots$ \\
\midrule
$a'_i$ & $a_1$ & $\ldots$ & $a_{k-1}$ & $a_{k+1}$ & $a_{k+2}$ & $a_{k+3}$ & $\ldots$ \\
\midrule
$b'_i$ & $b_1$ & $\ldots$ & $b_{k-1}$ & $a_k$ & $b_k$ & $b_{k+1}$ & $\ldots$ \\
\midrule
$c'_i$ & $c_1$ & $\ldots$ & $c_{k-1}$ & $c_k$ & $d_k$ & $c_{k+1}$ & $\ldots$ \\
\midrule
$d'_i$ & $d_1$ & $\ldots$ & $d_{k-1}$ & $d_{k+1}$ & $d_{k+2}$ & $d_{k+3}$ & $\ldots$ \\
\bottomrule
\end{tabular}
\]
and
\begin{itemize}
\item $c'_k=c_k>a_k=b'_k$, $c'_{k+1}=d_k>b_k=b'_{k+1}$;

\item $b'_k=a_k\geq d_k=c'_{k+1}$.
\end{itemize}
\end{enumerate}

\item[(2)] Suppose that $x=b_k$ for some $k$.
Let $M'=M\cup\binom{-}{b_k}$ and write $\Lambda_{M'}$ and
$\Lambda_{N'}$ as in (\ref{0710}) where $N'$ is defined as follows.
\begin{enumerate}
\item
Suppose that $b_k\geq d_k$.
If $c_k>a_k$,
then $\Lambda_{\binom{a_k}{b_k}}\in D_{Z'}$ and we get a contradiction,
and hence we must have $a_k\geq c_k$.
By Lemma~\ref{0210}, we have $d_{k-1}>a_k\geq c_k>b_k\geq d_k$, in particular, $c_k\in Z'_\rmI$.
Define $N'=(N\cup\binom{c_k}{-})\smallsetminus(N\cap\binom{c_k}{-})$.
Then we have
\[
\begin{tabular}{c|c|c|c|c|c|c|c}
\toprule
$i$ & $1$ & $\ldots$ & $k-1$ & $k$ & $k+1$ & $k+2$ & $\ldots$ \\
\midrule
$a'_i$ & $a_1$ & $\ldots$ & $a_{k-1}$ & $a_k$ & $b_k$ & $a_{k+1}$ & $\ldots$ \\
\midrule
$b'_i$ & $b_1$ & $\ldots$ & $b_{k-1}$ & $b_{k+1}$ & $b_{k+2}$ & $b_{k+3}$ & $\ldots$ \\
\midrule
$c'_i$ & $c_1$ & $\ldots$ & $c_{k-1}$ & $c_{k+1}$ & $c_{k+2}$ & $c_{k+3}$ & $\ldots$ \\
\midrule
$d'_i$ & $d_1$ & $\ldots$ & $d_{k-1}$ & $c_k$ & $d_k$ & $d_{k+1}$ & $\ldots$ \\
\bottomrule
\end{tabular}
\]
and
\begin{itemize}
\item $a'_k=a_k\geq c_k=d'_k$, $a'_{k+1}=b_k\geq d_k=d'_{k+1}$;

\item $d'_k=c_k>b_k=a'_{k+1}$.
\end{itemize}

\item
Suppose that $d_k>b_k$.
If $a_{k+1}\geq c_{k+1}$ or $c_{k+1}$ does not exist,
then $\Lambda_{\binom{a_{k+1}}{b_k}}\in D_{Z'}$ and we get a contradiction,
and hence we must have $c_{k+1}>a_{k+1}$.
Then $d_k>b_k\geq c_{k+1}>a_{k+1}\geq d_{k+1}$ by Lemma~\ref{0210}, in particular, $c_{k+1}\in Z'_\rmI$.
Define $N'=(N\cup\binom{c_{k+1}}{-})\smallsetminus(N\cap\binom{c_{k+1}}{-})$.
Then we have
\[
\begin{tabular}{c|c|c|c|c|c|c|c}
\toprule
$i$ & $1$ & $\ldots$ & $k-1$ & $k$ & $k+1$ & $k+2$ & $\ldots$ \\
\midrule
$a'_i$ & $a_1$ & $\ldots$ & $a_{k-1}$ & $a_k$ & $b_k$ & $a_{k+1}$ & $\ldots$ \\
\midrule
$b'_i$ & $b_1$ & $\ldots$ & $b_{k-1}$ & $b_{k+1}$ & $b_{k+2}$ & $b_{k+3}$ & $\ldots$ \\
\midrule
$c'_i$ & $c_1$ & $\ldots$ & $c_{k-1}$ & $c_k$ & $c_{k+2}$ & $c_{k+3}$ & $\ldots$ \\
\midrule
$d'_i$ & $d_1$ & $\ldots$ & $d_{k-1}$ & $d_k$ & $c_{k+1}$ & $d_{k+1}$ & $\ldots$ \\
\bottomrule
\end{tabular}
\]
and
\begin{itemize}
\item $a'_{k+1}=b_k\geq c_{k+1}=d'_{k+1}$;

\item $d'_k=d_k>b_k=a'_{k+1}$, $d'_{k+1}=c_{k+1}>a_{k+1}=a'_{k+2}$.
\end{itemize}
\end{enumerate}
\end{enumerate}
The above inequalities are all of which need to be checked.
We see that $(\Lambda_{M'},\Lambda_{N'})\in\overline\calb^+_{Z,Z'}$ for all four cases.
\end{proof}

\begin{lemma}
Suppose that $m'=m$, $Z$ is regular, and $D_Z=\{Z'\}$.
Let $N\subsetneq Z'_\rmI$ such that $\Lambda_N$ occurs in the relation $\overline\calb^+_{Z,Z'}$.
Suppose $x$ is an entry in $Z'_\rmI\smallsetminus N$ and $x$ is greater than any entry of $N$.
Then $\Lambda_{N\cup\{x\}}$ occurs in the relation $\overline\calb^+_{Z,Z'}$.
\end{lemma}
\begin{proof}
The proof is similar to that of Lemma~\ref{0901}.
\end{proof}

\begin{corollary}\label{0709}
Let $Z,Z'$ be special symbols of defects $1,0$ respectively.
\begin{enumerate}
\item[(i)] If $m'=m+1$, both $Z,Z'$ are regular,
and $D_{Z'}=\{Z\}$, then every symbol in $\overline\cals_Z$ occurs in the relation $\overline\calb^+_{Z,Z'}$.

\item[(ii)] If $m'=m$, both $Z,Z'$ are regular,
and $D_Z=\{Z'\}$, then every symbol in $\overline\cals_{Z'}$ occurs in the relation $\overline\calb^+_{Z,Z'}$.
\end{enumerate}
\end{corollary}
\begin{proof}
Suppose that $m'=m+1$, both $Z,Z'$ are regular, and $D_{Z'}=\{Z\}$.
Let $\Lambda_M$ be a symbol in $\overline\cals_Z$ for some subset $M\subset Z_\rmI$.
Define a chain of subsets
\[
\emptyset=M_0\subset M_1\subset\cdots\subset M_t=M
\]
where $M_0=\emptyset$ and, for $i\geq 1$, $M_i$ is obtained from $M_{i-1}$ by adding the smallest entry in $M\smallsetminus M_{i-1}$.
Now $\Lambda_{M_0}=\Lambda_\emptyset=Z$ occurs in the relation $\overline\calb^+_{Z,Z'}$.
If $\Lambda_{M_{i-1}}$ occurs in the relation,
then Lemma~\ref{0901} implies that $\Lambda_{M_i}$ also occurs in the relation.
Hence (i) is proved by induction.

The proof of (ii) is similar.
\end{proof}

\subsection{An isometry of inner product spaces}\label{0707}

Suppose that both $Z$ and $Z'$ are regular and write $Z,Z'$ as in (\ref{0311}).

\begin{enumerate}
\item[(1)] Suppose that $m'=m+1$.
Define
\begin{align*}
\theta\colon\{a_1,\ldots,a_{m+1}\}\sqcup\{b_1,\ldots,b_m\} &\rightarrow \{c_1,\ldots,c_{m+1}\}\sqcup\{d_1,\ldots, d_{m+1}\} \\
a_i &\mapsto d_i \\
b_i &\mapsto c_{i+1}
\end{align*}
for each $i$.
Note that $c_1$ is not in the image of $\theta$.
For $M\subset Z_\rmI$,
define $\theta(M)=\{\,\theta(x)\mid x\in M\,\}\subset Z'_\rmI$.
Then $\theta$ induces a mapping $\theta^+\colon\overline\cals_Z\rightarrow\overline\cals_{Z'}$ by
$\Lambda_M\mapsto\Lambda_{\theta(M)}$.
Moreover, if $\Lambda\in\cals_Z$,
then $\theta^+(\Lambda)\in\cals_{Z'}^+$.

\item[(2)] Suppose that $m'=m$.
Define
\begin{align*}
\theta\colon\{c_1,\ldots,c_m\}\sqcup\{d_1,\ldots,d_m\} &\rightarrow \{a_1,\ldots,a_{m+1}\}\sqcup\{b_1,\ldots, b_m\} \\
c_i &\mapsto b_i \\
d_i &\mapsto a_{i+1}
\end{align*}
for each $i$.
Note that $a_1$ is not in the image of $\theta$.
For $N\subset Z'_\rmI$,
define $\theta(N)=\{\,\theta(x)\mid x\in N\,\}\subset Z_\rmI$.
Then $\theta$ induces a mapping $\theta^+\colon\overline\cals_{Z'}\rightarrow\overline\cals_Z$ by
$\Lambda_N\mapsto\Lambda_{\theta(N)}$.
Moreover, if $\Lambda'\in\cals_{Z'}^+$,
then $\theta^+(\Lambda')\in\cals_Z$.
\end{enumerate}

\begin{remark}
Clearly, the maps $\theta^+$ defined in (\ref{0502}) and (\ref{0504})
are special cases of the maps $\theta^+$ in above (1) and (2) respectively.
\end{remark}

\begin{lemma}\label{0902}
Let $Z,Z'$ be given as in (\ref{0311}).
Suppose that both $Z$ and $Z'$ are regular,
and $\cald_{Z,Z'}$ is one-to-one.
\begin{enumerate}
\item[(i)] If $m'=m+1$, then $(\Lambda,\Lambda')\in\overline\calb^+_{Z,Z'}$ if and only if $\Lambda'=\theta^+(\Lambda)$.

\item[(ii)] If $m'=m$, then $(\Lambda,\Lambda')\in\overline\calb^+_{Z,Z'}$ if and only if $\Lambda=\theta^+(\Lambda')$.
\end{enumerate}
\end{lemma}
\begin{proof}
First suppose that $m'=m+1$.
By Corollary~\ref{0709} and Proposition~\ref{0810} we know that every $\Lambda\in\overline\cals_Z$ occurs in the relation
$\overline\calb^+_{Z,Z'}$ and $|B_\Lambda^+|=1$.
So it suffices to show that $\theta^+(\Lambda)\in B^+_\Lambda$.
Suppose that $\Lambda=\Lambda_M$ for some $M\subset Z_\rmI$.
We now want to prove the result by induction on $|M|$.
We know that $(Z,Z')$ is in $\overline\calb^+_{Z,Z'}$ by Proposition~\ref{0216}.
Hence the result is true if $|M|=0$ because $\theta^+(\Lambda_\emptyset)=\theta^+(Z)=Z'$.
Suppose that the result is true for $M$, i.e.,
$(\Lambda_M,\Lambda_{\theta(M)})\in\overline\calb^+_{Z,Z'}$.
If $a_k$ or $b_k$ is an entry in $Z_\rmI$ that is greater than any entry of $M$,
then from the proof of Lemma~\ref{0901} we know that
\[
(\Lambda_{\binom{a_k}{-}\cup M},\Lambda_{\binom{-}{d_k}\cup\theta(M)}),
(\Lambda_{\binom{-}{b_k}\cup M},\Lambda_{\binom{c_{k+1}}{-}\cup\theta(M)})\in\overline\calb^+_{Z,Z'}.
\]
It is clear that $\binom{-}{d_k}\cup\theta(M)=\theta(\binom{a_k}{-}\cup M)$ and
$\binom{c_{k+1}}{-}\cup\theta(M)=\theta(\binom{-}{b_k}\cup M)$.
Hence the lemma is proved for the case that $m'=m+1$.

The proof for the case that $m'=m$ is similar.
\end{proof}

Now we can show that the statement of Theorem~\ref{0310} holds when
$\epsilon=+$, both $Z,Z'$ are regular, and $\cald_{Z,Z'}$ is one-to-one.

\begin{proposition}\label{0806}
Let $(\bfG,\bfG')=(\Sp_{2n},\rmO^+_{2n'})$,
$Z,Z'$ special symbols of rank $n,n'$ and defect $1,0$ respectively.
Suppose that both $Z,Z'$ are regular, and $\cald_{Z,Z'}$ is one-to-one.
Then
\[
\frac{1}{2}\sum_{(\Sigma,\Sigma')\in\cald_{Z,Z'}} R^\bfG_\Sigma\otimes R^{\bfG'}_{\Sigma'}
=\sum_{(\Lambda,\Lambda')\in\calb^+_{Z,Z'}}\rho_\Lambda^\sharp\otimes\rho^\sharp_{\Lambda'}.
\]
\end{proposition}
\begin{proof}
Suppose that $Z$ and $Z'$ are regular, $D_Z=\{Z'\}$, and $D_{Z'}=\{Z\}$.
Write $Z,Z'$ as in (\ref{0311}) where $m'=m,m+1$, and let
\begin{equation}\label{0711}
Z_{(m)}=\binom{2m,2m-2,\ldots,0}{2m-1,2m-3,\ldots,1},\qquad
Z'_{(m')}=\binom{2m'-1,2m'-3,\ldots,1}{2m'-2,2m'-4,\ldots,0}.
\end{equation}
Then $Z_{(m)}$ (resp.~$Z'_{(m')}$) is a regular special symbol of rank $2m(m+1)$ (resp.~$2m'^2$) and
defect $1$ (resp.~$0$).

We define
\begin{align*}
h\colon\{a_1,\ldots,a_{m+1},b_1,\ldots,b_m\} &\rightarrow\{2m,2m-1,\ldots,1,0\} \\
a_i &\mapsto 2m+2-2i \\
b_i &\mapsto 2m+1-2i
\end{align*}
for each $i$.
Then the mapping $M\mapsto h(M)$ where $h(M)=\{\,h(x)\mid x\in M\,\}$ is bijective from the set of subsets of $Z$ to the
set of subsets of $Z_{(m)}$,
and hence induces a bijection $\bar h\colon\overline\cals_Z\rightarrow\overline\cals_{Z_{(m)}}$ by $\bar h(\Lambda_M)=\Lambda_{h(M)}$
for $M\subset Z_\rmI$.
Clearly, $\bar h$ maps $\cals_Z$ onto $\cals_{Z_{(m)}}$ and maps $\cals_{Z,1}$ onto $\cals_{Z_{(m)},1}$.
It is clearly that
$|M_1\cap M_2|=|h(M_1)\cap h(M_2)|$ for any $M_1,M_2\subset Z_\rmI$ and hence
\begin{equation}\label{0703}
\langle R_\Sigma,\rho_\Lambda\rangle=\langle R_{\bar h(\Sigma)},\rho_{\bar h(\Lambda)}\rangle
\end{equation}
for any $\Sigma\in\cals_{Z,1}$ and $\Lambda\in\cals_Z$ by Proposition~\ref{0306}.

Similarly, the bijection
\begin{align*}
h'\colon\{c_1,\ldots,c_{m'},d_1,\ldots,d_{m'}\} &\rightarrow\{2m'-1,2m'-2,\ldots,1,0\} \\
c_i &\mapsto 2m'+1-2i \\
d_i &\mapsto 2m'-2i
\end{align*}
for each $i$ induces a bijection
$\bar h'\colon\overline\cals_{Z'}\rightarrow\overline\cals_{Z'_{(m')}}$ by $\bar h'(\Lambda_N)=\Lambda_{h'(N)}$
for $N\subset Z'_\rmI$.
Clearly, $\bar h'$ maps $\cals^+_{Z'}$ onto $\cals^+_{Z'_{(m')}}$
and maps $\cals_{Z',0}\rightarrow\cals_{Z'_{(m')},0}$
such that
\begin{equation}\label{0704}
\langle R_{\Sigma'},\rho_{\Lambda'}\rangle=\langle R_{\bar h'(\Sigma')},\rho_{\bar h'(\Lambda')}\rangle
\end{equation}
for any $\Sigma'\in\cals_{Z',0}$ and $\Lambda'\in\cals^+_{Z'}$.

Then we see that
\begin{equation}\label{0705}
(\Lambda,\Lambda')\in\calb^+_{Z,Z'}\text{ \ if and only if \ }
(\bar h(\Lambda),\bar h'(\Lambda'))\in\calb^+_{Z_{(m)},Z'_{(m')}}
\end{equation}
by Example~\ref{0225}, Example~\ref{0226} and Lemma~\ref{0902}.

In terms of linear algebra,
(\ref{0342}) is to express the uniform projection of the vector
$\sum_{(\Lambda,\Lambda')\in\calb^+_{Z,Z'}}\rho_\Lambda\otimes\rho_{\Lambda'}\in\calv_Z\otimes\calv^+_{Z'}$
as a linear combination of the orthogonal basis
\[
\{\,R_\Sigma\otimes R_{\Sigma'}\mid\,\Sigma\in\cals_{Z,1},\Sigma'\in\bar\cals_{Z',0}\}
\]
of $\calv_Z^\sharp\otimes\calv^\sharp_{Z'}$ where $\bar\cals_{Z',0}$ is a complete set of representatives
of cosets $\{\Sigma',\Sigma'^\rmt\}$ in $\cals_{Z',0}$.
Now Proposition~\ref{0512} and Proposition~\ref{0513} mean that (\ref{0342}) holds if
$Z=Z_{(m)}$ and $Z'=Z'_{(m')}$ where $m'=m$ or $m+1$.
Now $(h,h')$ induces a vector space isomorphism
$\tilde h\otimes\tilde h'\colon\calv_Z\otimes\calv^+_{Z'}\rightarrow\calv_{Z_{(m)}}\otimes\calv^+_{Z'_{(m')}}$
such that $\tilde h(\rho_\Lambda)=\rho_{\bar h(\Lambda)}$, $\tilde h(R_\Sigma)=R_{\bar h(\Sigma)}$,
$\tilde h'(\rho_{\Lambda'})=\rho_{\bar h'(\Lambda')}$, and $\tilde h'(R_{\Sigma'})=R_{\bar h'(\Sigma')}$.
Therefore by (\ref{0703}), (\ref{0704}) and (\ref{0705}),
we have
\begin{align*}
\tilde h\otimes\tilde h'\Biggl(\Biggl(\sum_{(\Lambda,\Lambda')\in\calb^+_{Z,Z'}}\rho_\Lambda\otimes\rho_{\Lambda'}\Biggr)^\sharp\Biggr)
&=\Biggl(\tilde h\otimes\tilde h'\Biggl(\sum_{(\Lambda,\Lambda')\in\calb^+_{Z,Z'}}\rho_\Lambda\otimes\rho_{\Lambda'}\Biggr)\Biggr)^\sharp \\
&=\Biggl(\sum_{(\bar h(\Lambda),\bar h'(\Lambda'))\in\calb^+_{Z_{(m)},Z'_{(m')}}}\rho_{\bar h(\Lambda)}\otimes\rho_{\bar h'(\Lambda')}\Biggr)^\sharp \\
&=\frac{1}{2} \sum_{(\bar h(\Sigma),\bar h'(\Sigma'))\in\cald_{Z_{(m)},Z'_{(m')}}} R_{\bar h(\Sigma)}\otimes R_{\bar h'(\Sigma')} \\
&=\tilde h\otimes\tilde h'\Biggl(\frac{1}{2}\sum_{(\Sigma,\Sigma')\in\cald_{Z,Z'}} R_\Sigma\otimes R_{\Sigma'}\Biggr).
\end{align*}
Because now $\tilde h\otimes\tilde h'$ is a vector space isomorphism, we have
\[
\Biggl(\sum_{(\Lambda,\Lambda')\in\calb^+_{Z,Z'}}\rho_\Lambda\otimes\rho_{\Lambda'}\Biggr)^\sharp
=\frac{1}{2}\sum_{(\Sigma,\Sigma')\in\cald_{Z,Z'}} R_\Sigma\otimes R_{\Sigma'}.
\]
Thus the proposition is proved.
\end{proof}


\section{The Relation $\calb^+_{Z,Z'}$ II: General Case}\label{0802}
The purpose of this section (Proposition~\ref{0815}) is to show that the statement of
Theorem~\ref{0310} holds when $\epsilon=+$.
The strategy is to reduce the general case inductively via an isometry of inner product spaces
to the case that considered in the previous section.
Let
\begin{equation}\label{0914}
Z=\binom{a_1,a_2,\ldots,a_{m+1}}{b_1,b_2,\ldots,b_m},\qquad
Z'=\binom{c_1,c_2,\ldots,c_{m'}}{d_1,d_2,\ldots,d_{m'}}
\end{equation}
be special symbols of defect $1,0$ respectively such that $m'=m,m+1$.
We suppose that the relation $\cald_{Z,Z'}$ is nonempty.

\subsection{Derivative a special symbol}\label{0903}
We consider sets of pairs in $Z,Z'$ respectively according the following order
\begin{equation}\label{0814}
\textstyle
\{\binom{a_{m+1}}{b_m},\binom{c_{m'}}{d_{m'}}\},
\{\binom{a_m}{b_m},\binom{c_{m'}}{d_{m'-1}}\},
\{\binom{a_m}{b_{m-1}},\binom{c_{m'-1}}{d_{m'-1}}\},
\{\binom{a_{m-1}}{b_{m-1}},\binom{c_{m'-1}}{d_{m'-2}}\},\ldots.
\end{equation}
Suppose that $\{\binom{a_k}{b_l},\binom{c_{l'}}{d_{k'}}\}$ is the first set in the above order such that
at least one of the two pairs in the set is either a pair of doubles or in the core of the relation $\cald_{Z,Z'}$.
Such a set exists if $\cald_{Z,Z'}$ is not one-to-one or one of $Z,Z'$ is not regular.
It is easy to see that $l=k-1,k$; $l'=k',k'+1$, and 
\begin{equation}\label{0828}
\begin{cases}
k'=k-1\text{ and }l'=l, & \text{if $m'=m$};\\
k'=k\text{ and }l'=l+1, & \text{if $m'=m+1$}.
\end{cases}
\end{equation}

Now we want to construct new special symbols $Z^{(1)},Z'^{(1)}$ according to the following situations:
\begin{enumerate}
\item[(I)] If both $\binom{a_k}{b_l}$ and $\binom{c_{l'}}{d_{k'}}$ are pairs of doubles or both are in
the core of $\cald_{Z,Z'}$,
then we define
\begin{align*}
Z^{(1)} &=\binom{a_1-1,\ldots,a_{k-1}-1,a_{k+1},\ldots,a_{m+1}}{b_1-1,\ldots,b_{l-1}-1,b_{l+1},\ldots,b_{m}},\\
Z'^{(1)} &=\binom{c_1-1,\ldots,c_{l'-1}-1,c_{l'+1},\ldots,c_{m'}}{d_1-1,\ldots,d_{k'-1}-1,d_{k'+1},\ldots,d_{m'}}.
\end{align*}

\item[(II)] If $\binom{a_k}{b_l}$ is a pair of doubles or in the core,
but $\binom{c_{l'}}{d_{k'}}$ is neither,
then we define
\begin{align*}
Z^{(1)} &=\binom{a_1-1,\ldots,a_{k-1}-1,d_{k'+1},\ldots,d_{m'}}{b_1-1,\ldots,b_{l-1}-1,c_{l'+1},\ldots,c_{m'}},\\
Z'^{(1)} &=\binom{c_1,\ldots,c_{l'},b_{l+1},\ldots,b_m}{d_1,\ldots,d_{k'},a_{k+1},\ldots,a_{m+1}}.
\end{align*}

\item[(III)] If $\binom{c_{l'}}{d_{k'}}$ is a pair of doubles or in the core,
but $\binom{a_k}{b_l}$ is neither,
then we define
\begin{align*}
Z^{(1)} &=\binom{a_1,\ldots,a_k,d_{k'+1},\ldots,d_{m'}}{b_1,\ldots,b_l,c_{l'+1},\ldots,c_{m'}},\\
Z'^{(1)} &=\binom{c_1-1,\ldots,c_{l'-1}-1,b_{l+1},\ldots,b_m}{d_1-1,\ldots,d_{k'-1}-1,a_{k+1},\ldots,a_{m+1}}.
\end{align*}
\end{enumerate}
Note that if $\binom{a_k}{b_l}$ (resp.~$\binom{c_{l'}}{d_{k'}}$) is a pair of doubles or in the core,
then it is removed from $Z$ (resp.~$Z'$) to obtain the new symbol $Z^{(1)}$ (resp.~$Z'^{(1)}$).
The symbol $Z^{(1)}$ (resp.~$Z'^{(1)}$) can be regarded as the (first) ``\emph{derivative}'' of $Z$ (resp.~$Z'$)
and we can repeat the same process to obtain the second ``derivative'' $Z^{(2)}$ (resp.~$Z'^{(2)}$) and so on.

For the pairs $\binom{a_k}{b_l},\binom{c_{l'}}{d_{k'}}$ in sequence (\ref{0828}), 
we define their \emph{next pairs} $\binom{a_k}{b_l}^\flat,\binom{c_{l'}}{d_{k'}}^\flat$ 
as follows:
\begin{align*}
\textstyle\binom{a_k}{b_l}^\flat=\begin{cases}
\binom{a_{k-1}}{b_l}, & \text{if $\binom{a_k}{b_l}=\binom{a_k}{b_{k-1}}$};\\
\binom{a_k}{b_{l-1}}, & \text{if $\binom{a_k}{b_l}=\binom{a_k}{b_k}$},
\end{cases}\qquad\quad
\textstyle\binom{c_{l'}}{d_{k'}}^\flat=\begin{cases}
\binom{c_{l'-1}}{d_{k'}}, & \text{if $\binom{c_{l'}}{d_{k'}}=\binom{c_{k'+1}}{d_{k'}}$};\\
\binom{c_{l'}}{d_{k'-1}}, & \text{if $\binom{c_{l'}}{d_{k'}}=\binom{c_{k'}}{d_{k'}}$}.
\end{cases}
\end{align*}

\begin{lemma}\label{0827}
Keep the above notations and suppose that $\cald_{Z,Z'}\neq\emptyset$.
\begin{enumerate}
\item[(i)] If $m'=m+1$ and $\binom{a_k}{b_l}$ is either a pair of doubles or in the core,
then $\binom{c_{l'}}{d_{k'}}^\flat$ is neither.

\item[(ii)] If $m'=m$ and $\binom{c_{l'}}{d_{k'}}$ is either a pair of doubles or in the core,
then $\binom{a_k}{b_l}^\flat$ is neither.
\end{enumerate}
\end{lemma}
\begin{proof}
By Proposition~\ref{0216},
$\cald_{Z,Z'}\neq\emptyset$ implies that $(Z,Z')\in\cald_{Z,Z'}$.
We consider Case (i) first.
\begin{enumerate}
\item[(1)] Suppose that $m'=m+1$ and $l=k$.
Then $l'=k+1$ and $k'=k$ by (\ref{0828}), and we have $\binom{c_{l'}}{d_{k'}}^\flat=\binom{c_k}{d_k}$.

First suppose that $\binom{a_k}{b_k}$ is a pair of doubles, i.e., $a_k=b_k$.
If $\binom{c_k}{d_k}$ is also a pair of doubles,
then we have $a_k=b_k<c_k=d_k$ and $a_k\geq d_k$ by Lemma~\ref{0210}, and we get a contradiction.
If $\binom{c_k}{d_k}$ is in the core, then $d_k>b_k=a_k$ by Lemma~\ref{0301} and $a_k\geq d_k$ by Lemma~\ref{0210},
and we get a contradiction, again.

Next suppose $\binom{a_k}{b_k}$ is in the core.
If $\binom{c_k}{d_k}$ is a pair of doubles,
then we have $d_k=c_k>a_k$ by Lemma~\ref{0701} and $a_k\geq d_k$ by Lemma~\ref{0210},
and we get a contradiction.
If $\binom{c_k}{d_k}$ is also in the core,
then we have $a_k\geq c_k$ by Lemma~\ref{0301} and $c_k>a_k$ by Lemma~\ref{0701},
and we get a contradiction, again.

\item[(2)] Suppose that $m'=m+1$ and $l=k-1$.
Then $l'=k$ and $k'=k$ by (\ref{0828}),
and we have $\binom{c_{l'}}{d_{k'}}^\flat=\binom{c_k}{d_{k-1}}$.

First suppose that $\binom{a_k}{b_{k-1}}$ is a pair of doubles.
If $\binom{c_k}{d_{k-1}}$ is also a pair of doubles,
then we have $a_k=b_{k-1}\geq c_k=d_{k-1}$ and $d_{k-1}>a_k$ by Lemma~\ref{0210}, and we get a contradiction.
If $\binom{c_k}{d_{k-1}}$ is in the core, then $c_k>a_k=b_{k-1}$ by Lemma~\ref{0301} and $b_{k-1}\geq c_k$ by Lemma~\ref{0210},
and we get a contradiction, again.

Next suppose $\binom{a_k}{b_{k-1}}$ is in the core.
If $\binom{c_k}{d_{k-1}}$ is a pair of doubles,
then we have $d_{k-1}=c_k\leq a_k$ by Lemma~\ref{0701} and $d_{k-1}>a_k$ by Lemma~\ref{0210},
and we get a contradiction.
If $\binom{c_k}{d_{k-1}}$ is also in the core,
then we have $c_k>a_k$ by Lemma~\ref{0301} and $a_k\geq c_k$ by Lemma~\ref{0701},
and we get a contradiction, again.
\end{enumerate}
Therefore, (i) is proved.
The proof for (ii) is similar.
\end{proof}

\begin{lemma}\label{0905}
Keep the above notations.
\begin{enumerate}
\item[(i)] Case (II) happens only if $m'=m$.

\item[(ii)] Case (III) happens only if $m'=m+1$.
\end{enumerate}
\end{lemma}
\begin{proof}
Suppose that $m'=m+1$, $\binom{a_k}{b_l}$ is a pair of doubles or in the core,
but $\binom{c_{l'}}{d_{k'}}$ is neither.
Note that $l=k,k-1$.
\begin{enumerate}
\item[(1)] Suppose that $l=k$.
Consider the following two symbols
\[
X=\binom{a_k,a_{k+1},\ldots,a_{m+1}}{b_k,b_{k+1},\ldots,b_m},\qquad
X'=\binom{c_k,c_{k+1},\ldots,c_{m+1}}{d_k,d_{k+1},\ldots,d_{m+1}}.
\]
Then $X$ (resp.~$X'$) is a special symbol of size $(m+2-k,m+1-k)$ (resp.~$(m+2-k,m+2-k)$).
Because $(Z,Z')\in\cald_{Z,Z'}$, we see that $(X,X')\in\cald_{X,X'}$.

Now $\binom{a_k}{b_l}=\binom{a_k}{b_k}$ and $\binom{c_{l'}}{d_{k'}}=\binom{c_{k+1}}{d_k}$.
By Lemma~\ref{0827}, we see that $\binom{c_{l'}}{d_{k'}}^\flat=\binom{c_k}{d_k}$ is neither a pair of doubles and in the core,
so $X'$ is regular and $D_X=\{X'\}$ from our assumption on $\binom{a_k}{b_l}$, $\binom{c_{l'}}{d_{k'}}$.
Then by Lemma~\ref{0804} and Lemma~\ref{0805},
the symbol $X$ must be regular and $D_{X'}=\{X\}$.
But the assumption that $\binom{a_k}{b_l}$ is a pair of doubles
means that $X$ is not regular;
and the assumption that $\binom{a_k}{b_l}$ in the core means that
$D_{X'}\neq\{X\}$.
So we get a contradiction.

\item[(2)] Suppose $l=k-1$.
Let
\[
X=\binom{d_{k-1},d_k,\ldots,d_{m+1}}{c_k,c_{k+1},\ldots,c_{m+1}}\qquad
X'=\binom{b_{k-1},b_{k+1},\ldots,b_m}{a_k,a_{k+2},\ldots,a_{m+1}}.
\]
Then $X$ (resp.~$X'$) is a special symbol of size $(m+3-k,m+2-k)$ (resp.~$(m+2-k,m+2-k)$),
and $(X,X')\in\cald_{X,X'}$.

Now $\binom{a_k}{b_l}=\binom{a_k}{b_{k-1}}$, $\binom{c_{l'}}{d_{k'}}=\binom{c_k}{d_k}$.
By Lemma~\ref{0827}, we see that $\binom{c_{l'}}{d_{k'}}^\flat=\binom{c_k}{d_{k-1}}$ is neither a pair of doubles and in the core,
so $X$ is regular and $D_{X'}=\{X\}$ from our assumption on $\binom{a_k}{b_l}$, $\binom{c_{l'}}{d_{k'}}$.
Then by Lemma~\ref{0804} and Lemma~\ref{0805},
the symbol $X'$ is regular and $D_X=\{X'\}$.
We get a contradiction by the same argument in (1).
\end{enumerate}
Therefore, (i) is proved.

Suppose that $m'=m$, $\binom{c_{l'}}{d_{k'}}$ is a pair of doubles or in the core,
but $\binom{a_k}{b_l}$ is neither.
We can show that we will get a contradiction by the same argument for the case $m'=m+1$.
\end{proof}

\begin{lemma}\label{0821}
Keep the notations above.
Then
\[
{\rm size}(Z^{(1)})=\begin{cases}
(m,m-1) \\
(m,m-1) \\
(m+1,m)
\end{cases}\text{ and}\quad
{\rm size}(Z'^{(1)})=\begin{cases}
(m'-1,m'-1), & \text{for Case (I)};\\
(m',m'), & \text{for Case (II)};\\
(m'-1,m'-1), & \text{for Case (III)}.
\end{cases}
\]
In particular,
${\rm def}(Z^{(1)})={\rm def}(Z)=1$, ${\rm def}(Z'^{(1)})={\rm def}(Z')=0$
for all three cases.
\end{lemma}
\begin{proof}
Case (I) is trivial.

Now we consider Case (II).
From the choice of $\{\binom{a_k}{b_l},\binom{c_{l'}}{d_{k'}}\}$,
we see that
\begin{align*}
|\{a_{k+1},a_{k+2},\ldots,a_{m+1}\}| &=|\{d_{k'+1},d_{k'+2},\ldots,d_{m'}\}|,\\
|\{b_{l+1},b_{l+2},\ldots,b_m\}| &=|\{c_{l'+1},c_{l'+2},\ldots,c_{m'}\}|.
\end{align*}
Therefore, ${\rm size}(Z^{(1)})=(m,m-1)$ and ${\rm size}(Z'^{(1)})=(m',m')$.

The proof for Case (III) is similar to that of Case (II).
\end{proof}

\begin{remark}\label{0822}
Keep the notations above.
\begin{enumerate}
\item[(1)] For Case (I)
it is easy to see that
\begin{align*}
\deg(Z^{(1)}) &=\begin{cases}
\deg(Z), & \text{if $\binom{a_k}{b_l}$ is a pair of double};\\
\deg(Z)-1, & \text{if $\binom{a_k}{b_l}$ is in the core},
\end{cases} \\
\deg(Z'^{(1)}) &=\begin{cases}
\deg(Z'), & \text{if $\binom{c_{l'}}{d_{k'}}$ is a pair of double};\\
\deg(Z')-1, & \text{if $\binom{c_{l'}}{d_{k'}}$ is in the core}.
\end{cases}
\end{align*}

\item[(2)] For Case (II),
we have
\begin{align*}
\deg(Z^{(1)}) &=\begin{cases}
\deg(Z), & \text{if $\binom{a_k}{b_l}$ is a pair of double};\\
\deg(Z)-1, & \text{if $\binom{a_k}{b_l}$ is in the core},
\end{cases} \\
\deg(Z'^{(1)}) &= \deg(Z').
\end{align*}

\item[(3)] For Case (III),
we have
\begin{align*}
\deg(Z^{(1)}) &=\deg(Z);\\
\deg(Z'^{(1)}) &=\begin{cases}
\deg(Z'), & \text{if $\binom{c_{l'}}{d_{k'}}$ is a pair of double};\\
\deg(Z')-1, & \text{if $\binom{c_{l'}}{d_{k'}}$ is in the core}.
\end{cases}
\end{align*}
\end{enumerate}
\end{remark}

\begin{lemma}\label{0820}
Keep the notations above.
For all three cases, $Z^{(1)},Z'^{(1)}$ are special symbols.
\end{lemma}
\begin{proof}
We will only prove that $Z^{(1)}$ is a special symbol.
The proof for $Z'^{(1)}$ is similar.
Write $Z$ as in (\ref{0914}) and
\begin{equation}\label{0913}
Z^{(1)}=\binom{a_1',a'_2,\ldots,a'_{m_1}}{b_1',b'_2,\ldots,b'_{m_2}}.
\end{equation}

First we consider Case (I).
Then $m_1=m$, $m_2=m-1$, $l=k,k-1$ and
\[
Z^{(1)}=\binom{a_1-1,\ldots,a_{k-1}-1,a_{k+1},\ldots,a_{m+1}}{b_1-1,\ldots,b_{l-1}-1,b_{l+1},\ldots,b_{m}}.
\]
\begin{enumerate}
\item[(1)] Suppose that $l=k$.
We have
\begin{itemize}
\item $a'_i=a_i-1\geq b_i-1=b'_i$ if $i<k$; $a'_i=a_{i+1}\geq b_{i+1}=b'_i$ if $i\geq k$,

\item $b'_i=b_i-1\geq a_{i+1}-1=a'_{i+1}$ if $i<k-1$; $b'_{k-1}=b_{k-1}-1\geq b_k\geq a_{k+1}=a'_k$;
$b'_i=b_{i+1}\geq a_{i+2}=a'_{i+1}$ if $i\geq k$.
\end{itemize}

\item[(2)] Suppose that $l=k-1$.
We have
\begin{itemize}
\item $a'_i=a_i-1\geq b_i-1=b'_i$ if $i<k-1$; $a'_{k-1}=a_{k-1}-1\geq a_k\geq b_k=b'_{k-1}$; $a'_i=a_{i+1}\geq b_{i+1}=b'_i$ if $i\geq k$,

\item $b'_i=b_i-1\geq a_{i+1}-1=a'_{i+1}$ if $i<k-1$; $b'_i=b_{i+1}\geq a_{i+2}=a'_{i+1}$ if $i\geq k-1$.
\end{itemize}
\end{enumerate}

The proofs for Case (II) and (III) are similar.
\end{proof}

\subsection{Case (I)}\label{0826}
In this subsection,
we assume that we in the situation of Case (I) in the previous subsection.
Define a map $f\colon Z\smallsetminus\{a_k,b_l\}\rightarrow Z^{(1)}$ given by
\begin{align}\label{1001}
\begin{split}
a_i \mapsto\begin{cases}
a_i-1, & \text{if $1\leq i<k$;}\\
a_i, & \text{if $k<i\leq m+1$},
\end{cases}\qquad
b_i \mapsto\begin{cases}
b_i-1, & \text{if $1\leq i<l$;}\\
b_i, & \text{if $l<i\leq m$}.
\end{cases}
\end{split}
\end{align}

\begin{lemma}
Suppose that we are in Case (I) of Subsection~\ref{0903}.
The mapping $f$ given in (\ref{1001}) gives a bijection
from $Z_\rmI\smallsetminus\{a_k,b_l\}$ onto $Z^{(1)}_\rmI$.
\end{lemma}
\begin{proof}
Clearly $f$ gives a bijection between entries in $Z\smallsetminus\{a_k,b_l\}$
and entries in $Z^{(1)}$.
So we only need to check that $f$ maps singles to singles and non-singles to non-singles.

Write $Z^{(1)}$ as in (\ref{0913}).
Suppose that $\binom{a_i}{b_j}$ is a pair of doubles in $Z\smallsetminus\binom{a_k}{b_l}$.
From the requirement of $\binom{a_k}{b_l}$, we see that $i<k$ and $j<l$ and thus $f(a_i)=f(b_j)$.
Conversely, if $f(a_i)=f(b_j)$, then by the requirement for $\binom{a_k}{b_l}$
we must have $i<k$ and $j<l$, and thus we have $a_i-1=b_j-1$, i.e., $a_i=b_j$.
\qed
\end{proof}

Now we want to define a subspace $\calv^{(1)}_Z$ of $\calv_Z$ and an isomorphism
$\tilde f\colon\calv^{(1)}_Z\rightarrow\calv_{Z^{(1)}}$ according the following two cases:
\begin{enumerate}
\item[(1)] Suppose that $\binom{a_k}{b_l}$ is a pair of doubles.
Then $f$ induces a bijection $\bar f\colon\overline\cals_Z\rightarrow\overline\cals_{Z^{(1)}}$ by
$\bar f(\Lambda_M)=\Lambda_{f(M)}$ for $M\subset Z_\rmI=Z_\rmI\smallsetminus\{a_k,b_l\}$
and $f(M)=\{\,f(x)\mid x\in M\,\}$.
Clearly $\bar f(\cals_Z)=\cals_{Z^{(1)}}$.
Let $\calv^{(1)}_Z=\calv_Z$, $\rho^{(1)}_\Lambda=\rho_\Lambda$ and $R_\Sigma^{(1)}=R_\Sigma$
for $\Lambda\in\cals_Z$ and $\Sigma\in\cals_{Z,1}$.

\item[(2)] Suppose that $\binom{a_k}{b_l}$ is in the core of $\cald_{Z,Z'}$.
Now $f$ induces a bijection $\bar f\colon\overline\cals_Z^\Psi\rightarrow\overline\cals_{Z^{(1)}}$
by $\bar f(\Lambda_M)=\Lambda_{f(M)}$ for $M\subset Z_\rmI\smallsetminus\{a_k,b_l\}$
where $\Psi=\binom{a_k}{b_l}$.
Clearly $\bar f(\cals_Z^\Psi)=\cals_{Z^{(1)}}$.
Let $\calv^{(1)}_Z$ be the subspace of $\calv_Z$ spanned by vectors of the form
\begin{equation}\label{0930}
\rho_\Lambda^{(1)}:=\frac{1}{\sqrt 2}(\rho_\Lambda+\rho_{\Lambda^\natural})
\end{equation}
where $\Lambda=\Lambda_M$, $\Lambda^\natural=\Lambda_{M\cup\binom{a_k}{b_l}}$
for $M\subset Z_\rmI\smallsetminus\{a_k,b_l\}$.
It is clear that $\{\rho_\Lambda^{(1)}\mid\Lambda\in\cals_Z^\Psi\}$
forms an orthonormal basis for $\calv^{(1)}_Z$.
Moreover, by (\ref{0315}),
the uniform projection $(\calv^{(1)}_Z)^\sharp$ is spanned by vectors of the form
\[
R^{(1)}_\Sigma:=\frac{1}{\sqrt 2}(R_\Sigma+R_{\Sigma^\natural})
\]
where $\Sigma=\Lambda_N$ and $\Sigma^\natural=\Lambda_{N\cup\binom{a_k}{b_l}}$
for $N\subset Z_\rmI\smallsetminus\{a_k,b_l\}$ such that $|N^*|=|N_*|$.
\end{enumerate}
Then, for both (1) and (2), the linear transformation
$\tilde f\colon\calv^{(1)}_Z\rightarrow\calv_{Z^{(1)}}$ given by $\rho^{(1)}_{\Lambda}\mapsto
\rho_{\bar f(\Lambda)}$
is an inner product space isometry.

\begin{example}\label{0829}
Suppose that $Z=\binom{8,5,1}{6,3}$ and $Z'=\binom{8,6,2}{6,3,0}$.
Then $m=2$ and $m'=3$.
It is easy to see from the definition that $\calb^+_{Z,Z'}$ is the following:
\[
\begin{tabular}{c|cccc}
\toprule
& $\binom{8,6,2}{6,3,0}$ & $\binom{8,6,3}{6,2,0}$ & $\binom{6,2,0}{8,6,3}$ & $\binom{6,3,0}{8,6,2}$ \\
\midrule
$\binom{8,5,1}{6,3}$ & \checkmark & \checkmark \\
$\binom{8,3,1}{6,5}$ & \checkmark & \checkmark \\
$\binom{8,6,5}{3,1}$ & & & \checkmark & \checkmark \\
$\binom{8,6,3}{5,1}$ & & & \checkmark & \checkmark \\
\bottomrule
\end{tabular}
\]
where the top row (resp.~left column) are the symbols in $\cals^+_{Z'}$ (resp.~$\cals_Z$) which occur in the relation $\calb^+_{Z,Z'}$.
And it means that $(\Lambda,\Lambda')\in\calb^+_{Z,Z'}$ if there is a ``\checkmark'' in the position of row $\Lambda$
and column $\Lambda'$.

It is easy to see that the pairs $\binom{a_k}{b_l}=\binom{a_2}{b_2}=\binom{5}{3}$ and
$\binom{c_{l'}}{d_{k'}}=\binom{c_3}{d_2}=\binom{2}{3}$ are in the cores.
So we are in Case (I) and then $Z^{(1)}=\binom{7,1}{5}$ and $Z'^{(1)}=\binom{7,5}{5,0}$.
Then $\calb^+_{Z^{(1)},Z'^{(1)}}$ is one-to-one:
\[
\begin{tabular}{c|cc}
\toprule
& $\binom{7,5}{5,0}$ & $\binom{5,0}{7,5}$ \\
\midrule
$\binom{7,1}{5}$ & \checkmark & \\
$\binom{7,5}{1}$ & & \checkmark \\
\bottomrule
\end{tabular}
\]
but $Z'^{(1)}$ is not regular.
Now 
\[
\rho^{(1)}_{\binom{8,5,1}{6,3}}=\tfrac{1}{\sqrt 2}(\rho_{\binom{8,5,1}{6,3}}+\rho_{\binom{8,3,1}{6,5}}),\qquad
\rho^{(1)}_{\binom{8,6,5}{3,1}}=\tfrac{1}{\sqrt 2}(\rho_{\binom{8,6,5}{3,1}}+\rho_{\binom{8,6,3}{5,1}}),
\]
and the linear transformation $\tilde f\colon\calv^{(1)}_Z\rightarrow\calv_{Z^{(1)}}$ given by
$\rho^{(1)}_{\binom{8,5,1}{6,3}}\mapsto\rho_{\binom{7,1}{5}}$ and 
$\rho^{(1)}_{\binom{8,6,5}{3,1}}\mapsto\rho_{\binom{7,5}{1}}$ is an isometry.
\end{example}

\begin{lemma}\label{0904}
Suppose that we are in Case (I) of Subsection~\ref{0903}.
Then
\[
\langle R^{(1)}_\Sigma,\rho^{(1)}_\Lambda\rangle
=\langle R_{\bar f(\Sigma)},\rho_{\bar f(\Lambda)}\rangle
\]
for any $\Sigma=\Lambda_\Psi$ and $\Lambda=\Lambda_M$ where $\Psi,M\subset Z_\rmI\smallsetminus\{a_k,b_l\}$
such that $|\Psi^*|=|\Psi_*|$ and $|M^*|\equiv|M_*|\pmod 2$.
In particular, we have $\tilde f(R^{(1)}_\Sigma)=R_{\bar f(\Sigma)}$.
\end{lemma}
\begin{proof}
Let $\Sigma=\Lambda_\Psi$ and $\Lambda=\Lambda_M$ for $\Psi,M\subset Z_\rmI\smallsetminus\{a_k,b_l\}$
such that $|\Psi^*|=|\Psi_*|$ and $|M^*|\equiv|M_*|\pmod 2$.

First suppose that $\binom{a_k}{b_l}$ is a pair of doubles.
Then both entries $a_k,b_l$ are not in $Z_\rmI$.
It is clear that $x\in\Psi\cap M$ if and only if $f(x)\in f(\Psi)\cap f(M)$.
Therefore $|\Psi\cap M|=|f(\Psi)\cap f(M)|$.
Moreover, $Z$ and $Z^{(1)}$ have the same degree from Remark~\ref{0822},
so the lemma follows from Proposition~\ref{0306} and Proposition~\ref{0319}
immediately.

Next suppose that $\binom{a_k}{b_l}$ is in the core of $\cald_{Z,Z'}$.
Now
\begin{align*}
 \langle R^{(1)}_{\Sigma},\rho^{(1)}_{\Lambda}\rangle
&=\frac{1}{2}(\langle R_{\Sigma},\rho_{\Lambda}\rangle+\langle R_{\Sigma},\rho_{\Lambda^\natural}\rangle
+\langle R_{\Sigma^\natural},\rho_{\Lambda}\rangle+\langle R_{\Sigma^\natural},\rho_{\Lambda^\natural}\rangle) \\
&= \frac{1}{2^{\delta+1}}
((-1)^{\langle \Sigma,\Lambda\rangle}+(-1)^{\langle \Sigma,\Lambda^\natural\rangle}
+(-1)^{\langle \Sigma^\natural,\Lambda\rangle}+(-1)^{\langle \Sigma^\natural,\Lambda^\natural\rangle})
\end{align*}
by Proposition~\ref{0306} where $\delta$ is the degree of $Z$.
From the definition,
it is easy to see that
\begin{align*}
\textstyle|\Psi\cap M|=|\Psi\cap (M\cup\binom{a_k}{b_l})|=|(\Psi\cup\binom{a_k}{b_l})\cap M| &= |f(\Psi)\cap f(M)|,\\
\textstyle|(\Psi\cup\binom{a_k}{b_l})\cap(M\cup\binom{a_k}{b_l})| &= |f(\Psi)\cap f(M)|+2.
\end{align*}
Note that now $Z^{(1)}$ is of degree $\delta-1$ from Remark~\ref{0822}.
Therefore
\[
\langle R^{(1)}_{\Sigma},\rho^{(1)}_{\Lambda}\rangle
= \frac{1}{2^{\delta-1}} (-1)^{\langle\bar f(\Sigma),\bar f(\Lambda)\rangle} =\langle R_{\bar f(\Sigma)},\rho_{\bar f(\Lambda)}\rangle.
\]
\end{proof}

Similarly, we can define a bijection $f'$ from the set of subsets of $Z'_\rmI\smallsetminus\{c_{l'},d_{k'}\}$
onto the set of subsets of $Z'^{(1)}_\rmI$ given by
\begin{align}\label{0815}
\begin{split}
c_i \mapsto\begin{cases}
c_i-1, & \text{if $1\leq i<l'$;}\\
c_i, & \text{if $l'<i\leq m'$},
\end{cases}\qquad
d_i \mapsto\begin{cases}
d_i-1, & \text{if $1\leq i<k'$;}\\
d_i, & \text{if $k'<i\leq m'$}.
\end{cases}
\end{split}
\end{align}
We can define a bijection $\bar f\colon\overline\cals_{Z'}\rightarrow\overline\cals_{Z'^{(1)}}$
if $\binom{c_{l'}}{d_{k'}}$ is a pair of doubles;
a bijection $\bar f\colon\overline\cals_{Z'}^{\Psi'}\rightarrow\overline\cals_{Z'^{(1)}}$
for $\Psi'=\binom{c_{l'}}{d_{k'}}$ if $\binom{c_{l'}}{d_{k'}}$ is in the core of $\cald_{Z,Z'}$.

Let $\rho_{\Lambda'}^{(1)}=\rho_{\Lambda'}$ if $\binom{c_{l'}}{d_{k'}}$ is a pair of doubles,
$\rho_{\Lambda'}^{(1)}=\frac{1}{\sqrt 2}(\rho_{\Lambda'}+\rho_{\Lambda'^\natural})$ where
$\Lambda'^\natural=\Lambda_{N\cup\binom{c_{l'}}{d_{k'}}}$, $\Lambda'=\Lambda_N$ for
$N\subset Z'_\rmI\smallsetminus\{c_{l'},d_{k'}\}$ if $\binom{c_{l'}}{d_{k'}}$ is in the core.
Let $\calv^{+,(1)}_{Z'}$ be the subspace of $\calv^+_{Z'}$ spanned by these $\rho_{\Lambda'}^{(1)}$.
Then we have an isometry $\tilde f'\colon\calv^{+,(1)}_{Z'}\rightarrow\calv^+_{Z'^{(1)}}$ analogously.
And the analogue of Lemma~\ref{0904} is also true.

Define $\overline\calb^{+,(1)}_{Z,Z'}$ to be the subset consisting of
$(\Lambda,\Lambda')\in\overline\calb^+_{Z,Z'}$ where $\Lambda=\Lambda_M$
and $\Lambda'=\Lambda_N$ such that $M\subset Z_\rmI\smallsetminus\{a_k,b_l\}$
and $N\subset Z'_\rmI\smallsetminus\{c_{l'},d_{k'}\}$,
and define
\[\calb_{Z,Z'}^{+,(1)}=\overline\calb_{Z,Z'}^{+,(1)}\cap(\cals_Z\times\cals^+_{Z'})\quad\text{and}\quad
\cald^{(1)}_{Z,Z'}=\cald_{Z,Z'}\cap\overline\calb^{+,(1)}_{Z,Z'}.
\]
If both $\binom{a_k}{b_l}$ and $\binom{c_{l'}}{d_{k'}}$ are pairs of doubles,
then it is clear that $\calb^{+,(1)}_{Z,Z'}=\calb^+_{Z,Z'}$ and
$\cald^{(1)}_{Z,Z'}=\cald_{Z,Z'}$.

\begin{lemma}\label{0912}
For Case (I) of Subsection~\ref{0903},
\[
(\Lambda,\Lambda')\in\overline\calb^{+,(1)}_{Z,Z'}\quad\text{if and only if}\quad
(\bar f(\Lambda),\bar f'(\Lambda'))\in\overline\calb^+_{Z^{(1)},Z'^{(1)}}.
\]
\end{lemma}
\begin{proof}
Write $Z,Z'$ as in (\ref{0914}).
Suppose that $\Lambda=\Lambda_M$ and $\Lambda'=\Lambda_N$ for some $M\subset Z_\rmI\smallsetminus\{a_k,b_l\}$
and $N\subset Z'_\rmI\smallsetminus\{c_{l'},d_{k'}\}$.
Write
\begin{align*}
\Lambda &=\binom{a'_1,a'_2,\ldots,a'_{m_1}}{b'_1,b'_2,\ldots,b'_{m_2}}, &
\Lambda' &=\binom{c'_1,c'_2,\ldots,c'_{m'_1}}{d'_1,d'_2,\ldots,d'_{m'_2}}; \\
\bar f(\Lambda) &=\binom{a''_1,a''_2,\ldots,a''_{m_1-1}}{b''_1,b''_2,\ldots,b''_{m_2-1}}, &
\bar f'(\Lambda') &=\binom{c''_1,c''_2,\ldots,c''_{m_1'-1}}{d''_1,d''_2,\ldots,d''_{m'_2-1}}.
\end{align*}
Because $a_k,b_l\not\in M$ and $c_{l'},d_{k'}\not\in N$,
we see that $a_k$ (resp.~$b_l$) will be in the first (resp.~second) row of $\Lambda$,
similarly $c_{l'}$ (resp.~$d_{k'}$) will be in the first (resp.~second) row of $\Lambda'$.
Then there exit indices $i_0,j_0,j'_0,i'_0$ such that $a'_{i_0}=a_k$,
$b'_{j_0}=b_l$, $c'_{j'_0}=c_{l'}$ and $d'_{i'_0}=d_{k'}$.
Then from the definitions in (\ref{1001}) and (\ref{0815}), we see that
\begin{align*}
a''_i &=\begin{cases}
a'_i-1, & \text{if $i<i_0$};\\
a'_{i+1}, & \text{if $i\geq i_0$},
\end{cases} &\qquad
b''_i &=\begin{cases}
b'_i-1, & \text{if $i<j_0$};\\
b'_{i+1}, & \text{if $i\geq j_0$},
\end{cases} \\
c''_i &=\begin{cases}
c'_i-1, & \text{if $i<j'_0$};\\
c'_{i+1}, & \text{if $i\geq j'_0$},
\end{cases} &\qquad
d''_i &=\begin{cases}
d'_i-1, & \text{if $i<i'_0$};\\
d'_{i+1}, & \text{if $i\geq i'_0$}
\end{cases}
\end{align*}
and
\begin{align*}
a'_i &=\begin{cases}
a''_i+1, & \text{if $i<i_0$};\\
a_k, & \text{if $i=i_0$};\\
a''_{i-1}, & \text{if $i>i_0$},
\end{cases} &\qquad
b'_i &=\begin{cases}
b''_i+1, & \text{if $i<j_0$};\\
b_l, & \text{if $i=j_0$};\\
b''_{i-1}, & \text{if $i>j_0$},
\end{cases} \\
c'_i &=\begin{cases}
c''_i+1, & \text{if $i<j'_0$};\\
c_{l'}, & \text{if $i=j'_0$};\\
c''_{i-1}, & \text{if $i>j'_0$},
\end{cases} &\qquad
d'_i &=\begin{cases}
d''_i+1, & \text{if $i<i'_0$};\\
d_{k'}, & \text{if $i=i'_0$};\\
d''_{i-1}, & \text{if $i>i'_0$}.
\end{cases}
\end{align*}

Now we consider the following cases:
\begin{enumerate}
\item[(1)] Suppose that $m'=m$ and $l=k$.
Then $k'=k-1$, $l'=k$ by (\ref{0828}) and hence
$\max\{a_k,b_l\}=a_k$, $\min\{a_k,b_l\}=b_k$,
$\max\{c_{l'},d_{k'}\}=d_{k-1}$, $\min\{c_{l'},d_{k'}\}=c_k$.
Then $i_0+j_0=k+l=2k$ and
$i'_0+j'_0=l'+k'=2k-1$.
Note that $b_{k-1}$ is the smallest entry in $Z$ which is greater than both $a_k$ and $b_l$,
so we have $a'_{i_0-1}\geq b_{k-1}$ and $b'_{j_0-1}\geq b_{k-1}$.
Similarly, $d_k$ is the largest entry in $Z'$ which is less than both $c_{l'}$ and $d_{k'}$,
so we have $d_k\geq c'_{j'_0+1}$ and $d_k\geq d'_{i'_0+1}$.
\begin{enumerate}
\item If $\binom{c_{l'}}{d_{k'}}$ is a pair of doubles, i.e., $c_k=d_{k-1}$,
then $(Z,Z')\in\cald_{Z,Z'}$ implies that $b_{k-1}>c_k=d_{k-1}$ by Lemma~\ref{0210};
if $\binom{c_{l'}}{d_{k'}}$ is in the core,
then by Lemma~\ref{0301}, we have $b_{k-1}>d_{k-1}$, again.

\item If $\binom{a_k}{b_l}$ is a pair of doubles, i.e., $a_k=b_k$,
then $(Z,Z')\in\cald_{Z,Z'}$ implies $b_k=a_k>d_k$ by Lemma~\ref{0210};
if $\binom{a_k}{b_l}$ is in the core,
then by Lemma~\ref{0701}, we have $b_k>d_k$, again.
\end{enumerate}
Because $\cald_{Z,Z'}\neq\emptyset$, we know that $(Z,Z')\in\cald_{Z,Z'}$ by Proposition~\ref{0216} and we have
\[
a_i>d_i,\quad d_i\geq a_{i+1},\quad c_i\geq b_i,\quad b_i>c_{i+1}
\]
for each $i$ by Lemma~\ref{0210}.

Suppose that $(\Lambda,\Lambda')\in\overline\calb^{+,(1)}_{Z,Z'}$.
Then, by Lemma~\ref{0210}, we have
\begin{equation}\label{0915}
a'_i>d'_i,\quad d'_i\geq a'_{i+1},\quad c'_i\geq b'_i,\quad b'_i>c'_{i+1}
\end{equation}
for each $i$.
Then $a'_{i_0-1}\geq b_{k-1}>d_{k-1}=d'_{i'_0}$, so $i_0-1\leq i'_0$ by (\ref{0915}).
And $b'_{j_0-1}\geq b_{k-1}>d_{k-1}\geq c_k=c'_{j'_0}$, so $j_0-1<j'_0$, i.e., $j_0\leq j'_0$.
But we know that $i_0+j_0=i'_0+j'_0+1$,
so we conclude that $i_0=i'_0+1$ and $j_0=j'_0$.
Then
\begin{itemize}
\item $a''_i=a'_i-1>d'_i-1=d''_i$ if $i<i_0-1$;
$a''_{i_0-1}=a'_{i_0-1}-1>d'_{i_0-1}-1\geq d'_{i_0}=d''_{i_0-1}$;
$a''_i=a'_{i+1}>d'_{i+1}=d''_i$ if $i\geq i_0$,

\item $d''_i=d'_i-1\geq a'_{i+1}-1=a''_{i+1}$ if $i<i_0-1$;
$d''_i=d'_{i+1}\geq a'_{i+2}=a''_{i+1}$ if $i\geq i_0-1$,

\item $c''_i=c'_i-1\geq b'_i-1=b''_i$ if $i<j_0$;
$c''_i=c'_{i+1}\geq b'_{i+1}=b''_i$ if $i\geq j_0$,

\item $b''_i>b'_i-1>c'_{i+1}-1=c''_{i+1}$ if $i<j_0-1$;
$b''_{j_0-1}=b'_{j_0-1}-1\geq b'_{j_0}>c'_{j_0+1}=c''_{j_0}$;
$b''_i=b'_{i+1}>c'_{i+2}=c''_{i+1}$ if $i\geq j_0$.
\end{itemize}
Therefore, we conclude that $(\bar f(\Lambda),\bar f'(\Lambda'))\in\overline\calb^+_{Z^{(1)},Z'^{(1)}}$ by Lemma~\ref{0210}.

Conversely, suppose that $(\bar f(\Lambda),\bar f'(\Lambda'))\in\overline\calb^+_{Z^{(1)},Z'^{(1)}}$.
Then by Lemma~\ref{0210}, we have
\begin{equation}\label{0916}
a''_i>d''_i,\quad d''_i\geq a''_{i+1},\quad c''_i\geq b''_i,\quad b''_i>c''_{i+1}
\end{equation}
for each $i$.
Then $a''_{i_0-1}=a'_{i_0-1}-1\geq b_{k-1}-1>d_{k-1}-1=d'_{i'_0}-1\geq d'_{i'_0+1}=d''_{i'_0}$,
so $i_0-1\leq i'_0$.
And $b''_{j_0-1}=b'_{j_0-1}-1\geq b_{k-1}-1>d_{k-1}-1\geq c_k-1=c'_{j'_0}-1\geq c'_{j'_0+1}=c''_{j'_0}$,
so $j_0-1<j'_0$, i.e., $j_0\leq j'_0$.
Again, we conclude that $i_0=i'_0+1$ and $j_0=j'_0$.
\begin{itemize}
\item $a'_i=a''_i+1>d''_i+1=d'_i$ if $i<i_0-1$;
$a'_{i_0-1}\geq b_{k-1}>d_{k-1}=d_{k'}=d'_{i'_0}=d'_{i_0-1}$;
$a'_{i_0}=a_k>d_k\geq d'_{i'_0+1}=d'_{i_0}$;
$a'_i=a''_{i-1}>d''_{i-1}=d'_i$ if $i>i_0$,

\item $d'_i=d''_i+1\geq a''_{i+1}+1=a'_{i+1}$ if $i<i_0-1$;
$d'_{i_0-1}=d'_{i_0'}=d_{k'}=d_{k-1}\geq a_k=a'_{i_0}$;
$d'_i=d''_{i-1}\geq a''_i=a'_{i+1}$ if $i>i_0-1$,

\item $c'_i=c''_i+1\geq b''_i+1=b'_i$ if $i<j_0$;
$c'_{j_0}=c_{l'}=c_l\geq b_l=b'_{j_0}$;
$c'_i=c''_{i-1}\geq b''_{i-1}=b'_i$ if $i>j_0$,

\item $b'_i=b''_i+1>c''_{i+1}+1=c'_{i+1}$ if $i<j_0-1$;
$b'_{j_0-1}\geq b_{k-1}>c_k=c_{l'}=c'_{j'_0}=c'_{j_0}$;
$b'_{j_0}=b_l=b_k>d_k\geq c'_{j'_0+1}=c'_{j_0+1}$;
$b'_i=b''_{i-1}>c''_i=c'_{i+1}$ if $i>j_0$.
\end{itemize}
Therefore, we conclude that $(\Lambda,\Lambda')\in\overline\calb^{+,(1)}_{Z,Z'}$ by Lemma~\ref{0210}.

\item[(2)] Suppose that $m'=m$ and $l=k-1$.
Then $k'=k-1$, $l'=k-1$ by (\ref{0828}) and hence
$\max\{a_k,b_l\}=b_{k-1}$, $\min\{a_k,b_l\}=a_k$,
$\max\{c_{l'},d_{k'}\}=c_{k-1}$, $\min\{c_{l'},d_{k'}\}=d_{k-1}$.
Then $i_0+j_0=2k-1$ and
$i'_0+j'_0=2k-2$.
The remaining proof is similar to Case (1).

\item[(3)] Suppose that $m'=m+1$ and $l=k$.
Then $k'=k$, $l'=k+1$ by (\ref{0828}) and hence
$\max\{a_k,b_l\}=a_k$, $\min\{a_k,b_l\}=b_k$,
$\max\{c_{l'},d_{k'}\}=d_k$, $\min\{c_{l'},d_{k'}\}=c_{k+1}$.
Then $i_0+j_0=2k$ and
$i'_0+j'_0=2k+1$.
The remaining proof is similar to Case (1).

\item[(4)] Suppose that $m'=m+1$ and $l=k-1$.
Then $k'=k$, $l'=k$ by (\ref{0828}) and hence
$\max\{a_k,b_l\}=b_{k-1}$, $\min\{a_k,b_l\}=a_k$,
$\max\{c_{l'},d_{k'}\}=c_k$, $\min\{c_{l'},d_{k'}\}=d_k$.
Then $i_0+j_0=2k-1$ and
$i'_0+j'_0=2k$.
The remaining proof is similar to Case (1).
\end{enumerate}
\end{proof}

\begin{lemma}\label{0823}
Suppose we are in Case (I) of Subsection~\ref{0903}.
Then
\begin{align*}
\sum_{(\Lambda,\Lambda')\in\calb^+_{Z,Z'}}\rho_\Lambda\otimes\rho_{\Lambda'}
&= C\sum_{(\Lambda,\Lambda')\in\calb^{+,(1)}_{Z,Z'}}\rho^{(1)}_\Lambda\otimes\rho^{(1)}_{\Lambda'}\\
\sum_{(\Sigma,\Sigma')\in\cald_{Z,Z'}}R_\Sigma\otimes R_{\Sigma'}
&= C\sum_{(\Sigma,\Sigma')\in\cald^{(1)}_{Z,Z'}}R^{(1)}_\Sigma\otimes R^{(1)}_{\Sigma'}
\end{align*}
where
\[
C=\begin{cases}
1, & \text{if both $\binom{a_k}{b_l}$ and $\binom{c_{l'}}{d_{k'}}$ are pairs of doubles};\\
\sqrt 2, & \text{if exactly one of $\binom{a_k}{b_l}$ and $\binom{c_{l'}}{d_{k'}}$ is a pair of doubles};\\
2, & \text{if none of $\binom{a_k}{b_l}$ and $\binom{c_{l'}}{d_{k'}}$ is a pair of doubles}.
\end{cases}
\]
\end{lemma}
\begin{proof}
We will only prove the first identity.
The proof of the second identity is similar.
Now we consider the following four cases:
\begin{enumerate}
\item[(1)] Suppose that both $\binom{a_k}{b_l}$ and $\binom{c_{l'}}{d_{k'}}$ are pairs of doubles.
Now $\rho_\Lambda=\rho^{(1)}_\Lambda$, $\rho_{\Lambda'}=\rho^{(1)}_{\Lambda'}$ and
$\calb^+_{Z,Z'}=\calb^{+,(1)}_{Z,Z'}$, so the lemma is clearly true for this case.

\item[(2)] Suppose that $\binom{a_k}{b_l}$ is in the core and $\binom{c_{l'}}{d_{k'}}$ is a pair of doubles.
Now we have $\rho_\Lambda^{(1)}=\frac{1}{\sqrt 2}(\rho_\Lambda+\rho_{\Lambda^\natural})$,
$\rho^{(1)}_{\Lambda'}=\rho_{\Lambda'}$,
and by (\ref{0811}),
\begin{align*}
\calb^+_{Z,Z'}
&=\{\,(\Lambda+\Lambda_M,\Lambda')\mid(\Lambda,\Lambda')\in\calb^{+,(1)}_{Z,Z'},\ M\leq\{\textstyle\binom{a_k}{b_l}\}\,\} \\
&=\{\,(\Lambda,\Lambda'),(\Lambda+\Lambda_{\binom{a_k}{b_l}},\Lambda')\mid(\Lambda,\Lambda')\in\calb^{+,(1)}_{Z,Z'}\,\}.
\end{align*}
Therefore,
\begin{align*}
\sum_{(\Lambda,\Lambda')\in\calb^+_{Z,Z'}}\rho_\Lambda\otimes\rho_{\Lambda'}
&= \sum_{(\Lambda,\Lambda')\in\calb^{+,(1)}_{Z,Z'}}(\rho_\Lambda+\rho_{\Lambda+\Lambda_{\binom{a_k}{b_l}}})\otimes\rho_{\Lambda'} \\
&= \sqrt{2}\sum_{(\Lambda,\Lambda')\in\calb^{+,(1)}_{Z,Z'}}\rho^{(1)}_\Lambda\otimes\rho^{(1)}_{\Lambda'}.
\end{align*}

\item[(3)] Suppose that $\binom{a_k}{b_l}$ is a pair of double and $\binom{c_{l'}}{d_{k'}}$ is in the core.
The proof for this case is similar to that of (2).

\item[(4)] Suppose that both $\binom{a_k}{b_l}$ and $\binom{c_{l'}}{d_{k'}}$ are in the cores.
Now we have $\rho_\Lambda^{(1)}=\frac{1}{\sqrt 2}(\rho_\Lambda+\rho_{\Lambda^\natural})$,
$\rho^{(1)}_{\Lambda'}=\frac{1}{\sqrt 2}(\rho_{\Lambda'}+\rho_{\Lambda'^\natural})$,
and by (\ref{0811}),
\begin{align*}
\calb^+_{Z,Z'}
&=\{\,(\Lambda+\Lambda_M,\Lambda'+\Lambda_N)\mid(\Lambda,\Lambda')\in\calb^{+,(1)}_{Z,Z'},\ M\leq\{\textstyle\binom{a_k}{b_l}\},\ N\leq\{\textstyle\binom{c_{l'}}{d_{k'}}\}\,\} \\
&=\Bigl\{\,(\Lambda,\Lambda'),(\Lambda+\Lambda_{\binom{a_k}{b_l}},\Lambda'),(\Lambda,\Lambda'+\Lambda_{\binom{c_{l'}}{d_{k'}}}),
(\Lambda+\Lambda_{\binom{a_k}{b_l}},\Lambda'+\Lambda_{\binom{c_{l'}}{d_{k'}}})\,\\
&\qquad\qquad\qquad\qquad\qquad\qquad\qquad\qquad\qquad\qquad\quad\Bigl|\,(\Lambda,\Lambda')\in\calb^{+,(1)}_{Z,Z'}\,\Bigr\}.
\end{align*}
Therefore,
\begin{align*}
\sum_{(\Lambda,\Lambda')\in\calb^+_{Z,Z'}}\rho_\Lambda\otimes\rho_{\Lambda'}
&= \sum_{(\Lambda,\Lambda')\in\calb^{+,(1)}_{Z,Z'}}(\rho_\Lambda+\rho_{\Lambda+\Lambda_{\binom{a_k}{b_l}}})
\otimes(\rho_{\Lambda'}+\rho_{\Lambda'+\Lambda_{\binom{c_{l'}}{d_{k'}}}}) \\
&= 2\sum_{(\Lambda,\Lambda')\in\calb^{+,(1)}_{Z,Z'}}\rho^{(1)}_\Lambda\otimes\rho^{(1)}_{\Lambda'}.
\end{align*}
\end{enumerate}
\end{proof}

\begin{lemma}\label{0819}
Suppose we are in Case (I) of Subsection~\ref{0903}.
Then
\begin{align*}
\tilde f\otimes\tilde f'\Biggl(\sum_{(\Lambda,\Lambda')\in\calb^+_{Z,Z'}}\rho_\Lambda\otimes\rho_{\Lambda'}\Biggr)
&=C\sum_{(\Lambda^{(1)},\Lambda'^{(1)})\in\calb^+_{Z^{(1)},Z'^{(1)}}}\rho_{\Lambda^{(1)}}\otimes\rho_{\Lambda'^{(1)}} \\
\tilde f\otimes\tilde f'\Biggl(\sum_{(\Sigma,\Sigma')\in\cald_{Z,Z'}}R_\Sigma\otimes R_{\Sigma'}\Biggr)
&=C\sum_{(\Sigma^{(1)},\Sigma'^{(1)})\in\cald_{Z^{(1)},Z'^{(1)}}} R_{\Sigma^{(1)}}\otimes R_{\Sigma'^{(1)}}
\end{align*}
where $C$ is the constant given in Lemma~\ref{0823}.
\end{lemma}
\begin{proof}
The lemma follows from Lemma~\ref{0912} and Lemma~\ref{0823} directly.
\end{proof}

\subsection{Cases (II) and (III)}\label{0816}
Let $Z,Z'$ be given as in (\ref{0914}).

\begin{enumerate}
\item[(1)] Suppose now we are in Case (II) of Subsection~\ref{0903}.
Then $m'=m$ by Lemma~\ref{0905}.
Define a bijection $f\colon Z\smallsetminus\{a_k,b_l\}\rightarrow Z^{(1)}$ given by
\begin{align}\label{0906}
\begin{split}
a_i \mapsto\begin{cases}
a_i-1, & \text{if $1\leq i<k$;}\\
d_{i-1}, & \text{if $k<i\leq m+1$},
\end{cases}\qquad
b_i \mapsto\begin{cases}
b_i-1, & \text{if $1\leq i<l$;}\\
c_{i}, & \text{if $l<i\leq m$}.
\end{cases}
\end{split}
\end{align}
Similarly, define a bijection $f'\colon Z'\rightarrow Z'^{(1)}$ given by
\begin{align}\label{0907}
\begin{split}
c_i \mapsto\begin{cases}
c_i, & \text{if $1\leq i\leq l'$;}\\
b_{i}, & \text{if $l'<i\leq m'$},
\end{cases}\qquad
d_i \mapsto\begin{cases}
d_i, & \text{if $1\leq i\leq k'$;}\\
a_{i+1}, & \text{if $k'<i\leq m'$}.
\end{cases}
\end{split}
\end{align}

\item[(2)] Suppose now we are in Case (III) of Subsection~\ref{0903}.
Then we have $m'=m+1$ by Lemma~\ref{0905}.
Define a bijection $f\colon Z\rightarrow Z^{(1)}$ given by
\begin{align}\label{0918}
\begin{split}
a_i \mapsto\begin{cases}
a_i, & \text{if $1\leq i\leq k$;}\\
d_i, & \text{if $k<i\leq m+1$},
\end{cases}\qquad
b_i \mapsto\begin{cases}
b_i, & \text{if $1\leq i\leq l$;}\\
c_{i+1}, & \text{if $l<i\leq m$}.
\end{cases}
\end{split}
\end{align}
Similarly, define a bijection $f'\colon Z'\smallsetminus\{c_{l'},d_{k'}\}\rightarrow Z'^{(1)}$ given by
\begin{align}\label{0919}
\begin{split}
c_i \mapsto\begin{cases}
c_i-1, & \text{if $1\leq i<l'$;}\\
b_{i-1}, & \text{if $l'<i\leq m'$},
\end{cases}\qquad
d_i \mapsto\begin{cases}
d_i-1, & \text{if $1\leq i<k'$;}\\
a_i, & \text{if $k'<i\leq m'$}.
\end{cases}
\end{split}
\end{align}
\end{enumerate}

\begin{lemma}
\begin{enumerate}
\item[(i)] Suppose that we are in Case (II).
The map in (\ref{0906}) gives a bijection from $Z_\rmI\smallsetminus\{a_k,b_l\}$ onto $Z^{(1)}_\rmI$;
and the map in (\ref{0907}) gives a bijection from $Z'_\rmI$ onto $Z'^{(1)}_\rmI$.

\item[(ii)] Suppose that we are in Case (III).
The map in (\ref{0918}) gives a bijection from $Z_\rmI$ onto $Z^{(1)}_\rmI$;
and the map in (\ref{0919}) gives a bijection from $Z'_\rmI\smallsetminus\{c_{l'},d_{k'}\}$ onto $Z'^{(1)}_\rmI$.
\end{enumerate}
\end{lemma}
\begin{proof}
Suppose we are in Case (II).
Clearly $f$ gives a bijection between entries in $Z\smallsetminus\{a_k,b_l\}$
and entries in $Z^{(1)}$.
So we only need to check that $f$ maps singles to singles and non-singles to non-singles.

Write $Z^{(1)}$ as in (\ref{0913}).
Suppose that $\binom{a_i}{b_j}$ is a pair of doubles in $Z\smallsetminus\binom{a_k}{b_l}$.
From the assumption of $\binom{a_k}{b_l}$,
we see that $i<k$ and $j<l$ and thus $f(a_i)=f(b_j)$.
Conversely, if $f(a_i)=f(b_j)$, then we must have $i<k$ and $j<l$
and thus $a_i=b_j$.

Next suppose that $\binom{c_j}{d_i}$ is a pair of doubles in $Z'$.
From the assumption on $\binom{c_{l'}}{d_{k'}}$,
wee see that $i<k'$ and $j<l'$ and thus $f'(c_j)=f'(d_i)$.
Conversely, if $f'(c_j)=f'(d_i)$,
then we must have $i<k'$ and $j<l'$ and thus $c_j=d_i$.

The proof for Case (III) is similar.
\end{proof}

\begin{enumerate}
\item[(1)] Suppose now we are in Case (II) of Subsection~\ref{0903}.
Let $f\colon Z\smallsetminus\{a_k,b_l\}\rightarrow Z^{(1)}$ be the bijection given as in (\ref{0906}),
and let $f'\colon Z'\rightarrow Z'^{(1)}$ be given as in (\ref{0907}).
As in Case (I), $f$ induces a bijection $\bar f\colon\cals_Z\rightarrow\cals_{Z^{(1)}}$ if
$\binom{a_k}{b_l}$ is a pair of doubles;
a bijection $\bar f\colon\cals_Z^\Psi\rightarrow\cals_{Z^{(1)}}$ where $\Psi=\binom{a_k}{b_l}$ if $\binom{a_k}{b_l}$
is in the core.
Then we have a vector space isomorphism $\tilde f\colon\calv_Z^{(1)}\rightarrow\calv_{Z^{(1)}}$
given by $\rho_\Lambda^{(1)}\mapsto\rho_{\bar f(\Lambda)}$
where $\rho_\Lambda^{(1)}$ is given as in (\ref{0930}).
Moreover, we have a bijection $\bar f'\colon\cals_{Z'}^+\rightarrow\cals^+_{Z'^{(1)}}$
and a vector space isomorphism $\tilde f'\colon\calv^{+,(1)}_{Z'}\rightarrow\calv^+_{Z'^{(1)}}$
given by $\rho^{(1)}_{\Lambda'}\mapsto\rho_{\bar f'(\Lambda')}$
where $\rho_{\Lambda'}^{(1)}$ is defined to be $\rho_{\Lambda'}$.
Now $\calb^{+,(1)}_{Z,Z'}$ is defined to be the subset consisting of
$(\Lambda,\Lambda')\in\calb^+_{Z,Z'}$ where $\Lambda=\Lambda_M$
and $\Lambda'=\Lambda_N$ such that $M\subset Z_\rmI\smallsetminus\{a_k,b_l\}$
and $N\subset Z'_\rmI$,
and $\cald^{(1)}_{Z,Z'}=\cald_{Z,Z'}\cap\calb^{+,(1)}_{Z,Z'}$.

\item[(2)] Suppose now we are in Case (III) of Subsection~\ref{0903}.
Let $f\colon Z\rightarrow Z^{(1)}$ be the bijection given as in (\ref{0918}),
and let $f'\colon Z'\smallsetminus\{c_{l'},d_{k'}\}\rightarrow Z'^{(1)}$ be given as in (\ref{0919}).
Now $f$ induces a bijection $\bar f\colon\cals_Z\rightarrow\cals_{Z^{(1)}}$ and then we have a vector space isomorphism
$\tilde f\colon\calv^{(1)}_Z\rightarrow\calv_{Z^{(1)}}$ given by
$\rho^{(1)}_\Lambda\mapsto\rho_{\bar f(\Lambda)}$ where $\rho_{\Lambda}^{(1)}$ is defined to be $\rho_{\Lambda}$.
Moreover, as in Case (I), $f'$ induces a bijection $\bar f'\colon\cals_{Z'}^+\rightarrow\cals^+_{Z'^{(1)}}$ if
$\binom{c_{l'}}{d_{k'}}$ is a pair of doubles; a bijection $\bar f'\colon\cals^{+,\Psi'}_{Z'}\rightarrow\cals^+_{Z'^{(1)}}$
where $\Psi'=\binom{c_{l'}}{d_{k'}}$ if $\binom{c_{l'}}{d_{k'}}$ is in the core.
Then we have a vector space isomorphism $\tilde f'\colon\calv_{Z'}^{+,(1)}\rightarrow\calv^+_{Z'^{(1)}}$ given by
$\rho_{\Lambda'}^{(1)}\mapsto\rho_{\bar f'(\Lambda')}$ where $\rho_{\Lambda'}^{(1)}$ is defined as in (\ref{0930}).
Now $\calb^{+,(1)}_{Z,Z'}$ is the subset consisting of
$(\Lambda,\Lambda')\in\calb^+_{Z,Z'}$ where $\Lambda=\Lambda_M$
and $\Lambda'=\Lambda_N$ such that $M\subset Z_\rmI$ and $N\subset Z'_\rmI\smallsetminus\{c_{l'},d_{k'}\}$,
and $\cald^{(1)}_{Z,Z'}=\cald_{Z,Z'}\cap\calb^{+,(1)}_{Z,Z'}$.
\end{enumerate}

\begin{example}
Keep the notation in Example~\ref{0829}, i.e.,
$Z^{(1)}=\binom{7,1}{5}$ and $Z'^{(1)}=\binom{7,5}{5,0}$.
Next we see that $\binom{c'_2}{d'_1}=\binom{5}{5}$ is a pair of doubles,
and $\binom{a'_1}{b'_1}=\binom{7}{5}$ is neither a pair of doubles nor in the core.
So now we are in Case (III) and then
$Z^{(2)}=\binom{7,0}{5}$ and $Z'^{(2)}=\binom{6}{1}$.
Then finally $\calb^+_{Z^{(2)},Z'^{(2)}}$ is one-to-one:
\[
\begin{tabular}{c|cc}
\toprule
& $\binom{6}{1}$ & $\binom{1}{6}$ \\
\midrule
$\binom{7,0}{5}$ & \checkmark & \\
$\binom{7,5}{0}$ & & \checkmark \\
\bottomrule
\end{tabular}
\]
and both $Z^{(2)},Z'^{(2)}$ are regular.
\end{example}

\begin{lemma}\label{0909}
For Case (II) and Case (III),
\[
(\Lambda,\Lambda')\in\overline\calb^{+,(1)}_{Z,Z'}\quad\text{if and only if}\quad
(\bar f(\Lambda),\bar f'(\Lambda'))\in\overline\calb^+_{Z^{(1)},Z'^{(1)}}.
\]
\end{lemma}
\begin{proof}
Write $Z,Z'$ as in (\ref{0914}).
First, suppose that we are in Case (II).
Then $m'=m$ by Lemma~\ref{0905}.
Because $\cald_{Z,Z'}\neq\emptyset$, we know that $(Z,Z')\in\cald_{Z,Z'}$ by Proposition~\ref{0216} and we have
\[
a_i>d_i,\quad d_i\geq a_{i+1},\quad c_i\geq b_i,\quad b_i>c_{i+1}
\]
for each $i$ by Lemma~\ref{0210}.

Suppose that $\Lambda=\Lambda_M$, $\Lambda'=\Lambda_N$
for some $M\subset Z_\rmI\smallsetminus\{a_k,b_l\}$ and some $N\subset Z'_\rmI$,
and write
\begin{align*}
\Lambda &=\binom{a'_1,a'_2,\ldots,a'_{m_1}}{b'_1,b'_2,\ldots,b'_{m_2}}, &
\Lambda' &=\binom{c'_1,c'_2,\ldots,c'_{m'_1}}{d'_1,d'_2,\ldots,d'_{m'_2}}; \\
\bar f(\Lambda) &=\binom{a''_1,a''_2,\ldots,a''_{m_1-1}}{b''_1,b''_2,\ldots,b''_{m_2-1}}, &
\bar f'(\Lambda') &=\binom{c''_1,c''_2,\ldots,c''_{m_1'}}{d''_1,d''_2,\ldots,d''_{m'_2}}.
\end{align*}
We know that $m'_1=m_2$ and $m_2'=m_1-1$ by Lemma~\ref{0210}.
Because $a_k,b_l\not\in M$,
there exit indices $i_0,j_0$ such that $a'_{i_0}=a_k$,
$b'_{j_0}=b_l$.
And there are indices $j'_0,i'_0$ such that
\begin{align}\label{0817}
&\begin{cases}
c'_i\geq\min\{c_{l'},d_{k'}\}, & \text{if $i\leq j'_0$};\\
c'_i<\min\{c_{l'},d_{k'}\}, & \text{if $i>j'_0$},
\end{cases}, &\qquad
&\begin{cases}
d'_i\geq\min\{c_{l'},d_{k'}\}, & \text{if $i\leq i'_0$};\\
d'_i<\min\{c_{l'},d_{k'}\}, & \text{if $i>i'_0$}.
\end{cases}
\end{align}
It is clear that $i_0+j_0=k+l$ and $i'_0=j'_0=l'+k'$.
Similarly, there exist indices $i_1,j_1,j'_1,i'_1$ such that
\begin{align*}
&\begin{cases}
a''_i\geq\max\{a_k,b_l\}, & \text{if $i\leq i_1$};\\
a''_i<\min\{a_k,b_l\}, & \text{if $i>i_1$},
\end{cases} &\qquad
&\begin{cases}
b''_i\geq\max\{a_k,b_l\}, & \text{if $i\leq j_1$};\\
b''_i<\min\{a_k,b_l\}, & \text{if $i>j_1$}.
\end{cases} \\
&\begin{cases}
c''_i\geq\min\{c_{l'},d_{k'}\}, & \text{if $i\leq j'_1$};\\
c''_i<\min\{c_{l'},d_{k'}\}, & \text{if $i>j'_1$},
\end{cases} &\qquad
&\begin{cases}
d''_i\geq\min\{c_{l'},d_{k'}\}, & \text{if $i\leq i'_1$};\\
d''_i<\min\{c_{l'},d_{k'}\}, & \text{if $i>i'_1$}.
\end{cases}
\end{align*}
It is clear that $i_1+j_1=(k-1)+(l-1)=k+l-2$ and $i'_1+j'_1=l'+k'$.

Then from (\ref{0906}) and (\ref{0907}),
we have
\begin{align*}
a''_i &=\begin{cases}
a'_i-1, & \text{if $i<i_0$};\\
d'_i, & \text{if $i\geq i_0$},
\end{cases}, &\qquad
b''_i &=\begin{cases}
b'_i-1, & \text{if $i<j_0$};\\
c'_{i+1}, & \text{if $i\geq j_0$},
\end{cases} \\
c''_i &=\begin{cases}
c'_i, & \text{if $i\leq j'_0$};\\
b'_{i}, & \text{if $i>j'_0$},
\end{cases}, &\qquad
d''_i &=\begin{cases}
d'_i, & \text{if $i\leq i'_0$};\\
a'_{i+1}, & \text{if $i>i'_0$}
\end{cases}
\end{align*}
and
\begin{align*}
a'_i &=\begin{cases}
a''_i+1, & \text{if $i\leq i_1$};\\
a_k, & \text{if $i=i_1+1$};\\
d''_{i-1}, & \text{if $i> i_1+1$},
\end{cases} &\qquad
b'_i &=\begin{cases}
b''_i+1, & \text{if $i\leq j_1$};\\
b_l, & \text{if $i=j_1+1$};\\
c''_i, & \text{if $i> j_1+1$},
\end{cases} \\
c'_i &=\begin{cases}
c''_i, & \text{if $i\leq j'_1$};\\
b''_{i-1}, & \text{if $i>j'_1$},
\end{cases} &\qquad
d'_i &=\begin{cases}
d''_i, & \text{if $i\leq i'_1$};\\
a''_i, & \text{if $i>i'_1$}.
\end{cases}
\end{align*}

Now we consider the following cases:
\begin{enumerate}
\item[(1)] Suppose that $l=k$.
Then $k'=k-1$, $l'=l=k$, $\max\{a_k,b_l\}=a_k$, $\min\{a_k,b_l\}=b_k$, $\max\{c_{l'},d_{k'}\}=d_{k-1}$ and $\min\{c_{l'},d_{k'}\}=c_k$.
Then $i_0+j_0=2k$, $i'_0+j'_0=2k-1$, $i_1+j_1=2k-2$ and $i'_1+j'_1=2k-1$.

\begin{enumerate}
\item If $\binom{a_k}{b_l}$ is a pair of doubles,
then $c_k\geq b_k=a_k$;

\item if $\binom{a_k}{b_l}$ is in the core of $\cald_{Z,Z'}$,
then we still have $c_k\geq a_k$ by Lemma~\ref{0701}.
\end{enumerate}

Suppose that $(\Lambda,\Lambda')\in\overline\calb^{+,(1)}_{Z,Z'}$.
Then we have
\[
a'_i>d'_i,\quad d'_i\geq a'_{i+1},\quad c'_i\geq b'_i,\quad b'_i>c'_{i+1}
\]
for each $i$ by Lemma~\ref{0210}.
Then $d'_{i'_0}\geq c_k\geq a_k=a'_{i_0}$,
so $i'_0\leq i_0-1$.
And $c'_{j'_0}\geq c_k\geq b_k=b'_{j_0}$,
so $j'_0\leq j_0$.
But now $i'_0+j'_0=i_0+j_0-1$, so $i'_0=i_0-1$ and $j'_0=j_0$.
Note that now $Z^{(1)}$ is of size $(m,m-1)$ and $Z'^{(1)}$ is of size $(m,m)$.
We have
\begin{itemize}
\item $a''_i=a'_i-1\geq d'_i=d''_i$ if $i\leq i_0$;
$a''_i=d'_i\geq a'_{i+1}=d''_i$ if $i>i_0$,

\item $d''_i=d'_i>a'_{i+1}-1=a''_{i+1}$ if $i<i_0$;
$d''_{i_0}=d'_{i_0}>d'_{i_0+1}=a''_{i_0+1}$;
$d''_i=a'_{i+1}>d'_{i+1}=a''_{i+1}$ if $i>i_0$,

\item $c''_i=c'_i>b'_i-1=b''_i$ if $i\leq j_0$;
$c''_{j_0+1}=c'_{j_0+1}>c'_{j_0+2}=b''_{j_0+1}$;
$c''_i=b'_i>c'_{i+1}=b''_i$ if $i>j_0$,

\item $b''_i=b'_i-1\geq c'_{i+1}=c''_{i+1}$ if $i\leq j_0$;
$b''_i=c'_{i+1}\geq b'_{i+1}=c''_{i+1}$ if $i>j_0$.
\end{itemize}
These conditions imply that $(\bar f(\Lambda),\bar f'(\Lambda'))\in\overline\calb^+_{Z^{(1)},Z'^{(1)}}$.

Conversely, suppose that $(\bar f(\Lambda),\bar f'(\Lambda'))\in\overline\calb^+_{Z^{(1)},Z'^{(1)}}$.
So we have
\[
a''_i\geq d''_i,\quad d''_i>a''_{i+1},\quad c''_i>b''_i,\quad b''_i\geq c''_{i+1}
\]
for each $i$ by Lemma~\ref{0210}.
Now $d''_{i'_1}\geq c_k\geq b_k\geq a''_{i_1+1}$,
so $i'_1\leq i_1$.
And $c''_{j'_1}\geq c_k\geq b_k\geq b''_{j_1+1}$,
so $j'_1\leq j_1+1$.
But now $i'_1+j'_1=i_1+j_1+1=2k-1$, so $i'_1=i_1$ and $j'_1=j_1+1$.
Then we have
\begin{itemize}
\item $a'_i=a''_i+1>d''_i=d'_i$ if $i\leq i_1$;
$a'_{i_1+1}=a_k>d_k\geq d''_{i_1+1}=d'_{i_1+1}$;
$a'_i=d''_{i-1}>a''_i=d'_i$ if $i>i_1+1$,

\item $d'_i=d''_i \geq a''_i+1=a'_{i+1}$ if $i<i_1$;
$d'_{i_1}=d''_{i_1}\geq c_k\geq a_k=a'_{i_1+1}$;
$d'_i=a''_i\geq d''_i=a'_{i+1}$ if $i>i_1$,

\item $c'_i=c''_i\geq b''_i+1=b'_i$ if $i\leq j_1$;
$c'_{j_1+1}=b''_{j_1}\geq b_l=b'_{j_1+1}$;
$c'_i=b''_{i-1}\geq c''_i=b'_i$ if $i>j_1+1$,

\item $b'_i=b''_i+1>c''_{i+1}=c'_{i+1}$ if $i<j_1$;
$b'_{j_1}=b''_{j_1}+1>b''_{j_1}=c'_{j_1+1}$;
$b'_{j_1+1}=b_l>b''_{j_1+1}=c'_{j_1+2}$;
$b'_i=c''_i>b''_i=c'_{i+1}$ if $i>j_1+1$.
\end{itemize}
Note that $Z$ is of size $(m+1,m)$ and $Z'$ is of size $(m,m)$.
Then the above inequalities imply $(\Lambda,\Lambda')\in\overline\calb^{+,(1)}_{Z,Z'}$ by Lemma~\ref{0210}.

\item[(2)] Suppose that $l=k-1$.
Then $k'=k-1$, $l'=l=k-1$, $\max\{a_k,b_l\}=b_{k-1}$, $\min\{a_k,b_l\}=a_k$,
$\max\{c_{l'},d_{k'}\}=c_{k-1}$ and $\min\{c_{l'},d_{k'}\}=d_{k-1}$.
Then $i_0+j_0=2k-1$, $i'_0+j'_0=2k-2$, $i_1+j_1=2k-3$ and $i'_1+j'_1=2k-2$.

\begin{enumerate}
\item If $\binom{a_k}{b_l}$ is a pair of doubles,
then $d_{k-1}\geq a_k=b_{k-1}$;

\item if $\binom{a_k}{b_l}$ is in the core,
then by Lemma~\ref{0701}, we still have $d_{k-1}\geq b_{k-1}$.
\end{enumerate}
Then $d'_{i'_0}\geq d_{k-1}\geq a_k=a'_{i_0}$ and $c'_{j'_0}\geq d_{k-1}\geq b_{k-1}=b'_{j_0}$.
As in (1), the assumption $(\Lambda,\Lambda')\in\overline\calb^{+,(1)}_{Z,Z'}$ will imply
that $i'_0\leq i_0-1$ and $j'_0\leq j_0$.
Hence $i'_0=i_0-1$ and $j'_0=j_0$.
Then we can show that $(\bar f(\Lambda),\bar f'(\Lambda'))\in\overline\calb^+_{Z^{(1)},Z'^{(1)}}$ as in (1).

Conversely, we have
$d''_{i'_1}\geq d_{k-1}\geq a_k\geq a''_{i_1+1}$ and
$c''_{j'_1}\geq d_{k-1}\geq a_k\geq b''_{j_1+1}$.
As in (1), the assumption $(\bar f(\Lambda),\bar f'(\Lambda'))\in\overline\calb^+_{Z^{(1)},Z'^{(1)}}$ will imply
that $i'_1<i_1+1$ and $j'_1\leq j_1+1$.
Hence $i'_1=i_1$ and $j'_1=j_1+1$.
Then we can show that $(\Lambda,\Lambda')\in\overline\calb^{+,(1)}_{Z,Z'}$ as in (1).
\end{enumerate}

The proof for Case (III) is similar.
\end{proof}

\begin{lemma}\label{0824}
Suppose we are in Case (II) or Case (III) of Subsection~\ref{0903}.
Then
\begin{align*}
\sum_{(\Lambda,\Lambda')\in\calb^+_{Z,Z'}}\rho_\Lambda\otimes\rho_{\Lambda'}
&= C\sum_{(\Lambda,\Lambda')\in\calb^{+,(1)}_{Z,Z'}}\rho^{(1)}_\Lambda\otimes\rho^{(1)}_{\Lambda'}\\
\sum_{(\Sigma,\Sigma')\in\cald_{Z,Z'}}R_\Sigma\otimes R_{\Sigma'}
&= C\sum_{(\Sigma,\Sigma')\in\cald^{(1)}_{Z,Z'}}R^{(1)}_\Sigma\otimes R^{(1)}_{\Sigma'}
\end{align*}
where
\[
C=\begin{cases}
1, & \text{if $\binom{a_k}{b_l}$ (resp.~$\binom{c_{l'}}{d_{k'}}$) is a pair of doubles for Case (II) (resp.~Case (III))};\\
\sqrt 2, & \text{if $\binom{a_k}{b_l}$ (resp.~$\binom{c_{l'}}{d_{k'}}$) is in the core for Case (II) (resp.~Case (III))}.
\end{cases}
\]
\end{lemma}
\begin{proof}
The proof is similar to that of Lemma~\ref{0823}.
\end{proof}

\begin{lemma}
Suppose we are in Cases (II) or (III) of Subsection~\ref{0903}.
We have
\begin{align*}
\tilde f\otimes\tilde f'\Biggl(\sum_{(\Lambda,\Lambda')\in\calb^+_{Z,Z'}}\rho_\Lambda\otimes\rho_{\Lambda'}\Biggr)
&=C\sum_{(\Lambda^{(1)},\Lambda'^{(1)})\in\calb^+_{Z^{(1)},Z'^{(1)}}}\rho_{\Lambda^{(1)}}\otimes\rho_{\Lambda'^{(1)}}, \\
\tilde f\otimes\tilde f'\Biggl(\sum_{(\Sigma,\Sigma')\in\cald_{Z,Z'}}R_\Sigma\otimes R_{\Sigma'}\Biggr)
&=C\sum_{(\Sigma^{(1)},\Sigma'^{(1)})\in\cald_{Z^{(1)},Z'^{(1)}}} R_{\Sigma^{(1)}}\otimes R_{\Sigma'^{(1)}}
\end{align*}
where $C$ is the constant given in Lemma~\ref{0824}.
\end{lemma}
\begin{proof}
The lemma follows from Lemma~\ref{0909} and Lemma~\ref{0824} directly.
\end{proof}

\subsection{$\cald_{Z,Z'}$ general}\label{0911}
Now we prove that Theorem~\ref{0310} holds when $\epsilon=+$.

\begin{proposition}\label{0813}
Let $(\bfG,\bfG')=(\Sp_{2n},\rmO^+_{2n'})$,
$Z,Z'$ special symbols of rank $n,n'$ and defect $1,0$ respectively.
Then
\[
\frac{1}{2}\sum_{(\Sigma,\Sigma')\in\cald_{Z,Z'}} R^\bfG_\Sigma\otimes R^{\bfG'}_{\Sigma'}
=\sum_{(\Lambda,\Lambda')\in\calb^+_{Z,Z'}}\rho_\Lambda^\sharp\otimes\rho^\sharp_{\Lambda'}.
\]
\end{proposition}
\begin{proof}
Suppose that $\epsilon=+$,
and let $Z,Z'$ be special symbols of defect $1,0$ respectively
such that $\cald_{Z,Z'}\neq\emptyset$.
If $\cald_{Z,Z'}$ is not one-to-one or one of $Z,Z'$ is not regular,
we are in one of the three cases (I),(II), (III) of Subsection~\ref{0903}.

After Applying the construction in the previous two subsections several times,
we can find special symbols $Z^{(t)},Z'^{(t)}$ of defects $1,0$ respectively and bijections
\[
\bar f_t\colon\cals_Z^{\Psi_0}\longrightarrow\cals_{Z^{(t)}}\qquad\text{and}\qquad
\bar f'_t\colon\cals_{Z'}^{+,\Psi'_0}\longrightarrow\cals^+_{Z'^{(t)}}
\]
such that
\begin{itemize}
\item both $Z^{(t)},Z'^{(t)}$ are regular,

\item $\deg(Z\smallsetminus\Psi_0)=\deg(Z^{(t)})$ and $\deg(Z'\smallsetminus\Psi'_0)=\deg(Z'^{(t)})$,

\item $\cald_{Z^{(t)},Z'^{(t)}}$ is one-to-one,
in particular, $\cald_{Z^{(t)},Z'^{(t)}}\neq\emptyset$,

\item $(\Lambda,\Lambda')\in\calb^{+,\natural}_{Z,Z'}$ if and only if
$(\bar f_t(\Lambda),\bar f'_t(\Lambda'))\in\calb^+_{Z^{(t)},Z'^{(t)}}$.
In particular, $(\Sigma,\Sigma')\in\cald_{Z,Z'}^\natural$ if and only if
$(\bar f_t(\Sigma),\bar f'_t(\Sigma'))\in\cald_{Z^{(t)},Z'^{(t)}}$
\end{itemize}
where $\Psi_0,\Psi'_0$ are the cores of $\cald_{Z,Z'}$ in $Z_\rmI,Z'_\rmI$ respectively,
and $\calb_{Z,Z'}^{+,\natural},\cald_{Z,Z'}^\natural$ are defined in Subsection~\ref{0712}.
The bijection $\bar f_t\times\bar  f'_t$ induces a vector space isomorphism $\tilde f_t\otimes\tilde f'_t$ from some subspace
$\calv_Z^{(t)}\otimes\calv^{+,(t)}_{Z'}$ of $\calv_Z\otimes\calv^+_{Z'}$ onto $\calv_{Z^{(t)}}\otimes\calv^+_{Z'^{(t)}}$ such that
\begin{align*}
\tilde f_t\otimes\tilde f'_t\Biggl(\sum_{(\Lambda,\Lambda')\in\calb^+_{Z,Z'}}\rho_\Lambda\otimes\rho_{\Lambda'}\Biggr)
&=\sum_{(\Lambda^{(t)},\Lambda'^{(t)})\in\calb^+_{Z^{(t)},Z'^{(t)}}}\rho_{\Lambda^{(t)}}\otimes\rho_{\Lambda'^{(t)}} \\
\tilde f_t\otimes\tilde f'_t\Biggl(\sum_{(\Sigma,\Sigma')\in\cald_{Z,Z'}}R_\Sigma\otimes R_{\Sigma'}\Biggr)
&=\sum_{(\Sigma^{(t)},\Sigma'^{(t)})\in\cald_{Z^{(t)},Z'^{(t)}}} R_{\Sigma^{(t)}}\otimes R_{\Sigma'^{(t)}}.
\end{align*}

Now the special symbols $Z^{(t)},Z'^{(t)}$ are regular and the relation $\cald_{Z^{(t)},Z'^{(t)}}$ is one-to-one.
So by Lemma~\ref{0213}, we see that
either $\deg(Z'^{(t)})=\deg(Z^{(t)})$ or $\deg(Z'^{(t)})=\deg(Z^{(t)})+1$.
From the results in Lemma~\ref{0902} and bijections $\bar f_t,\bar f'_t$,
we have the following two situations:
\begin{enumerate}
\item[(1)] Suppose that $\deg(Z'^{(t)})=\deg(Z^{(t)})$.
Then we have a mapping $\theta^+\colon\cals^{+,\Psi'_0}_{Z'}\rightarrow\cals_Z^{\Psi_0}$
such that $(\Lambda,\Lambda')\in\calb^{+,\natural}_{Z,Z'}$ if and only if
$\Lambda=\theta^+(\Lambda')$.

\item[(2)] Suppose that $\deg(Z'^{(t)})=\deg(Z^{(t)})+1$.
Then we have a mapping $\theta^+\colon\cals_Z^{\Psi_0}\rightarrow\cals^{+,\Psi'_0}_{Z'}$
such that $(\Lambda,\Lambda')\in\calb^{+,\natural}_{Z,Z'}$ if and only if
$\Lambda'=\theta^+(\Lambda)$.
\end{enumerate}
Hence by Proposition~\ref{0806}, we have
\[
\Biggl(\sum_{(\Lambda^{(t)},\Lambda'^{(t)})\in\calb^+_{Z^{(t)},Z'^{(t)}}}\rho_{\Lambda^{(t)}}\otimes\rho_{\Lambda'^{(t)}}\Biggr)^\sharp
=\frac{1}{2}\sum_{(\Sigma^{(t)},\Sigma'^{(t)})\in\cald_{Z^{(t)},Z'^{(t)}}} R_{\Sigma^{(t)}}\otimes R_{\Sigma'^{(t)}}.
\]
Therefore, we have
\begin{align*}
\tilde f_t\otimes\tilde f'_t\Biggl(\Biggl(\sum_{(\Lambda,\Lambda')\in\calb^+_{Z,Z'}}\rho_\Lambda\otimes\rho_{\Lambda'}\Biggr)^\sharp\Biggr)
&=\Biggl(C_t\sum_{(\Lambda^{(t)},\Lambda'^{(t)})\in\calb^+_{Z^{(t)},Z'^{(t)}}}
\rho_{\Lambda^{(t)}}\otimes\rho_{\Lambda'^{(t)}}\Biggr)^\sharp \\
&=\frac{1}{2} C_t\sum_{(\Sigma^{(t)},\Sigma'^{(t)})\in\cald_{Z^{(t)},Z'^{(t)}}} 
R_{\Sigma^{(t)}}\otimes R_{\Sigma'^{(t)}} \\
&=\tilde f_t\otimes\tilde f'_t\Biggl(\frac{1}{2}\sum_{(\Sigma,\Sigma')\in\cald_{Z,Z'}} R_\Sigma\otimes R_{\Sigma'}\Biggr)
\end{align*}
for some nonzero constant $C_t$.
Because now $\tilde f_t\otimes\tilde f'_t$ is a vector space isomorphism, we have
\[
\Biggl(\sum_{(\Lambda,\Lambda')\in\calb^+_{Z,Z'}}\rho_\Lambda\otimes\rho_{\Lambda'}\Biggr)^\sharp
=\frac{1}{2}\sum_{(\Sigma,\Sigma')\in\cald_{Z,Z'}} R_\Sigma\otimes R_{\Sigma'}.
\]
\end{proof}

\begin{remark}\label{0825}
If $\cald_{Z,Z'}$ is one-to-one, i.e., $\Psi_0=\Psi'_0=\emptyset$,
then $\binom{a_k}{b_l}$ or $\binom{c_{l'}}{d_{k'}}$ can only be pairs of doubles.
Then $\deg(Z)=\deg(Z^{(t)})$ and $\deg(Z')=\deg(Z'^{(t)})$ in the proof of Proposition~\ref{0813}.
Then we have either $\deg(Z')=\deg(Z)$ or $\deg(Z')=\deg(Z)+1$.
\end{remark}


\section{The Relation $\calb^-_{Z,Z'}$}\label{0803}
The purpose of this section (Proposition~\ref{0921}) is to show that the statement of
Theorem~\ref{0310} holds when $\epsilon=-$.
The strategy is parallel to that used in Section~\ref{0706} and Section~\ref{0802}.
In this section,
we assume that $\epsilon=-$ and $\cald_{Z,Z'}\neq\emptyset$.
Let $Z,Z'$ be special symbols of defect $1,0$,
and $\Psi_0,\Psi_0'$ the cores in $Z_\rmI,Z'_\rmI$ of $\cald_{Z,Z'}$ respectively.

\subsection{The set $\calb^-_{Z,Z'}$}

\begin{lemma}\label{0924}
Let $\Lambda\in\overline\cals_Z$ and $\Lambda'\in\overline\cals_{Z'}$.
Then
\[
(\Lambda,\Lambda')\in\overline\calb^+_{Z,Z'}\quad\text{if and only if}\quad
(\Lambda^\rmt,\Lambda'^\rmt)\in\overline\calb^-_{Z,Z'}.
\]
\end{lemma}
\begin{proof}
Because ${\rm def}(\Lambda^\rmt)=-{\rm def}(\Lambda)$ for any $\Lambda$,
it is clear that
${\rm def}(\Lambda')=-{\rm def}(\Lambda)+1$ if and only if
${\rm def}(\Lambda'^\rmt)=-{\rm def}(\Lambda^\rmt)-1$.
Write
\begin{align*}
\Lambda &=\binom{a_1,a_2,\ldots,a_{m_1}}{b_1,b_2,\ldots,b_{m_2}}, &\qquad
\Lambda' &=\binom{c_1,c_2,\ldots,c_{m'_1}}{d_1,d_2,\ldots,d_{m'_2}} \\
\Lambda^\rmt &=\binom{a'_1,a'_2,\ldots,a'_{m_2}}{b'_1,b'_2,\ldots,b'_{m_1}}, &\qquad
\Lambda' &=\binom{c'_1,c'_2,\ldots,c'_{m'_2}}{d'_1,d'_2,\ldots,d'_{m'_1}},
\end{align*}
i.e., $a'_i=b_i$, $b'_i=a_i$, $c'_i=d_i$, $d'_i=c_i$ for each $i$.
\begin{enumerate}
\item[(1)] Suppose that $m'=m$.
If $(\Lambda,\Lambda')\in\overline\calb^+_{Z,Z'}$, then
\begin{itemize}
\item $d'_i=c_i\geq b_i=a'_i$;

\item $a_i'=b_i>c_{i+1}=d'_{i+1}$;

\item $b'_i=a_i>d_i=c'_i$;

\item $c'_i=d_i\geq a_{i+1}=b'_{i+1}$
\end{itemize}
for each $i$ and hence $(\Lambda^\rmt,\Lambda'^\rmt)\in\overline\calb^-_{Z,Z'}$ by Lemma~\ref{0210}.
Conversely, by the same argument it is easy to see that
$(\Lambda^\rmt,\Lambda'^\rmt)\in\overline\calb^-_{Z,Z'}$ implies that
$(\Lambda,\Lambda')\in\overline\calb^+_{Z,Z'}$.

\item[(2)] The proof for $m'=m+1$ is similar.
\end{enumerate}
\end{proof}

Then we have the following analogue of Proposition~\ref{0810}:

\begin{proposition}
Suppose that $\cald_{Z,Z'}\neq\emptyset$.
\begin{enumerate}
\item[(i)] If $\Lambda'_0\in B^-_\Lambda$ for $\Lambda\in\overline\cals_Z$,
then $B^-_\Lambda=\{\,\Lambda'_0+\Lambda'\mid\Lambda'\in D_Z\,\}$.

\item[(ii)] If $\Lambda_0\in B^-_{\Lambda'}$ for $\Lambda'\in\overline\cals_{Z'}$,
then $B^-_{\Lambda'}=\{\,\Lambda_0+\Lambda\mid\Lambda\in D_{Z'}\,\}$.
\end{enumerate}
\end{proposition}
\begin{proof}
Suppose that $(\Lambda,\Lambda_0')\in\overline\calb^-_{Z,Z'}$, i.e., $\Lambda'_0\in B^-_\Lambda$.
By Lemma~\ref{0924}, we see that $\Lambda_0'^\rmt\in B^+_{\Lambda^\rmt}$.
Then $\Lambda_0'^\rmt+\Lambda'\in B^+_{\Lambda^\rmt}$ for any $\Lambda'\in D_Z$ by Proposition~\ref{0810}.
Then $(\Lambda_0'^\rmt+\Lambda')^\rmt\in B^-_\Lambda$ by Lemma~\ref{0924}, again.
Then by Lemma~\ref{0231}, we have $\Lambda_0'+\Lambda'=(\Lambda_0'^\rmt+\Lambda')^\rmt\in B^-_\Lambda$.
Thus we have shown that $\{\,\Lambda'_0+\Lambda'\mid\Lambda'\in D_Z\,\}\subset B^-_\Lambda$.

Conversely, let $\Lambda''\in B^-_\Lambda$.
Then $\Lambda''^\rmt\in B^+_{\Lambda^\rmt}$ by Lemma~\ref{0924} and hence $\Lambda''^\rmt=\Lambda_0'^\rmt+\Lambda'$ for some
$\Lambda'\in D_Z$ by Proposition~\ref{0810}.
Then $\Lambda''=\Lambda_0'+\Lambda'$ by Lemma~\ref{0231}.
Thus we have shown that $B^-_\Lambda\subset\{\,\Lambda'_0+\Lambda'\mid\Lambda'\in D_Z\,\}$.

The proof of (ii) is similar.
\end{proof}

For $\Psi\leq\Psi_0$ and $\Psi'\leq\Psi'_0$,
we define
\[
\calb_{Z,Z'}^{-,\Psi,\Psi'}=\calb^-_{Z,Z'}\cap(\cals_Z^\Psi\times\cals_{Z'}^{-,\Psi'})
\]
where $\cals_Z^\Psi$ and $\cals_{Z'}^{-,\Psi'}$ are defined in Subsection~\ref{0609}.
Then as in (\ref{0811}),
we have
\begin{equation}\label{1002}
\calb^-_{Z,Z'}=\{\,(\Lambda+\Lambda_M,\Lambda'+\Lambda_N)\mid(\Lambda,\Lambda')\in\calb^{-,\Psi,\Psi'}_{Z,Z'},\ M\leq\Psi,\ N\leq\Psi'\,\}
\end{equation}
If $\Psi=\Psi_0$ and $\Psi'=\Psi'_0$,
then $\calb_{Z,Z'}^{-,\Psi,\Psi'}$ is denoted by $\calb_{Z,Z'}^{-,\natural}$.
Moreover, we know that $\calb_{Z,Z'}^{-,\natural}$ is one-to-one if it is non-empty.

\subsection{The case for $Z,Z'$ regular and $\cald_{Z,Z'}$ one-to-one}\label{0922}

\begin{lemma}
Let $Z,Z'$ be special symbols of defects $1,0$ respectively.
\begin{enumerate}
\item[(i)] If $m'=m+1$, both $Z$ and $Z'$ are regular,
and $D_{Z'}=\{Z\}$, then every symbol in $\cals_Z$ occurs in the relation $\calb^-_{Z,Z'}$.

\item[(ii)] If $m'=m$, both $Z$ and $Z'$ are regular,
and $D_Z=\{Z'\}$, then every symbol in $\cals^-_{Z'}$ occurs in the relation $\calb^-_{Z,Z'}$.
\end{enumerate}
\end{lemma}
\begin{proof}
First suppose that $m'=m+1$.
By Corollary~\ref{0709}, every $\Lambda\in\overline\cals_Z$ occurs in the relation
$\overline\calb^+_{Z,Z'}$.
Note that taking transpose $\Lambda\mapsto\Lambda^\rmt$ is a bijection from $\overline\cals_Z$
onto itself, so, by Lemma~\ref{0924},
every $\Lambda\in\overline\cals_Z$ occurs in the relation $\overline\calb^-_{Z,Z'}$.
Hence every $\Lambda\in\cals_Z$ occurs in the relation $\calb^-_{Z,Z'}$.

The proof for the case that $m'=m$ is similar.
\end{proof}

Suppose that both $Z,Z'$ are regular, and let $\theta$ be given in Subsection~\ref{0707}.
\begin{enumerate}
\item[(1)] Suppose that $m'=m$.
Note that $a_1$ is not in the image of $\theta$ for this case.
Then $\theta$ induces a mapping $\theta^-\colon\overline\cals_{Z'}\rightarrow\overline\cals_Z$ by
$\Lambda_N\mapsto\Lambda_{\binom{a_1}{-}\cup\theta(N)}$.
Moreover, if $\Lambda\in\cals_{Z'}^-$, then $\theta^-(\Lambda')\in\cals_Z$.

\item[(2)] Suppose that $m'=m+1$.
Note that $c_1$ is not in the image of $\theta$.
Then $\theta$ induces a mapping $\theta^-\colon\overline\cals_Z\rightarrow\overline\cals_{Z'}$ by
$\Lambda_M\mapsto\Lambda_{\binom{c_1}{-}\cup\theta(M)}$.
Moreover, if $\Lambda\in\cals_Z$, then $\theta^-(\Lambda)\in\cals_{Z'}^-$.
\end{enumerate}

\begin{lemma}\label{0926}
Let $Z,Z'$ be special symbols of defect $1,0$ respectively.
Suppose that both $Z$ and $Z'$ are regular,
and $\cald_{Z,Z'}$ is one-to-one.
\begin{enumerate}
\item[(i)] If $m'=m+1$, then $(\Lambda,\Lambda')\in\calb^-_{Z,Z'}$ if and only if $\Lambda'=\theta^-(\Lambda)$.

\item[(ii)] If $m'=m$, then $(\Lambda,\Lambda')\in\calb^-_{Z,Z'}$ if and only if $\Lambda=\theta^-(\Lambda')$.
\end{enumerate}
\end{lemma}
\begin{proof}
First suppose that $m'=m+1$.
It is known that $(\Lambda_M,\Lambda')\in\overline\calb^+_{Z,Z'}$ if and only
if $\Lambda'=\theta^+(\Lambda)=\Lambda_{\theta(M)}$ by Lemma~\ref{0902}.
Then by Lemma~\ref{0924},
we see that $((\Lambda_M)^\rmt,\Lambda'')\in\overline\calb^-_{Z,Z'}$ if and only
if $\Lambda''=(\Lambda_{\theta(M)})^\rmt$.
Now $(\Lambda_M)^\rmt=\Lambda_{Z_\rmI\smallsetminus M}$ and
$(\Lambda_{\theta(M)})^\rmt=\Lambda_{Z'_\rmI\smallsetminus\theta(M)}=\Lambda_{\binom{c_1}{-}\cup\theta(Z_\rmI\smallsetminus M)}=\theta^-((\Lambda_M)^\rmt)$.

The proof for the case that $m'=m$ is similar.
\end{proof}

Now we prove that the statement in Theorem~\ref{0310} holds
when $\epsilon=-$, $Z$ and $Z'$ are regular, and the relation $\cald_{Z,Z'}$ is one-to-one.

\begin{proposition}\label{1003}
Let $(\bfG,\bfG')=(\Sp_{2n},\rmO^-_{2n'})$,
$Z,Z'$ special symbols of rank $n,n'$ and defect $1,0$ respectively.
Suppose that both $Z,Z'$ are regular, and $\cald_{Z,Z'}$ is one-to-one.
Then
\[
\frac{1}{2}\sum_{(\Sigma,\Sigma')\in\cald_{Z,Z'}} R^\bfG_\Sigma\otimes R^{\bfG'}_{\Sigma'}
=\sum_{(\Lambda,\Lambda')\in\calb^-_{Z,Z'}}\rho_\Lambda^\sharp\otimes\rho^\sharp_{\Lambda'}.
\]
\end{proposition}
\begin{proof}
The proof is similar to that of Proposition~\ref{0806}.
Suppose that $Z$ and $Z'$ are regular, and $\cald_{Z,Z'}$ is one-to-one.
Write $Z,Z'$ as in (\ref{0311}) and let $Z_{(m)},Z'_{(m')}$ be given in (\ref{0711}).
Let $h,h'$ be defined as in the proof of Proposition~\ref{0806}.
Then $h$ induces a bijection $\bar h\colon\cals_Z\rightarrow\cals_{Z_{(m)}}$ and $\cals_{Z,1}\rightarrow\cals_{Z_{(m)},1}$
such that
$\langle R_\Sigma,\rho_\Lambda\rangle=\langle R_{\bar h(\Sigma)},\rho_{\bar h(\Lambda)}\rangle$,
and $h'$ induces bijections
$\bar h'\colon\cals^-_{Z'}\rightarrow\cals^-_{Z'_{(m)}}$ and $\cals_{Z',0}\rightarrow\cals_{Z'_{(m)},0}$
such that
$\langle R_{\Sigma'},\rho_{\Lambda'}\rangle=\langle R_{\bar h'(\Sigma')},\rho_{\bar h'(\Lambda')}\rangle$.
We see that $(\Lambda,\Lambda')\in\calb^-_{Z,Z'}$ if and only if
$(\bar h(\Lambda),\bar h'(\Lambda'))\in\calb^-_{Z_{(m)},Z'_{(m')}}$ by Example~\ref{0225}, Example~\ref{0226}
and Lemma~\ref{0926}.
Hence the proposition is proved by the same argument in the proof of Proposition~\ref{0806}.
\end{proof}

\subsection{General case}
In this subsection, we assume that $\cald_{Z,Z'}$ is non-empty.
Suppose that $\cald_{Z,Z'}$ is not one-to-one or one of $Z,Z'$ is not regular.
So one of Cases (I), (II), (III) occurs.
Write $Z,Z'$ as in (\ref{0914}) and let $\binom{a_k}{b_l}$ and $\binom{c_{l'}}{d_{k'}}$
be given as in Subsection~\ref{0903}.
Let $Z^{(1)},Z'^{(1)}$ be defined as in Subsection~\ref{0903} according to Case (I),(II),(III) respectively.

\begin{enumerate}
\item[(1)] Suppose that we are in Case (I) of Subsection~\ref{0903}.
Now $f$ is a bijection from $Z_\rmI\smallsetminus\{a_k,b_l\}$ onto $Z^{(1)}_\rmI$,
and induces a bijection $\bar f\colon\overline\cals_Z\rightarrow\overline\cals_{Z^{(1)}}$
by $\bar f(\Lambda_M)=\Lambda_{f(M)}$ for $M\subset Z_\rmI$ if $\binom{a_k}{b_l}$ is a pair of doubles;
and induces a bijection $\bar f\colon\overline\cals_Z^\Psi\rightarrow\overline\cals_{Z^{(1)}}$
by $\bar f(\Lambda_M)=\Lambda_{f(M)}$ for $M\subset Z_\rmI\smallsetminus\{a_k,b_l\}$ and $\Psi=\binom{a_k}{b_l}$ if
$\binom{a_k}{b_l}$ is in the core.
Let $\rho^{(1)}_\Lambda$, $\calv_Z^{(1)}$ and $\calv_{Z^{(1)}}$ are defined as there.

Similarly, $f'$ is a bijection from $Z'_\rmI\smallsetminus\{c_{l'},d_{k'}\}$ onto $Z'^{(1)}_\rmI$,
and induces a bijection
$\bar f'\colon\overline\cals_{Z'}\rightarrow\overline\cals_{Z'^{(1)}}$
by $\bar f'(\Lambda_N)=\Lambda_{f'(N)}$ for $N\subset Z'_\rmI$ if $\binom{c_{l'}}{d_{k'}}$ is a pair of doubles;
and induces a bijection $\bar f'\colon\overline\cals_{Z'}^{\Psi'}\rightarrow\overline\cals_{Z'^{(1)}}$
by $\bar f'(\Lambda_N)=\Lambda_{f(N)}$ for $N\subset Z'_\rmI\smallsetminus\{c_{l'},d_{k'}\}$ and $\Psi'=\binom{c_{l'}}{d_{k'}}$ if
$\binom{c_{l'}}{d_{k'}}$ is in the core.
For $\Lambda'\in\cals_{Z'}^-$,
let $\rho_{\Lambda'}^{(1)}$ be defined as in Subsection~\ref{0826} and let
$\calv_{Z'}^{-,(1)}$ be the subspace of $\calv_{Z'}^-$ spanned by $\rho_{\Lambda'}^{(1)}$
for $\Lambda'\in\cals_{Z'}^-$.

\item[(2)] Suppose that we are in Case (II) of Subsection~\ref{0903}.
Now $f$ is given as in (1) and $f'$ is a bijection from $Z'_\rmI$ onto $Z'^{(1)}_\rmI$.
Let $\bar f$, $\calv_Z^{(1)}$, $\calv_{Z^{(1)}}$, $\bar f'$, $\calv^{-,(1)}_{Z'}$, $\calv^-_{Z'^{(1)}}$ be defined analogously.

\item[(3)] Suppose that we are in Case (III) of Subsection~\ref{0903}.
Now $f$ is a bijection from $Z_\rmI$ onto $Z^{(1)}_\rmI$, and
$f'$ is a bijection from $Z'_\rmI\smallsetminus\{c_{l'},d_{k'}\}$ onto $Z'^{(1)}_\rmI$.
Let $\bar f$, $\calv_Z^{(1)}$, $\calv_{Z^{(1)}}$, $\bar f',\calv^{-,(1)}_{Z'}$, $\calv^-_{Z'^{(1)}}$ be defined analogously.
\end{enumerate}
For all three cases above, we have vector space isomorphisms
\[
\tilde f\colon\calv_Z^{(1)}\longrightarrow\calv_{Z^{(1)}}\qquad \text{and}\qquad
\tilde f'\colon\calv_{Z'}^{-,(1)}\longrightarrow\calv_{Z'^{(1)}}^-.
\]
Then let $\calb_{Z,Z'}^{-,(1)}$ be defined similarly.

\begin{lemma}
For Cases (I),(II),(III) in Subsection~\ref{0903},
we have $(\Lambda,\Lambda')\in \calb^{-,(1)}_{Z,Z'}$ if and only if
$(\bar f(\Lambda),\bar f(\Lambda'))\in\calb^-_{Z^{(1)},Z'^{(1)}}$.
\end{lemma}
\begin{proof}
Suppose that $(\Lambda,\Lambda')\in\calb_{Z,Z'}^{-,(1)}\subset\calb^-_{Z,Z'}\subset\overline\calb^-_{Z,Z'}$ for
$\Lambda=\Lambda_M$ and $\Lambda'=\Lambda_N$
where $M\subseteq Z_\rmI\smallsetminus\{a_k,b_l\}$ for Cases (I) and (II);
$M\subseteq Z_\rmI$ for Case (III),
and $N\subseteq Z'_\rmI\smallsetminus\{c_{l'},d_{k'}\}$ for Cases (I) and (III);
$M\subseteq Z'_\rmI$ for Case (II).

By Lemma~\ref{0924}, we see that
$((\Lambda_M)^\rmt,(\Lambda_N)^\rmt)\in\overline\calb^+_{Z,Z'}$.
Now $(\Lambda_M)^\rmt=\Lambda_{Z_\rmI\smallsetminus M}$ and $(\Lambda_N)^\rmt=\Lambda_{Z'_\rmI\smallsetminus N}$.
Define
\[
M'=\begin{cases}
M\cup\binom{a_k}{b_l}, & \text{if $\binom{a_k}{b_l}$ is the core};\\
M, & \text{otherwise},
\end{cases}\qquad
N'=\begin{cases}
N\cup\binom{c_{l'}}{d_{k'}}, & \text{if $\binom{c_{l'}}{d_{k'}}$ is the core};\\
N, & \text{otherwise}.
\end{cases}
\]
By Proposition~\ref{0810}, we see that
$(\Lambda_{Z_\rmI\smallsetminus M'},\Lambda_{Z'_\rmI\smallsetminus N'})\in\overline\calb_{Z,Z'}^+$.
\begin{enumerate}
\item[(1)] For Case (I), (II), if $\binom{a_k}{b_l}$ is a pair of doubles,
then both $a_k,b_l$ are not in $Z_\rmI$,
hence $(Z_\rmI\smallsetminus M')=(Z_\rmI\smallsetminus M)\subset Z_\rmI\smallsetminus\{a_k,b_l\}$ for
any $M\subset Z_\rmI\smallsetminus\{a_k,b_l\}$;
if $\binom{a_k}{b_l}$ is in the core,
then $(Z_\rmI\smallsetminus M')\subset Z_\rmI\smallsetminus\{a_k,b_l\}$, again.

\item[(2)] For Case (I) or (III), if $\binom{c_{l'}}{d_{k'}}$ is a pair of doubles,
then both $c_{l'},d_{k'}$ are not in $Z'_\rmI$,
hence $(Z'_\rmI\smallsetminus N')=(Z'_\rmI\smallsetminus N)\subset Z'_\rmI\smallsetminus\{c_{l'},d_{k'}\}$ for
any $N\subset Z'_\rmI\smallsetminus\{c_{l'},d_{k'}\}$;
if $\binom{c_{l'}}{d_{k'}}$ is in the core,
then $(Z'_\rmI\smallsetminus N')\subset Z'_\rmI\smallsetminus\{c_{l'},d_{k'}\}$, again.
\end{enumerate}
So we have $(\Lambda_{Z_\rmI\smallsetminus M'},\Lambda_{Z'_\rmI\smallsetminus N'})\in\overline\calb_{Z,Z'}^{+,(1)}$.
Then by Lemma~\ref{0912} and Lemma~\ref{0909}, we have
$(\Lambda_{f(Z_\rmI\smallsetminus M')},\Lambda_{f'(Z'_\rmI\smallsetminus N')})\in\overline\calb_{Z^{(1)},Z'^{(1)}}^+$.
Then
\[
((\Lambda_{f(Z_\rmI\smallsetminus M')})^\rmt,(\Lambda_{f'(Z'_\rmI\smallsetminus N')})^\rmt)
\in\overline\calb_{Z^{(1)},Z'^{(1)}}^-
\]
by Lemma~\ref{0924}.
It is easy to see that $Z^{(1)}\smallsetminus f(Z_\rmI\smallsetminus M')=f(M)$ and
$Z'^{(1)}\smallsetminus f(Z'_\rmI\smallsetminus N')=f(N)$.
Then $(\Lambda_{f(Z_\rmI\smallsetminus M')})^\rmt=\Lambda_{Z^{(1)}\smallsetminus f(Z_\rmI\smallsetminus M')}=\Lambda_{f(M)}$
and $(\Lambda_{f'(Z'_\rmI\smallsetminus N')})^\rmt=\Lambda_{Z'^{(1)}\smallsetminus f'(Z'_\rmI\smallsetminus N')}=\Lambda_{f'(N)}$.
Therefore, we conclude that $(\bar f(\Lambda),\bar f(\Lambda'))\in\overline\calb^-_{Z^{(1)},Z'^{(1)}}$.
Because now $\bar f(\Lambda)\in\cals_{Z^{(1)}}$,
we have in fact $(\bar f(\Lambda),\bar f(\Lambda'))\in\calb^-_{Z^{(1)},Z'^{(1)}}$.

Conversely, we can also show that $(\bar f(\Lambda),\bar f(\Lambda'))\in\calb^-_{Z^{(1)},Z'^{(1)}}$ will imply
that $(\Lambda,\Lambda')\in \calb^{-,(1)}_{Z,Z'}$.
\end{proof}

\begin{lemma}
For Cases (I),(II),(III) in Subsection~\ref{0903},
we have
\begin{align*}
\sum_{(\Lambda,\Lambda')\in\calb^-_{Z,Z'}}\rho_\Lambda\otimes\rho_{\Lambda'}
&= C\sum_{(\Lambda,\Lambda')\in\calb^{(1),-}_{Z,Z'}}\rho^{(1)}_\Lambda\otimes\rho^{(1)}_{\Lambda'}\\
\sum_{(\Sigma,\Sigma')\in\cald_{Z,Z'}}R_\Sigma\otimes R_{\Sigma'}
&= C\sum_{(\Sigma,\Sigma')\in\cald^{(1)}_{Z,Z'}}R^{(1)}_\Sigma\otimes R^{(1)}_{\Sigma'}
\end{align*}
where
\[
C=\begin{cases}
1, & \text{if none of $\binom{a_k}{b_l},\binom{c_{l'}}{d_{k'}}$ is in the core};\\
\sqrt 2, & \text{if exactly one of $\binom{a_k}{b_l},\binom{c_{l'}}{d_{k'}}$ is in the core};\\
2, & \text{if both $\binom{a_k}{b_l},\binom{c_{l'}}{d_{k'}}$ are in the core}.
\end{cases}
\]
\end{lemma}
\begin{proof}
The proof is similar to those of Lemma~\ref{0823} and Lemma~\ref{0824}.
\end{proof}

\begin{corollary}
Suppose that $\epsilon=-$.
If the core of $\cald_{Z,Z'}$ in $Z'_\rmI$ is $Z'_\rmI$ itself,
then
\[
\sum_{(\Sigma,\Sigma')\in\cald_{Z,Z'}}R_\Sigma\otimes R_{\Sigma'}=0.
\]
\end{corollary}
\begin{proof}
Let $\Psi_0,\Psi_0'$ be the cores in $Z,Z'$ respectively.
Suppose that $\Psi_0'=Z'_\rmI$.
From (\ref{1002}) we know that
\[
\cald_{Z,Z'}=\{\,(\Sigma+\Lambda_M,\Sigma'+\Lambda_N)\mid(\Sigma,\Sigma')\in\cald^\natural_{Z,Z'},\ M\leq\Psi_0,\ N\leq\Psi'_0\,\}.
\]
Now our assumption $\Psi'_0=Z'_\rmI$ means that $\Sigma'=Z'$ and $\Sigma'+\Lambda_N=\Lambda_N$.
Because $\cald^\natural_{Z,Z'}$ is one-to-one,
we must have $\Sigma=Z$, and hence $\Sigma+\Lambda_M=\Lambda_M$.
Note that $(\Lambda_N)^\rmt=\Lambda_{Z_\rmI'\smallsetminus N}=\Lambda_{\Psi_0'\smallsetminus N}$,
so $\Lambda_N\in\cals_{Z',\Psi_0'}$ if and only if $(\Lambda_N)^\rmt\in\cals_{Z',\Psi_0'}$.
Let $X$ denote a set of representatives of cosets $\{\Lambda_N,(\Lambda_N)^\rmt\}$ in $\cals_{Z',\Psi'_0}$.
Then
\begin{align*}
\sum_{(\Sigma,\Sigma')\in\cald_{Z,Z'}}R_\Sigma\otimes R_{\Sigma'}
&= \biggl(\sum_{M\leq\Psi_0}R_{\Lambda_M}\biggr)\otimes\biggl(\sum_{\Lambda_N\in X}(R_{\Lambda_N}+R_{(\Lambda_N)^\rmt})\biggr)
\end{align*}
Because now $\epsilon=-$,
we have $R_{\Lambda_N}+R_{(\Lambda_N)^\rmt}=0$ for any $\Lambda_N\in X$.
\end{proof}

Now we prove that Theorem~\ref{0310} holds when $\epsilon=-$.

\begin{proposition}\label{0921}
Let $(\bfG,\bfG')=(\Sp_{2n},\rmO^-_{2n'})$,
$Z,Z'$ special symbols of rank $n,n'$ and defect $1,0$ respectively.
Then
\[
\frac{1}{2}\sum_{(\Sigma,\Sigma')\in\cald_{Z,Z'}} R^\bfG_\Sigma\otimes R^{\bfG'}_{\Sigma'}
=\sum_{(\Lambda,\Lambda')\in\calb^-_{Z,Z'}}\rho_\Lambda^\sharp\otimes\rho^\sharp_{\Lambda'}.
\]
\end{proposition}
\begin{proof}
The proof is similar to that of Proposition~\ref{0813}.
Suppose that $\epsilon=-$,
and let $Z,Z'$ be special symbols of defect $1,0$ respectively.
Suppose that $\cald_{Z,Z'}$ is not one-to-one or one of $Z,Z'$ is not regular.
Let $\Psi_0$ and $\Psi'_0$ be the cores of the relation $\cald_{Z,Z'}$ in
$Z_\rmI$ and $Z'_\rmI$ respectively.
As in the proof of Proposition~\ref{0813},
we can find special symbols $Z^{(t)},Z'^{(t)}$ of defects $1,0$ respectively and bijections
\[
\bar f_t\colon\cals_Z^{\Psi_0}\longrightarrow\cals_{Z^{(t)}}\qquad\text{and}\qquad
\bar f'_t\colon\cals_{Z'}^{-,\Psi'_0}\longrightarrow\cals^-_{Z'^{(t)}}
\]
such that
\begin{itemize}
\item both $Z^{(t)},Z'^{(t)}$ are regular,

\item $\deg(Z\smallsetminus\Psi_0)=\deg(Z^{(t)})$ and $\deg(Z'\smallsetminus\Psi'_0)=\deg(Z'^{(t)})$,

\item $\cald_{Z^{(t)},Z'^{(t)}}$ is one-to-one,
in particular, $\cald_{Z^{(t)},Z'^{(t)}}\neq\emptyset$,

\item $(\Sigma,\Sigma')\in\cald_{Z,Z'}^\natural$ if and only if
$(\bar f_t(\Sigma),\bar f'_t(\Sigma'))\in\cald_{Z^{(t)},Z'^{(t)}}$,

\item $(\Lambda,\Lambda')\in\calb^{-,\natural}_{Z,Z'}$ if and only if
$(\bar f_t(\Lambda),\bar f'_t(\Lambda'))\in\calb^-_{Z^{(t)},Z'^{(t)}}$.
\end{itemize}
Moreover, we have
\begin{align*}
\tilde f_t\otimes\tilde f'_t\Biggl(\sum_{(\Lambda,\Lambda')\in\calb^-_{Z,Z'}}\rho_\Lambda\otimes\rho_{\Lambda'}\Biggr)
&=\sum_{(\Lambda^{(t)},\Lambda'^{(t)})\in\calb^-_{Z^{(t)},Z'^{(t)}}}\rho_{\Lambda^{(t)}}\otimes\rho_{\Lambda'^{(t)}} \\
\tilde f_t\otimes\tilde f'_t\Biggl(\sum_{(\Sigma,\Sigma')\in\cald_{Z,Z'}}R_\Sigma\otimes R_{\Sigma'}\Biggr)
&=\sum_{(\Sigma^{(t)},\Sigma'^{(t)})\in\cald_{Z^{(t)},Z'^{(t)}}} R_{\Sigma^{(t)}}\otimes R_{\Sigma'^{(t)}}.
\end{align*}

Now the special symbols $Z^{(t)},Z'^{(t)}$ are regular and the relation $\cald_{Z^{(t)},Z'^{(t)}}$ is one-to-one.
So by Lemma~\ref{0213}, we see that
either $\deg(Z'^{(t)})=\deg(Z^{(t)})$ or $\deg(Z'^{(t)})=\deg(Z^{(t)})+1$.
From the results in Lemma~\ref{0926} and bijections $f_t,f'_t$,
we have the following two situations:
\begin{enumerate}
\item[(1)] Suppose that $\deg(Z'^{(t)})=\deg(Z^{(t)})$.
Then we have a mapping $\theta^-\colon\cals^{-,\Psi'_0}_{Z'}\rightarrow\cals_Z^{\Psi_0}$
such that $(\Lambda,\Lambda')\in\calb^{-,\natural}_{Z,Z'}$ if and only if
$\Lambda=\theta^-(\Lambda')$.

\item[(2)] Suppose that $\deg(Z'^{(t)})=\deg(Z^{(t)})+1$.
Then we have a mapping $\theta^-\colon\cals_Z^{\Psi_0}\rightarrow\cals^{-,\Psi'_0}_{Z'}$
such that $(\Lambda,\Lambda')\in\calb^{-,\natural}_{Z,Z'}$ if and only if
$\Lambda'=\theta^-(\Lambda)$.
\end{enumerate}
Hence by Proposition~\ref{1003}, we have
\[
\Biggl(\sum_{(\Lambda^{(t)},\Lambda'^{(t)})\in\calb^-_{Z^{(t)},Z'^{(t)}}}\rho_{\Lambda^{(t)}}\otimes\rho_{\Lambda'^{(t)}}\Biggr)^\sharp
=\frac{1}{2}\sum_{(\Sigma^{(t)},\Sigma'^{(t)})\in\cald_{Z^{(t)},Z'^{(t)}}} R_{\Sigma^{(t)}}\otimes R_{\Sigma'^{(t)}}.
\]
Therefore as in the proof of Proposition~\ref{0813}, we have
\[
\tilde f_t\otimes\tilde f'_t\Biggl(\Biggl(\sum_{(\Lambda,\Lambda')\in\calb^-_{Z,Z'}}\rho_\Lambda\otimes\rho_{\Lambda'}\Biggr)^\sharp\Biggr)
=\tilde f_t\otimes\tilde f'_t\Biggl(\frac{1}{2}\sum_{(\Sigma,\Sigma')\in\cald_{Z,Z'}} R_\Sigma\otimes R_{\Sigma'}\Biggr).
\]
Because now $\tilde f_t\otimes\tilde f'_t$ is a vector space isomorphism, we conclude that
\[
\Biggl(\sum_{(\Lambda,\Lambda')\in\calb^-_{Z,Z'}}\rho_\Lambda\otimes\rho_{\Lambda'}\Biggr)^\sharp
=\frac{1}{2}\sum_{(\Sigma,\Sigma')\in\cald_{Z,Z'}} R_\Sigma\otimes R_{\Sigma'}.
\]
\end{proof}

\bibliography{refer}
\bibliographystyle{amsalpha}

\end{document}